\documentclass[letterpaper,11pt]{article}
\usepackage[latin1]{inputenc}
\usepackage[T1]{fontenc}
\usepackage[spanish]{babel} 
\usepackage[pdftex]{hyperref}
\usepackage{amsmath,amsfonts,amssymb}
\usepackage{amsthm}
\usepackage{multicol} 
\usepackage{graphicx}
\usepackage{fancyhdr}

\newtheorem{lema}{Lema}[section]
\newtheorem{teo}{Teorema}[section]
\newtheorem{coro}{Corolario}[section]
\newtheorem{conje}{Conjetura}[section]
\newtheorem{defi}{Definici\'on}[section]
\newtheorem{eje}{Ejemplo}[section]
\author{\sc
		Dimas N. T. Tejada\\
		Facultad de Ciencias Naturales y Matemática \\
		Universidad de El Salvador, El Salvador \\
		\textit{E-mail:} {\tt dimas.tejada@ues.edu.sv } }
\title{\bf An\'alisis Global de la Ecuaci\'on de Lyness}
\date{ Enero de 2023}

\hypersetup{colorlinks=true,linkcolor=black}

\begin{document}

\maketitle

\begin{abstract}
En este trabajo de divulgación se describe el comportamiento global de la ecuación de Lyness. Se obtienen resultados para diferentes valores del parámetro $\alpha$.

\vspace*{.2cm}
\centerline{\bf Abstract}
In this work of divulgation,  the global behavior of the Lyness equation is described. Results are obtained for different values of the parameter $\alpha$.
\end{abstract}


\begin{center}\bf Introducci\'on\end{center} 

Originalmente, en 1942 R. C. Lyness present\'o a modo de divertimiento matem\'atico en \emph{The Mathematical Gazetta} las siguientes ecuaciones que son \emph{3-ciclos, 4-ciclos, 6-ciclos y 5-ciclos} respectivamente (ecuaciones en las que todas las soluciones son constantes o peri\'odicas del mismo per\'iodo)
\begin{align*}
x_{n+1}x_{n-1}+a&=-x_n(x_{n-1}+x_{n+1})\\
x_{n+1}x_{n-1}+a&=x_n(x_{n-1}+x_n+x_{n+1})\\
x_{n+1}x_{n-1}+a&=x_n(x_{n-1}+2x_n+x_{n+1})\\
x_{n+1}x_{n-1}-a^2&=ax_n
\end{align*}

\noindent y se preguntaba si pod\'ia existir ecuaciones que fueran \emph{7-ciclos}.\\

En otro art\'iculo en \emph{The Mathematical Gazetta} de 1945, Lyness analiza con m\'as profundidad el ejemplo de \emph{5-ciclo}  anterior y llega a la conclusi\'on que el  \emph{5-ciclo} m\'as sencillo a considerar es $$x_{n+1}x_{n-1}-1=x_n$$ que desde entonces se conoce como \emph{la ecuaci\'on de Lyness}. Si tomamos el punto $(x_1,x_2)=(x,y)$ como punto inicial de una soluci\'on de $\{x_n\}_{n\geq1}$, tenemos
\begin{align*}
x_3&=\frac{1+y}{x},\\ 
x_4&=\frac{1+\frac{1+y}{x}}{y}=\frac{1+x+y}{xy},\\
x_5&=\frac{1+\frac{1+x+y}{xy}}{\frac{1+y}{x}}=\frac{1+x+y+xy}{y(1+y)}=\frac{(1+x)(1+y)}{y(1+y)}=\frac{1+x}{y},\\
x_6&=\frac{1+\frac{1+x}{y}}{\frac{1+x+y}{xy}}=\frac{(1+x+y)x}{1+x+y}=x,\\
x_7&=\frac{1+x}{\frac{1+x}{y}}=y
\end{align*}
de donde se deduce que es un \emph{5-ciclo}.\\
 
En un art\'iculo en \emph{Acta Arithmetica} de 1971, H. S. Coxeter puso en actualidad la idea de los \emph{frieze patterns} que ya hab\'ian aparecido en el siglo XVI y que Gauss hab\'ia incorporado en su \emph{pentagramma mirificum}. Un \emph{frieze pattern} es una disposici\'on de n\'umeros enteros formando rombos de tal forma que si los n\'umeros $a,b,c,d$ ocupan sus extremos, se verifica la propiedad $ad-bc=1$. La conclusi\'on sorprendente es que cada uno de tales patrones es peri\'odico. La prueba realizada por Coxeter de este hecho es esencialmente la comprobaci\'on de que ciertas ecuaciones en diferencias sean globalmente peri\'odicas. La \emph{ecuaci\'on de Lyness} es uno de estos ejemplos.\\

Se denomina \emph{ecuaci\'on generalizada de Lyness}, a la ecuaci\'on 

$$x_{n+1}x_{n-1}-x_n=\alpha$$

\noindent que generalmente la escribimos

\begin{equation}\label{sub:ecintro}
x_{n+1}=\frac{\alpha+x_n}{x_{n-1}} 
\end{equation}
 
\noindent donde $\alpha \geq 0$ y $x_{n-1}\neq0$. Si las condiciones iniciales $x_1,x_2$ en esta ecuaci\'on son n\'umeros positivos, entonces tambi\'en lo ser\'a la sucesi\'on $\{x_n\}_{n\geq1}$, y este es el caso que nos ocupara.\\

Si $\alpha\neq1$, es f\'acil ver que la ecuaci\'on anterior no es un \emph{5-ciclo}. Para estudiar sus propiedades extendemos el problema a un difeomorfismo en el primer cuadrante de $\mathbb{R}^2$ dado por $$f(x,y)=(y,\frac{\alpha+y}{x}).$$

Este procedimiento se conoce como un desplegamiento de la ecuaci\'on en diferencias. Las \'orbitas generadas por la iteraci\'on de $f$ producen  soluciones de la ecuaci\'on anterior proyectando la \'orbita en el eje  de las abscisas. El desplegamiento del problema a una dimensi\'on m\'as, nos permite analizar el comportamiento de la ecuaci\'on utilizando las t\'ecnicas de los \emph{Sistemas Din\'amicos Discretos} y en algunos teoremas la \emph{Geometr\'ia Algebraica}.\\

Para analizar el comportamiento del difeomorfismo anterior, se comprueba que $$V(x,y)=(\frac1x+1)(\frac1y+1)(\alpha+x+y)$$ es un invariante de $f$, es decir,  $V(f(x,y))=V(x,y)$ para $\forall(x,y) \in \mathbb{R}^2_+$, de aqu\'i que las \'orbitas de $f$ est\'an sobre las curvas de nivel de $V$, dichas curvas de nivel son obtenidas por la  ecuaci\'on $$V(x,y)=(x+1)(y+1)(x+y+\alpha)=vxy$$ donde $v \in \mathbb{R}$. Se comprueba f\'acilmente tambi\'en, que \'estas  son cerradas, acotadas y \emph{llenan} el primer cuadrante, es decir, por cada punto de $\mathbb{R}^2_+$ pasa una \'unica curva de la familia; por tanto, al proyectarlas en el eje de las abscisas se obtiene que todas las soluciones de la ecuaci\'on (~\ref{sub:ecintro}) son acotadas. \\

$V$ es una funci\'on  que tiene un \'unico punto cr\'itico $F_\alpha=(\omega_\alpha,\omega_\alpha)$ donde $\omega_\alpha=\frac{1+\sqrt{1+4\alpha}}{2}$ y es no degenerado, tiene un m\'inimo en dicho punto cr\'itico $v_\alpha=\min V=V(F_\alpha)=\frac{(\omega_\alpha+1)^3}{\omega_\alpha}$, y para cada $v>v_\alpha$ la curva de nivel asociada $C^v_\alpha$ es difeomorfa a la circunferencia unidad. \\

En 1996, C. Zeeman prob\'o que la funci\'on $f|_{C^v_\alpha}$ es conjugada a una rotaci\'on de la circunferencia unidad y que el n\'umero de rotaci\'on $\rho^v_\alpha$ asociado a esta funci\'on depende derivablemente de $v$. La consecuencia de esto es que las \'orbitas de $f$ son todas peri\'odicas o densas dependiendo si el n\'umero de rotaci\'on es racional o irracional.\\

Las propiedades de $\rho^v_\alpha$ han sido estudiadas. En 1997, F. Beukers y R. Cushman probaron que tal funci\'on es anal\'itica real en el intervalo $[v_\alpha,\infty)$ y que es estrictamente creciente si $0<\alpha<1$ o decreciente si $1<\alpha<\infty$.\\

Finalmente, si $\alpha=1$ se obtiene que todas las \'orbitas son de per\'iodo 5 (resultado que obtuvimos directamente) por lo que en todas ellas el n\'umero de rotaci\'on es $\frac15$. Pero cuando $\alpha\neq1$, se verifica que $$\lim_{v\to\infty} \rho^v_\alpha=\frac15$$

\noindent a trav\'es de funciones mon\'otonas de la variable $v$.



\section{Notas hist\'oricas: \emph{frieze patterns}}
\subsection{Introducci\'on}
La idea de un  \emph{frieze patterns}(patrones o modelos con los que se construyen los frisos cl\'asicos de los templos) es m\'as f\'acil de explicar por medio de un ejemplo, como el que sigue que es un \emph{pattern} de orden 7:\\
\begin{table}[h!]
\begin{tabular}{p{.1cm}p{.1cm}p{.1cm}p{.1cm}p{.1cm}p{.1cm}p{.1cm}p{.1cm}p{.1cm}p{.1cm}p{.1cm}p{.1cm}p{.1cm}p{.1cm}p{.1cm}p{.1cm}p{.1cm}p{.1cm}p{.1cm}p{.cm}p{.1cm}p{.1cm}p{.1cm}p{.1cm}}
0      & &0& &0& &0& &0& &0& &0& &0& &0& &0& &0& &0&       \\ 
       &1& &1& &1& &1& &1& &1& &1& &1& &1& &1& &1& &\ldots \\ 
\ldots & &1& &2& &2& &3& &1& &2& &4& &1& &2& &2& &3&       \\
       & & &1& &3& &5& &2& &1& &7& &3& &1& &3& &5& &\ldots  \\
\ldots & &2& &1& &7& &3& &1& &3& &5& &2& &1& &7& &3&        \\ 
       & & &1& &2& &4& &1& &2& &2& &3& &1& &2& &4& &\ldots  \\
\ldots & &1& &1& &1& &1& &1& &1& &1& &1& &1& &1& &1&       \\ 
       & & &0& &0& &0& &0& &0& &0& &0& &0& &0& &0& &\ldots     
\end{tabular}

\caption[Ejemplo de \emph{frieze patterns}.]{Ejemplo de \emph{frieze patterns}. \label{sub:pattern1} }
\end{table}\\

Se colocan ceros y unos en los bordes y los dem\'as n\'umeros se colocan fomando cuadrados
$$
\begin{tabular}{ccc}
 &b& \\
a& &d\\
 &c&
\end{tabular}
$$
tales que $ad-bc=1$, donde adem\'as, todos los n\'umeros (excepto los bordes de ceros) tienen que ser positivos. A esta disposici\'on se le denomina un patr\'on (el que se muestra en el ejemplo es de orden 7). La sorprendente conclusi\'on que se obtiene es que todos estos patrones  son peri\'odicos. M\'as precisamente, son sim\'etricos por desplazamientos: el producto de una traslaci\'on horizontal y una reflexi\'on horizontal. Es evidente y obligado dar condiciones necesarias y suficientes para que un \emph{frieze pattern} consista de enteros.

\subsection{\emph{Frieze patterns} de orden 5.}
\begin{figure}[h!]
\begin{center}
\includegraphics[scale=0.4]{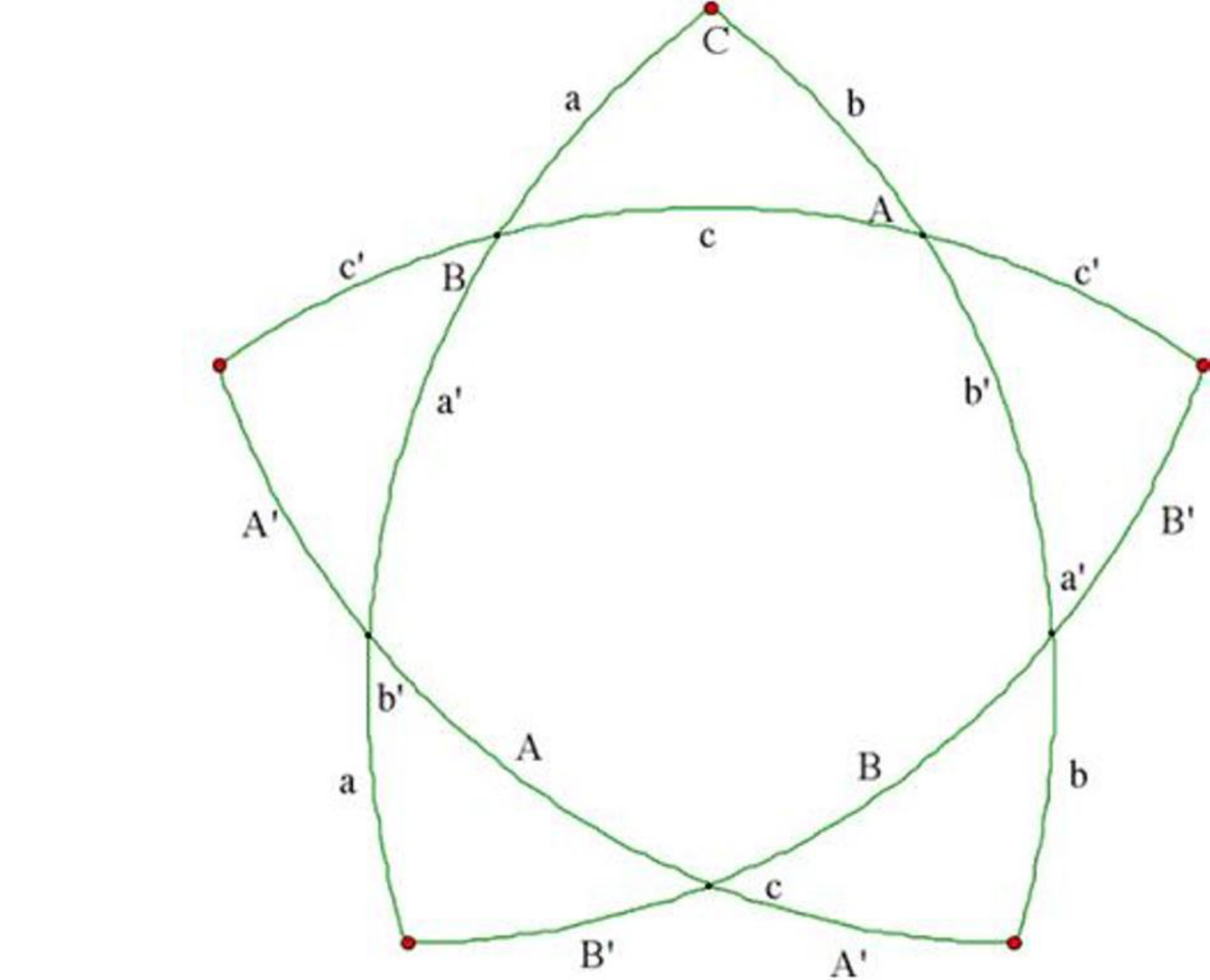}
\end{center}
\caption[\emph{Pentagramma Mirificum}]{\emph{Pentagramma Mirificum}.\label{sub:pentagrama}}
\end{figure}

La historia comienza en 1602, cuando Nathaniel Torporley (1564-1632) comenz\'o a investigar los cinco lados y \'angulos representados por $a, A, b, B, c$ de un tri\'angulo rect\'angulo esf\'erico (con el \'angulo recto en $C$ como en la figura ~\ref{sub:pentagrama}). Torporley anticip\'o en unos doce a\~nos las famosas reglas de Napier que Gauss incluy\'o en su \emph{pentagramma mirificum} (Las reglas de Napier permiten derivar una gran parte de las f\'ormulas de la trigonometr\'ia esf\'erica con muy poco esfuerzo). Por otra parte, Gauss utiliz\'o la identidad n\'umeros reales
$$(1+\gamma)(1+\beta-\delta\epsilon)-(1+\beta)(1+\gamma-\epsilon\alpha)=\epsilon\left[(1+\alpha-\gamma\delta)-(1+\delta-\alpha\beta)\right]$$
para probar que tres cualesquiera de las cuatro relaciones num\'ericas
\begin{equation}\label{sub:cincorelaciones}1+\alpha=\gamma\delta,1+\beta=\delta\epsilon,1+\gamma=\epsilon\alpha,1+\delta=\alpha\beta \end{equation}
implican la restante y que adem\'as $1+\alpha=\beta\gamma$. Esta observaci\'on establece la periodicidad del \emph{frieze pattern}
$$
\begin{tabular}{ccccccccccccccccccc}
0& &0       &        &0       &          &0       &        &0       &        &0         &        &0       &        &0      &          &        &        &    \\ 
 &1&        &1       &        &1         &        &1       &        &1       &          &1       &        &1       &       &1         &        &        &     \\ 
 & &$\alpha$&        &$\beta$ &          &$\gamma$&        &$\delta$&        &$\epsilon$&        &$\alpha$&        &$\beta$&          &$\gamma$&        &     \\
 & &        &$\delta$&        &$\epsilon$&        &$\alpha$&        &$\beta$ &          &$\gamma$&        &$\delta$&       &$\epsilon$&        &$\alpha$& \ldots\\
 & &        &        &1       &          &1       &        &1       &        &1         &        &1       &        &1      &          &1       &        &     \\ 
 & &        &        &        &0         &        &0       &        &0       &          &0       &        &0       &       &0         &        &        &    
\end{tabular}
$$

El \emph{pentagramma mirificum} (ver figura~\ref{sub:pentagrama}) es un pentagrama esf\'erico formado por cinco arcos de circunferencia m\'aximo. El ``n\'ucleo'' del pentagrama es un pent\'agono cuyos v\'ertices son (obviamente) polos de los cinco arcos; decimos entonces que es un pent\'agono auto-polar(esto quiere decir que si un v\'ertice del pent\'agono interior se pone en el polo norte de la esfera, el lado opuesto queda sobre el ecuador). La figura completa puede ser derivada del tri\'angulo rect\'angulo $ABC$ (dibujado en la parte superior) extendiendo los lados y dibujando los arcos polares de circunferencia m\'axima de los v\'ertices $A$ y $B$. Usando $x'$ para indicar el complemento
$$x' = 1/2 \pi - x,$$
vemos f\'acilmente que los lados del pent\'agono auto-polar son
$$A, B, b', c, a',$$
mientras que los restantes arcos y \'angulos son como indica la figura~\ref{sub:pentagrama}. Claramente, cualquier f\'ormula que conecte la cinco letras, sigue siendo v\'alida cuando se permutan en forma c\'iclica. \'Esta es la observaci\'on de Napier (salvo que \'el us\'o el ciclo alternado $a', B, c, A, b'$). Gauss us\'o los siguientes n\'umeros

 $$\alpha = \tan^2 A, \beta = \tan^2 B, \gamma = \tan^2 b', \delta = \tan^2 c, \epsilon = \tan^2 a'$$
 
\noindent y us\'o una de las relaciones cl\'asicas (como por ejemplo $\cos c = \cot A \cot B$) para deducir

$$\sec^2 A = \gamma\delta,\sec^2 B = \delta\epsilon,\sec^2 b' = \epsilon\alpha,\sec^2 c = \alpha\beta,\sec^2 a' = \beta\gamma.$$
Estas igualdades implican que $1 +\alpha = \gamma\delta$, y as\'i con las dem\'as igualdades.\\

La obvia igualdad $$\frac{1}{\alpha\beta}=\frac{\epsilon}{\beta\gamma}\frac{\gamma}{\epsilon\alpha}$$
da una relaci\'on esf\'erica an\'aloga al teorema de Pit\'agoras:
$$\cos c=\sin a'\sin b'=\cos a \cos b.$$

Lobachevsky mostr\'o que para cualquier tri\'angulo recto $ABC$ en el plano hiperb\'olico correspondiente al tri\'angulo esf\'erico con hipotenusa $\Pi(a)$ y cateto $B$ y $\Pi(c)$ opuesto los \'angulos $A'$ y $\Pi(b)$,  $\cot \Pi(a)=\sinh a$, y por tanto $\tan \Pi(a)=\frac{1}{\sinh a}$, con esto obtenemos las siguientes interpretaciones hiperb\'olicas para los \emph{frieze patterns} de Gauss:

$$\cot^2 A=\alpha,\frac{1}{\sinh^2b}=\beta,\sinh^2c=\gamma,\frac{1}{\sinh^2a}=\delta,\cot^2b=\epsilon,$$
$$\csc^2A=\gamma\delta,\coth^2b=\delta\epsilon,\cosh^2c=\epsilon\alpha,\coth^2a=\alpha\beta,\frac{1}{\sinh^2 B}=\beta\gamma.$$
Cambiando la notaci\'on 
$$u_1=\alpha,u_2=\delta,u_3=\beta,u_4=\epsilon,u_5=\gamma$$
tenemos el \emph{frieze pattern}
$$
\begin{tabular}{p{.2cm}p{.2cm}p{.2cm}p{.2cm}p{.2cm}p{.2cm}p{.2cm}p{.2cm}p{.2cm}p{.2cm}p{.2cm}p{.2cm}p{.2cm}p{.2cm}p{.2cm}p{.2cm}p{.2cm}p{.2cm}p{.2cm}p{.2cm}p{.2cm}p{.2cm}}
0& &0    &     &0       &          &0       &        &0       &        &0         &        &0       &        &0      &          &        &        &    \\ 
 &1&     &1    &        &1         &        &1       &        &1       &          &1       &        &1       &       &1         &        &        &     \\ 
 & &$u_1$&     &$u_3$   &          &$u_5$   &        &$u_2$   &        &$u_4$     &        &$u_1$   &        &$u_3$  &          &$u_5$   &        &     \\
 & &     &$u_2$&        &$u_4$     &        &$u_1$   &        &$u_3$   &          &$u_5$   &        &$u_2$   &       &$u_4$     &        &$u_1$   & \ldots\\
 & &     &     &1       &          &1       &        &1       &        &1         &        &1       &        &1      &          &1       &        &     \\ 
 & &     &     &        &0         &        &0       &        &0       &          &0       &        &0       &       &0         &        &        &    
\end{tabular}
$$
y las relaciones $$u_1u_3=1+u_2, \quad u_2u_4=1+u_3, \quad u_3u_5=1+u_4, \ldots,$$
y por tanto esto implica que $u_{r+5}=u_r$ para todo $r$ (siempre que $u_r\neq0$).  Estos $5-ciclos$ fueron conocidos en  matem\'aticas por un largo tiempo. La siguiente sencilla prueba fue dada por \emph{Israel Halperin}. Dados
$$u_1u_3=1+u_2, \quad u_2u_4=1+u_3, \quad u_3u_5=1+u_4,\quad u_4u_6=1+u_5$$
y $u_3u_4\neq0$, tenemos que
$$u_1u_3u_4=(1+u_2)u_4=u_4+1+u_3=u_3u_5+u_3=u_3(1+u_5)=u_3u_4u_6;$$
por lo tanto $u_1=u_6$. En las siguientes secciones estudiaremos la periodicidad para un \emph{frieze pattern} de orden $n$ con una nueva notaci\'on.

\subsection{La simbolog\'ia de los dos d\'igitos $(r,s)$.}
Para investigar la generalidad del \emph{frieze pattern}, buscamos una simbolog\'ia conveniente para representar estos elementos por s\'imbolos $(r,s)$ como sigue:\\
\begin{table}[h!]
\begin{tabular}{p{1cm}p{1.3cm}p{1.3cm}p{1.3cm}p{1.2cm}p{1.2cm}p{.6cm}p{1.2cm}}
(0,0)  &        &(1,1)   &        &  (2,2)   &       &(3,3)    & \ldots  \\
       &(0,1)   &        &(1,2)   &          &(2,3)  &         &(3,4)   \\
(-1,1) &        &(0,2)   &        &  (1,3)   &       &(2,4)    &          \\ 
\ldots &\ldots  &\ldots  &\ldots  & \ldots   &\ldots &\ldots   &\ldots  \\
       &(-1,n-3)&        &(0,n-2) &          &(1,n-1)&         &(2,n)      \\
       &        &(-1,n-2)&        & (0,n-1)  &       &(1,n)    & \ldots \\ 
       &        &        &(-1,n-1)&          &(0,n)  &         &(1,n+1)     
\end{tabular}      
\caption[Una nueva notaci\'on para los\emph{frieze patterns}.]{Una nueva notaci\'on para los \emph{frieze patterns}. \label{pattern3} }
\end{table}\\

Los bordes de ceros y unos son obtenidos por las ecuaciones siguientes:
\begin{equation}\label{sub:ceros}
(r,r)=0, \quad \quad (r,r+n)=0
\end{equation}
\begin{equation}\label{sub:unos}
(r,r+1)=1, \quad \quad (r+1,r+n)=1
\end{equation}
Tambi\'en tenemos que $(r,s)>0$ para $r<s<r+n$, porque trabajamos solo con n\'umeros positivos ($r$ y $s$ pueden ser negativos, pero el valor $(r,s)$ es un elemento del \emph{frieze pattern} que es siempre positivo); la relaci\'on
\begin{equation}\label{sub:regla}
(r-1,s)(r,s+1)-(r,s)(r-1,s+1)=1,
\end{equation}
que implica
\begin{equation}\label{sub:unosnegativos}
(r,r-1)=-1, \quad \quad (r-1,r+n)=-1
\end{equation}
Claramente, todo el \emph{ frieze pattern} es determinado por los elementos en una diagonal, esto es $$(0,0),(0,1),(0,2),\ldots,(0,n-1),(0,n)$$
(recordando que comienza con ceros y unos y termina con unos y ceros). Por ejemplo, $(1,2)=1$, y tomando a $r=1$ en la ecuaci\'on (~\ref{sub:regla}), obtenemos que
$$(1,3)=\frac{(1,2)(0,3)+1}{(0,2)}, \quad (1,4)=\frac{(1,3)(0,4)+1}{(0,3)} $$ etc\'etera. En consecuencia, adoptaremos la notaci\'on auxiliar 
\begin{equation}\label{sub:notacionfg}
f_s=(-1,s), \quad g_s=(0,s),
\end{equation}
de modo que 
\begin{align*}
f_{-2}&=g_{-1}=f_n=g_{n+1}=-1,\\
f_{-1}&=g_0=f_{n-1}=g_n=0,\\
f_0&=g_1=f_{n-2}=g_{n-1}=1,\\
g_2&=\frac{g_1f_2+1}{f_1},\\
g_3&=\frac{g_2f_3+1}{f_2},\\
\vdots&\\
g_{n-2}&=\frac{g_{n-3}f_{n-2}+1}{f_{n-3}}
\end{align*}
Podemos verificar que, para valores v\'alidos de $r$ y $s$, es
\begin{equation}\label{sub:valoresvalidos}
(r,s)=f_rf_s-f_sg_r.
\end{equation}
De hecho, esta definici\'on de $(r,s)$ (para todos los enteros $r$ y $s$) implica
\begin{align}\label{sub:relacion33}
(r,s)(t,u)+(r,t)(u,s)+(r,u)(s,t)=&0,\\
(s,r)=&-(r,s),
\end{align}

y $(r-1,s)(r,s+1)-(r,s)(r-1,s+1)-(r-1,r)(s,s+1)=0$. \\

La \'ultima relaci\'on est\'a en concordancia con la ecuaci\'on (~\ref{sub:regla}), siendo $(r-1,r)=(s,s+1)=1$.

\subsection{La periodicidad de los \emph{frieze patterns.}}
Usando las ecuaciones (~\ref{sub:ceros}),(~\ref{sub:unos}),(~\ref{sub:unosnegativos}),(~\ref{sub:relacion33}), tenemos que 
\begin{equation}\label{sub:relacion34}
(r,s)+(r,s+n)=-(r,s)(s-1,s+n)+(r,s+n)(s-1,s)=-(r,s-1)(s,s+n)=0.
\end{equation}
En particular, $$f_{s+n}=-f_s, \quad g_{s+n}=-g_s.$$ Estos valores de $f_s$ para el \emph{pattern} en la tabla~\ref{sub:pattern1} son
$$0,1,1,1,1,2,1,0,-1,-1,-1,-1,-2,-1,0,1,\ldots$$
De la ecuaci\'on (~\ref{sub:relacion33}) y (~\ref{sub:relacion34}) tenemos que 
\begin{equation}\label{sub:relacion35}
(r,s)=(s,r+n)=(r+n,s+n)
\end{equation}
Esta ecuaci\'on muestra que los \emph{frieze patterns} son sim\'etricos por desplazamientos, y por lo tanto tambi\'en por traslaci\'on horizontal. As\'i, si definimos $$a_r=(r-1,r+1),$$ tenemos que 
\begin{equation}\label{sub:relacion40}a_{r+kn}=a_r\end{equation} 
para todo entero $k$. Adem\'as, los $n$ n\'umeros $a_0,a_1,\ldots,a_{n-1}$ no son independientes. Luego veremos que, en lugar de derivar todo el \emph{pattern} de $f_1,f_2,\ldots,f_{n-s}$, igualmente podemos derivarlo de $a_0,a_1,\ldots,a_{n-4}$.\\

De la igualdad
$$(r,s)=(s-2,s)(r,s-1)-(r,s-2)(s-1,s)=a_{s-1}(r,s-1)-(r,s-2),$$ 
sigue por inducci\'on que 
\begin{align}
(r,s)=&\begin{vmatrix} 
a_{r+1} & 1       &  0    & \dots &0   & 0 \\ 
1       & a_{r+2} &1      & \dots &0   & 0 \\
0       &1        &a_{r+3}&\dots  &0   & 0 \\
      \hdotsfor[3.5]{4}          \\ 
0       &0        &0      &\dots &1    & a_{s-1} 
\end{vmatrix} 
\end{align}
En particular, la ecuaci\'on (~\ref{sub:unos}) (con $r=-2$) muestra que $a_{n-3}$ puede derivarse de $a_0,a_1,\ldots,a_{n-4}$ por la ecuaci\'on lineal
\begin{align}
\begin{vmatrix} 
a_0 & 1   &  0& \dots &0   & 0 \\ 
1   & a_1 &1  & \dots &0   & 0 \\
0   &1    &a_2&\dots  &0   & 0 \\
      \hdotsfor[3.5]{4}          \\ 
0   &0    &0  &\dots  &1   & a_{n-3} 
\end{vmatrix}&=1, \tag*{$n>3$.} 
\end{align}
En consecuencia, la periodicidad en la ecuaci\'on (~\ref{sub:relacion40}) puede ser expresada como sigue: 
\emph{Sean $n-3$  n\'umeros positivos $a_0,a_1,\ldots,a_{n-4}$ dados, y sea la sucesi\'on $\{a_r\}$ construida siguiendo las reglas tal que
\begin{align}
\begin{vmatrix} 
a_r & 1       &  0    & \dots &0   & 0 \\ 
1   & a_{r+1} &1      & \dots &0   & 0 \\
0   &1        &a_{r+2}&\dots  &0   & 0 \\
      \hdotsfor[3.5]{4}          \\ 
0   &0        &0      &\dots  &1   & a_{r+n-3} 
\end{vmatrix}&=1, \tag*{$r=0,1,\ldots$.} 
\end{align}
Entonces la sucesi\'on es peri\'odica: $a_r=a_{r+n}$.}

\subsection{\emph{Frieze patterns} de enteros.}
Si todos lo n\'umeros $(r,s)$ son enteros (como el ejemplo en la tabla~\ref{sub:pattern1}), la regla muestra que cada $(r,s)$ es coprimo respecto a sus cuatro m\'as cercanos $$(r\pm 1,s) \quad y \quad (r,s\pm1).$$
En particular, para cualquiera seis de los n\'umeros ordenados como sigue:
$$
\begin{tabular}{cccc}
 &b& & \\
a& &d& \\
 &c& &f\\
 & &e&
\end{tabular}
$$
$c$ y $d$ son coprimos. Luego tenemos que, $ad-bc=1=cf-de$, as\'i
\begin{equation}\label{sub:relacion36}
(a+e)d=(b+f)c,
\end{equation}
entonces $c$ divide a $a+e$, y $d$ divide a $b+f$. En otras palabras, \emph{tres enteros consecutivos  sobre la diagonal son tal que el del centro divide la suma de los otros dos.} Inversamente, si en un \emph{frieze pattern}, $a,b,c,d,e$ son enteros tal que $c$ divide a $a+e$, entonces la ecuaci\'on (~\ref{sub:relacion36}) muestra que $b+f$ es un entero, y por lo tanto $f$ es un entero. Volviendo a la notaci\'on de pares, inferimos que, si $f_0,f_1,\ldots,f_{n-2}$ es una sucesi\'on de enteros iniciando y finalizando con $1$, y si 
\begin{equation} \label{sub:ecuacion10}
f_s\,\, divide\,\, a \,\,f_{s-1}+f_{s+1}
\end{equation}
entonces todos los $g_s$ son enteros y por lo tanto, por  ecuaci\'on (~\ref{sub:valoresvalidos}), todo los $(r,s)$ son enteros. En otras palabras, un \emph{frieze pattern} consiste solo de enteros si y s\'olo si la sucesi\'on generadora $f_0,f_1,\ldots,f_{n-2}$(iniciando y terminando con 1) consite solo de enteros satisfaciendo la ecuaci\'on (~\ref{sub:ecuacion10}).\\ 

Por ejemplo, $f_0,f_1,\ldots,f_{n-2}$ puede ser $1,1,1,\dots,1$, o cualquier sucesi\'on de n\'umeros que cumpla la ecuaci\'on (~\ref{sub:ecuacion10}). Adem\'as, varias sucesiones adecuadas pueden yuxtaponerse para hacer una nueva; por ejemplo, $$1,2,5,3,1 \quad y \quad 1,2,3,4,1$$ 
pueden combinarse para obtener $$1,2,5,3,1,1,2,3,4,1,$$
y por supuesto cada 1 puede ser representado por un cadena de cualquier cantidad de 1. Despu\'es eligiendo $f_{s-1}$ y $f_s$, podemos tomar $mf_s-f_{s-1}$ para cualquier entero $m>\frac{f_{s-1}}{f_s}$. La \'unica dificultad consiste en asegurarse que, para un \emph{pattern} de orden $n$, $f_{n-2}=1$.\\

El siguiente ejemplo (basado en la sucesi\'on de arriba) ilustra el hecho de que un \emph{pattern} de enteros no necesariamente incluye una diagonal consistiendo completamente de los n\'umeros 1 y 2:
$$
\begin{tabular}{p{.1cm}p{.1cm}p{.1cm}p{.1cm}p{.1cm}p{.1cm}p{.1cm}p{.1cm}p{.1cm}p{.1cm}p{.1cm}p{.1cm}p{.1cm}p{.1cm}p{.1cm}p{.1cm}p{.1cm}p{.1cm}p{.1cm}p{.1cm}p{.cm}p{.1cm}p{.1cm}p{.1cm}p{.1cm}}
 &0& &0& &0& &0&  &0 &  &0 & &0& &0& &0&  &0 &  &0 & &0&    \\ 
1& &1& &1& &1& &1 &  &1 &  &1& &1& &1& &1 &  &1 &  &1& &\dots  \\ 
 &2& &3& &1& &2&  &5 &  &2 & &2& &1& &4&  &2 &  &2 & &3&       \\
 & &5& &2& &1& &9 &  &9 &  &3& &1& &3& &7 &  &3 &  &5& &\ldots  \\
 & & &3& &1& &4&  &16&  &13& &1& &2& &5&  &10&  &7 & &3&        \\ 
 & & & &1& &3& &7 &  &23&  &4& &1& &3& &7 &  &23&  &4& &\ldots  \\
 & & &1& &2& &5&  &10&  &7 & &3& &1& &4&  &16&  &13& &1&       \\ 
 & & & &1& &3& &7 &  &3 &  &5& &2& &1& &9 &  &9 & &3 & &\ldots \\
 & & & & &1& &4&  &2 &  &2 & &3& &1& &2&  &5 &  &2&  &2&       \\  
 & & & & & &1& &1 &  &1 &  &1& &1& &1& &1 &  &1 &  &1& &\dots  \\ 
 & & & & &0& &0&  &0 &  &0 & &0& &0& &0&  &0 &  &0 & &0.&     
\end{tabular}
$$


\section{Resultados sobre Sistemas Din\'amicos Discretos}
\sectionmark{Res. sobre Sistemas Din\'amicos}

El objetivo fundamental de la teor\'ia de sistemas din\'amicos discretos es entender la evoluci\'on de un proceso iterativo dado por $(X,f)$ donde $X$ es un espacio topol\'ogico y $f$ una funci\'on continua, es decir, comprender el comportamiento de la sucesi\'on de puntos del espacio de fases $\{f^n(x)\}^{n=\infty}_{n=0}$ cuando $n$ es grande. Informalmente hablando, nos preguntaremos ¿Hacia d\'onde evoluciona esta sucesi\'on de puntos y qu\'e hace cuando llega a alg\'un lugar?. Esta forma de analizar las funciones tienen una gran aplicaci\'on en resoluci\'on de problemas en matem\'atica.

\subsection{\'Orbitas y periodicidad}

\begin{defi}
La \'orbita hacia delante de $x$ por medio de una funci\'on $f$ es el conjunto de puntos $x,f(x),f^2(x),\ldots,f^n(x),\dots$ y es denotada por $O^+_f(x)$. Si $f$ es un homeomorfismo, se puede definir la \'orbita completa\footnote{En todo este trabajo, usaremos como ``\'orbita'' solo a la ``\'orbita hacia delante'' salvo que se indique lo contrario.} de de $x$, denotada por $O_f(x)$, que es el conjunto de puntos $f^n(x)$ para $n\in \mathbb{Z}$, y la \'orbita hacia atr\'as de $x$, denotada por $O^-_f(x)$, es el conjunto de puntos  $x,f^{-1}(x),f^{-2}(x),\ldots,f^{-n}(x),\dots$. Donde $f^n=\underbrace{f\circ f\circ f\circ \ldots \circ f}_{n-veces}$ con $n \in \mathbb{Z}$.
\end{defi}

\begin{defi}
El punto $x$ es un punto fijo para $f$ si $f(x)=x$. El punto $x$ es peri\'odico de per\'iodo $n$ si $f^n(x)=x$. El menor n\'umero positivo $n$  para el que $f^n(x)=x$ es llamado el per\'iodo principal de $x$.
\end{defi}

Se denota por $Per_n(f)$ al conjunto de puntos peri\'odicos (no necesariamente principales) de per\'iodo $n$. El conjunto de todas las iteraciones de un punto peri\'odico forma una \'orbita peri\'odica.

\begin{defi}
Un punto $x$ es eventualmente peri\'odico de per\'iodo $n$ si $x$ no es peri\'odico, pero existe un $m>0$ tal que $f^{n+i}(x)=f^i(x)$ para todo $i\geq m$. Esto es, $f^i(x)$ es peri\'odica para todo $i \geq m$.
\end{defi}

\begin{defi}
Sea $p$ un punto peri\'odico de per\'iodo $n$. Un punto $x$ es asint\'otico hacia $p$ si $\lim_{x \to \infty} f^{in}(x)=p$. El conjunto estable de $p$, denotado por $W^s(p)$, es el conjunto de puntos  asit\'oticos a $p$.
\end{defi}


\subsection{Transformaciones  de la circunferencia unidad en s\'i misma.}
\subsectionmark{Transformaciones  de la circunferencia.}

En esta secci\'on presentamos propiedades  importantes del sistema din\'amico $(S^1,f)$ donde con $S^1$ representamos la circunferencia unidad\footnote{En lo que sigue, hablaremos indistintamente circunferencia unidad o simplemente circunferencia.} y $f$ es un endomorfismo continuo de $S^1$ en s\'i misma. Para ello comenzamos introduciendo los denominados \emph{levantamientos}\footnote{Traducci\'on de la palabra en ingles ``liftings''} para transformaciones de la circunferencia unidad en general. Luego, trabajamos algunas propiedades del n\'umero de rotaci\'on asociado a las transformaciones de la circunferencia, pero en este caso centraremos nuestra atenci\'on en difeomorfismos de $S^1$ que preservan la orientaci\'on, es decir, difeomorfismos $f:S^1\rightarrow S^1$ que preservan el orden sobre la circunferencia. Lo anterior lo haremos por simplicidad y porque nuestro inter\'es son las rotaciones de la circunferencia que son evidentemente transformaciones que preservan la orientaci\'on (ver \cite{devaney}).
\subsection{\emph{Levantamientos} y sus propiedades.}
Dada la circunferencia unidad $S^1=\{z\in \mathbb{C}:|z|=1\}$, definamos como proyecci\'on natural $\pi_1:R\rightarrow S^1$ $$\pi_1(x)=e^{2\pi ix}.$$
Esta funci\'on es continua y sobreyectiva. Utilizaremos la proyecci\'on natural $\pi_1$  para estudiar propiedades topol\'ogicas de transformaciones de la circunferencia. Esta transformaci\'on nos permite definir el concepto de un \emph{levantamiento} de transformaciones continuas de $S^1$ en s\'i misma. Para esto necesitamos alguna notaci\'on, definiciones y los siguientes dos \emph{lemas} auxiliares.\\

Sea $A$ y $B$ dos subconjuntos de $\mathbb{R}$. Para $x\in \mathbb{R}$ definimos $A+x=\{a+x:a\in A\}$, y $A+B=\{a+b:a\in A,b\in B\}$. Es inmediato el siguiente resultado,
\begin{lema} \label{sub:lemadeconjuntos}
Sea $v \in \mathbb{R}$. Entonces las siguientes afirmaciones son ciertas:
\begin{itemize}
\item[a)] La funci\'on $\bar{\pi}_1:(v,v+1)\rightarrow S^1\backslash \{\pi_1(v)\}$, obtenida por la restricci\'on de $\pi_1$ al intervalo $(v,v+1)$, es un homeomorfismo.
\item[b)] Sea $B$ cualquier subconjunto de $S^1\backslash \{\pi_1(v)\}$, y $A=[v,v+1]\cap \pi^{-1}_1(B)$. Entonces $\pi^{-1}_1(B)$ es la uni\'on de los conjuntos $A+n,(n\in \mathbb{Z})$, cada uno de estos conjuntos es abierto en $\pi^{-1}_1(B)$ y $\pi_1$ induce un homeomorfismo de cada $A+n$ sobre $B$.
\end{itemize}
\end{lema}

Sea $J$ un subintevalo de $[0,1]$ y sea $g:J \rightarrow S^1$ una funci\'on continua. Consideremos el problema: ¿Existir\'a  una funci\'on continua $F:J \rightarrow \mathbb{R}$ tal que $g=\pi_1\circ F$? Si esta $F$ la encontramos, la llamaremos un \emph{levantamiento} en el intervalo $J$  de $g$ y decimos que $g$ puede ser \emph{levantada}.\\

\begin{lema}\label{sub:lemadelifting}
Cualquier funci\'on continua $g:[0,1] \rightarrow S^1$ tiene un \emph{levantamiento} $F:[0,1] \rightarrow \mathbb{R}$ en el intervalo $[0,1]$, que es \'unico excepto por traslaciones por un entero. As\'i, si $y \in \mathbb{R}$  con $\pi_1(y)=g(0)$, existe un \'unico \emph{levantamiento} $F$ en el intervalo $[0,1]$ con $F(0)=y$.
\end{lema}

\begin{proof}
Para cada $x \in [0,1]$, la continuidad de $g$ en $x$ implica que podemos encontrar un $\epsilon_x>0$ tal que $g([x-\epsilon_x,x+\epsilon_x]\cap [0,1])$ sea un subcunjunto conexo de $S^1$. Por ser $[0,1]$ compacto, podemos encontrar una colecci\'on finita de  conjuntos  de la forma $(x-\epsilon_x,x+\epsilon_x)\cap [0,1]$ que cubran a $[0,1]$. Denotemos el conjunto de todos los puntos finales de cada intervalo por $0=x_0<x_1< \ldots <x_n =1$. Entonces cada $[x_{i-1},x_i]$ es contenido en alg\'un intervalo $[x-\epsilon_x,x+\epsilon_x]\cap [0,1]$, as\'i cada imagen por $g$ es un conjunto conexo $S_i$ de $S^1$.\\

Ahora probaremos, por inducci\'on sobre $i$, que hay un \'unico \emph{levantamiento} en el intervalo $[0,1]$ de $g|_{[0,x_i]}$ para la funci\'on $F:[0,x_i] \rightarrow \mathbb{R}$ con $F(0)=y$. Para $i=0$ es obvio. Supongamos que sea cierto para $i-1$. Sea $u_i$ un punto de $S^1$ tal que $u_i \in S_i$ (que es  un conjunto conexo), $v_i \in \mathbb{R}$ tal que $\pi_1(v_i)=u_i$ y $A_i=\pi^{-1}_1(S_i)\cap[v_i,v_{i+1}]$. Entonces por el Lema~\ref{sub:lemadeconjuntos} sabemos  que $\pi^{-1}_1(S_i)$ es una  uni\'on disjunta de conjuntos abiertos $A_i+n$ en $\pi^{-1}_1(S_i)$, cada funci\'on es homeomorfica sobre $S_i$. Sea $n_i$ un entero tal que $F(x_{i-1}) \in A_i+n_i$, y sea $\bar{\pi}_i: A_i+n_i \rightarrow S_i$ el homeomorfismo inducido por $\pi_1$. Entonces podemos definir  $F$ sobre $[0,x_i]$ como $\bar{\pi}^{-1}_iog$. Es f\'acil ver  que esto combinado con  $F$ sobre $[0,x_{i-1}]$ obtenemos una funci\'on continua sobre $[0,x_i]$. Claramente, $F:[0,x_i] \rightarrow \mathbb{R}$ es un \emph{levantamiento} en el intervalo $[0,1]$ de $g|_{[0,x_i]}$.\\

Inversamente, sea $\bar{F}$ un \emph{levantamiento} continuo cualquiera de $g$ en el intervalo $[0,1]$ de $[x_{i-1},x_i]$ en la uni\'on  de los conjuntos $A_i+n$ con $n\in \mathbb{Z}$. Sabiendo que  $[x_{i-1},x_i]$  es conexo, entonces todos las funciones $\bar{F}$ llegan a un \'unico $A_i+n$ ($A_i+n$ son conjuntos abiertos disjuntos en $\pi^{-1}_1(S_i)$), y como $x_{i-1}$ es transformado por medio de $\bar{F}$ dentro  $A_i+n_i$, as\'i todo elemento del intervalo $[x_{i-1},x_i]$ debe ser transformado por $\bar{F}$ al conjunto $A_i+n_i$. De esta manera cualquier \emph{levantamiento} $\bar{F}$  en el intervalo $[x_{i-1},x_i]$ queda determinado de manera \'unica. Con esto completamos el \'ultimo paso de inducci\'on. 
\end{proof}

Sea $f$ una transformaci\'on de la circunferencia. Por el Lema~\ref{sub:lemadelifting}, $f\circ (\pi_1|_{[0,1]})$ tiene como  \emph{levantamiento} una funci\'on $F:[0,1] \rightarrow \mathbb{R}$, como consecuencia que el siguiente diagrama conmuta:
\begin{center}
\begin{tabular}{ccc} 
$[0,1]$ & $\stackrel{F}{\rightarrow}$ & $\mathbb{R}$ \\ 
 $\pi_1|_{[0,1]}\,\,\,\downarrow \,\,\,\,\,\,\,\,\,\,\,\,\,\,\,\,\,\,\,\,$ &  & $\,\,\,\,\,\,\,\,\downarrow \,\,\,\pi_1$ \\ 
$S^1$ & $\stackrel{f}{\rightarrow}$ & $S^1.$ 
\end{tabular}
\end{center}
Como $\pi_1(0)=\pi_1(1)=1$, tenemos que $\pi_1(F(1))=f(\pi_1(1))=f(\pi_1(0))=\pi_1(F(0))$. Por lo tanto, $F(1)-F(0)$ es un entero. Este entero es llamado el \emph{grado} de $f$, es denotado por  $deg\,f$\footnote{Abreviatura por la palabra en ingles ``degree''.}. Cualquier otro \emph{levantamiento} $\bar{F}$ de $f$ en el intervalo $[0,1]$  es obtenido de $F$ por una traslaci\'on por un entero $n$ por el Lema~\ref{sub:lemadelifting}. Por lo tanto $$\bar{F}(1)-\bar{F}(0)=F(1)+n-(F(0)+n)=F(1)-F(0),$$ y se tiene que el grado es independiente de la elecci\'on del \emph{levantamiento} en el intervalo $[0,1]$.\\

Hablando sin mucha precisi\'on, el grado de la transformaci\'on $f$ es el m\'inimo n\'umero de veces que la imagen de $S^1$ por $f$  cubre completamente a $S^1$, en sentido contrario de las agujas del reloj si es positivo y en el sentido de las agujas del reloj si es negativo.\\

Dado un \emph{levantamiento} $F:[0,1] \rightarrow \mathbb{R}$ de $f\circ (\pi_1|_{[0,1]})$ para una transformaci\'on de la circunferencia $f$, podemos extender \'este a una funci\'on  $F':\mathbb{R} \rightarrow \mathbb{R}$ por 
\begin{equation} \label{sub:ecuacionlifting} 
F'(x)=F(x-E[x])+E[x].deg\,f,
\end{equation} 
Donde $E[x]$ es la funci\'on parte entera. De ahora en adelante, llamaremos a $F'$ un \emph{levantamiento} de $f$ porque el diagrama
\begin{center}
\begin{tabular}{ccc} 
$\mathbb{R}$ & $\stackrel{F'}{\rightarrow}$ & $\mathbb{R}$ \\ 
 $\pi_1\,\downarrow \,\,\,\,\,\,$ &  & $\,\,\,\,\,\,\,\,\downarrow \,\pi_1$ \\ 
$S^1$ & $\stackrel{f}{\rightarrow}$ & $S^1$ 
\end{tabular}
\end{center}
 conmuta, es decir, 
 \begin{align*}
 (\pi_1\circ F')(x)&=\pi_1(F(x-E[x])+E[x].deg\,f)\\
                 &=\pi_1(F(x-E[x])\\
                 &=f\circ \pi_1|_{[0,1]}(x-E[x])\\
                 &=(f\circ \pi_1)(x).
 \end{align*}

Para un entero $d$ dado, denotaremos por $\Gamma_d$ al conjunto de todos los \emph{levantamientos} de transformaciones de la circunferencia de grado $d$.

\begin{lema}\label{sub:lemaconjuntos}
Sea $d$ un entero, sea $Y\subset \mathbb{R}$ un conjunto tal que $Y+Z=Y$y sea $F:Y \rightarrow Y$ una funci\'on continua tal que $F(x+1)=F(x)+d$ para todo $x\in Y$. Entonces $F(x+k)=F(x)+kd$ para todo $x \in Y$ y $k$ un entero.
\end{lema} 
\begin{proof}
Si $k=0$, es obvio.  Si $k>0$ entonces tenemos $$F(x+k)=F(x+k-1)+d=F(x+k-2)+2d=\ldots=F(x)+kd.$$
Si $k<0$ entonces $F(x)=F(x+k)-k)=F(x+k)-kd$, por ser $-k>0$. As\'i, se obtiene que $F(x+k)=F(x)+kd$ tambi\'en en este caso.
\end{proof}

\begin{teo}\label{sub:teolifting1}
Una funci\'on continua $F:\mathbb{R} \rightarrow \mathbb{R}$  pertenece a $\Gamma_d$ si y s\'olo si $F(x+1)=F(x)+d$ para todo $x \in \mathbb{R}$.
\end{teo}
\begin{proof}
Supongamos primero que $F \in \Gamma_d$. Entonces $E[x+1]=E[x]+1$ y de la ecuaci\'on (~\ref{sub:ecuacionlifting}) obtenemos que $F(x+1)=F(x)+d$.\\

Supongamos ahora que  $F(x+1)=F(x)+d$ para todo $x \in \mathbb{R}$. Si $y \in S^1$ y $\pi_1(z)=\pi_1(z')=y$ entonces $z-z' \in \mathbb{Z}$ y por el Lema~\ref{sub:lemaconjuntos}, $F(z')=F(z)+n$ para $n=(z'-z)d \in \mathbb{Z}$. Por lo tanto, $\pi_1(F(z'))=\pi_1(F(z))$, es decir, $\pi_1(F(z))$ es independiente de la elecci\'on de $z \in \pi^{-1}_1(y)$. Por eso, podemos definir una funci\'on $f:S^1 \rightarrow S^1$ por $f(y)=\pi_1(F(z)),z \in \pi^{-1}_1(y)$. Si $y \in S^1$, $z \in \pi^{-1}_1(y)$ y $0<\epsilon<\frac12$ entonces $f|_{\pi_1([z-\epsilon,z+\epsilon])}=\pi_1\circ F\circ (\pi_1|_{([z-\epsilon,z+\epsilon])})^{-1}$ es continua por ser composici\'on de tres funciones continuas. De esta forma, $f$  es continua en cada $y \in S^1$, por lo que es continua sobre todo $S^1$. Por definici\'on, $F$  es un \emph{levantamiento} de $f$ y $deg\,f=F(1)-F(0)=d$. Por lo tanto, $F \in \Gamma_d$.
\end{proof}

En el siguiente teorema describimos algunas de las  propiedades b\'asicas  que cumplen los \emph{levantamientos}. Por $F+k$ vamos a denotar la funci\'on definida por $(F+k)(x)=F(x)+k$.

\begin{teo}\label{sub:propiedadeslifting}
Sea $f$  una transformaci\'on de la circunferencia de grado $d$ y sea $F$ un \emph{levantamiento} de $f$. Entonces las siguientes propiedades son ciertas:
\begin{itemize}
\item[a)] Si $F'$ es otro \emph{levantamiento} de $f$, entonces $F=F'+k$ para alg\'un entero $k$.
\item[b)] Si $k \in \mathbb{Z}$  entonces $F+k$ es tambi\'en \emph{levantamiento} de $f$.
\item[c)] $F^n(x+k)=F^n(x)+kd^n$ para todo $x \in \mathbb{R}$, $k \in \mathbb{Z}$ y $n\geq 0$.
\item[d)] $(F+k)^n(x)=F^n(x)+k(1+d+d^2+\ldots+d^{n-1})$ para todo $x \in \mathbb{R}$, $k \in \mathbb{Z}$, y $n\geq 0$.
\end{itemize}
\end{teo}
\begin{proof}
Sea $F'$ otro \emph{levantamiento} de $f$. Como  $\pi_1(F(x))=f(\pi_1(x))=\pi_1(F'(x))$,  obtenemos que $(F-F')(x) \in \mathbb{Z}$ para todo $x \in \mathbb{R}$. La funci\'on $F-F'$ es continua, y por lo tanto hay alg\'un $k \in \mathbb{Z}$ tal que $F-F'=k$. As\'i, a) est\'a probado. Probar b) simplemente es $f(\pi_1(x))=\pi_1(F(x))=\pi_1(F(x)+k)=\pi_1((F+k)(x))$.\\

Por el Teorema~\ref{sub:teolifting1}  y el Lema~\ref{sub:lemaconjuntos} tenemos que $F(z+k)=F(z)+kd$ para todo $z \in \mathbb{R}$ y $k \in \mathbb{Z}$. Por lo tanto
\begin{align*}
F^n(x+k)&=F^{n-1}(F(x+k))=F^{n-1}(F(x)+kd)\\
        &=F^{n-2}(F(F(x)+kd))=F^{n-2}(F^2(x)+kd^2)\\
        &\vdots \\
        &=F^n(x)+kd^n,
\end{align*}
y c) est\'a probado. Usando b) y c) obtenemos que 
\begin{align*}
(F+k)^n(x)&=(F+k)^{n-1}(F(x)+k)=(F+k)^{n-1}(F(x))+kd^{n-1}\\
         &=(F+k)^{n-2}(F^2(x)+k)+kd^{n-1}\\
         &=(F+k)^{n-2}(F^2(x))+kd^{n-2}+kd^{n-1}\\
         &\vdots\\
         &=F^n(x)+k(1+d+d^2+\ldots +d^{n-1}).
\end{align*}
As\'i d) est\'a probado.
\end{proof}

\subsection{N\'umero de rotaci\'on.}
El  invariante m\'as importante asociado a las trasnformaciones de la circunferencia es el que denominamos \emph{n\'umero de rotaci\'on}. Este n\'umero, entre $0$ y $1$, nos facilita la informaci\'on esencial del comportamiento de la transformaci\'on.\\

Si consideramos un difeomorfismo $f$ de la circunferencia que preserva el orden y le asociamos un \emph{levantamiento} $F$, entonces por la definici\'on de \emph{levantamiento} tenemos que $e^{2\pi i F(x)}=f(e^{2\pi ix})$ y por tanto la funci\'on $F$ tiene que ser creciente. Por ser $f$ un difeomorfismo sabemos que si $d=1$ y $k=1$ por el Teorema~\ref{sub:propiedadeslifting}  tenemos que $F(x+1)=F(x)+1$, luego $F(x+1)-(x+1)=F(x)+1-(x+1)=F(x)-x$. As\'i, la funci\'on $F-id$ (donde $id(x)=x$ es la funci\'on identidad) es una funci\'on peri\'odica con per\'iodo 1.\\

Adem\'as
\begin{align*}
\pi_1\circ F^n&=(\pi_1\circ F)\circ F^{n-1}=(f\circ \pi_1)\circ F^{n-1} \\
         &=f\circ (\pi_1 \circ F)\circ F^{n-2}=f^2\circ (\pi_1\circ F^{n_2})\\
         &\vdots \\
         &=f^n\circ \pi_1
\end{align*}

\noindent por lo tanto $F^n$ es un \emph{levantamiento} de $f^n$ y se tiene tambi\'en que $F^n-id$ es peri\'odica de per\'iodo 1.\\

Por otra parte, si $|x-y|<1$, entonces por ser $F$ creciente, por la periodicidad de $F^n-id$ y suponiendo sin perdidad de generalidad que $y<x$ ($x=y+\delta$ donde $0<\delta<1$) tenemos, 
\begin{equation} \label{sub:desigualdad1}| F^n(x)-F^n(y)|=F^n(x)-F^n(y)= F^n(y+\delta)-F^n(y)\leq F^n(y+1)-F^{n}(y)=1\end{equation} donde la primera igualdad y la desigualdad se verifican debido a la monoton\'ia de  $F$ y por tanto de $F^n$. \\

\begin{defi}
Sea $f:S^1\rightarrow  S^1$ un difeomorfismo que preserva la orientaci\'on y elijamos cualquier \emph{levantamiento} $F$ de $f$. Definamos 
$$\bar{\rho}(F)=\lim_{n\to \infty} \frac{F^n(x)}{n}.$$ 
\end{defi}

$\bar{\rho}(F)$ depende  del \emph{levantamiento} que tomemos, si es $F_0$ otro \emph{levantamiento} de $f$, entonces por el Teorema~\ref{sub:propiedadeslifting} es $F_0=F+k$  para un cierto entero $k$ y tambi\'en $(F+k)^n(x)=F^n(x)+nk$, as\'i
$$\bar{\rho}(F_0)=\lim_{n\to \infty} \frac{F_0^n(x)}{n}=\lim_{n\to \infty} \frac{(F+k)^n(x)}{n}=\lim_{n\to \infty} \frac{F^n(x)+nk}{n}=\bar{\rho}(F)+k$$
y por tanto difieren en un valor entero. Como veremos a continuaci\'on lo que nos  interesa es la parte fraccionaria de $\bar{\rho}(F)$ y por tanto no es importante el \emph{levantamiento} que tomemos.

\begin{defi}
El \emph{n\'umero de rotaci\'on} del difeomorfismo $f$, $\rho(f)$, es la parte fraccionaria de $\bar{\rho}(F)$ para cualquier \emph{levantamiento} $F$ de $f$. Esto es, $\rho(f)$ es el \'unico n\'umero en $[0,1)$ tal que $\bar{\rho}(F)-\rho(f)$ es entero.
\end{defi}

\begin{teo}\label{sub:definicionro}
Sea $f:S^1 \rightarrow S^1$ un difeomorfismo que preserva la orientaci\'on y $F$ un \emph{levantamiento} de $f$. Entonces $$\bar{\rho}(F)=\lim_{n\to \infty} \frac{F^n(x)}{n}$$ existe y es independiente de $x$. Por consiguiente, el n\'umero de rotaci\'on $\rho(f)$ est\'a bien definido y solo depende de la funci\'on.
\end{teo}

\begin{proof}
Veamos la existencia del  l\'imite. Dividamoslo en dos casos, cuando $f$ tiene puntos peri\'odicos y cuando no los tiene. Supongamos primero que $f$ tiene un punto peri\'odico de per\'iodo $m$, es decir, $f^m(\theta)=\theta$ con $\theta \in S^1$, luego existe $x \in \mathbb{R}$ tal que $\pi_1(x)=\theta$. Entonces para un \emph{levantamiento} cualquiera $F$ de $f$ sabemos que $F^n$ es tambi\'en un \emph{levantamiento} de $f^n$, as\'i
\begin{align*}
     \pi_1(F^m(x))&=(\pi_1\circ F^m)(x)\\
                  &=(f^m\circ \pi_1)(x)\\
                  &=f^m(\theta)  \\
                  &=\theta \\
                  &=\pi_1(x)              
\end{align*}
entonces $F^m(x)=x+k$ para alg\'un entero $k$. Luego, $F^{jm}(x)=\underbrace{F^m\circ F^m\circ \ldots \circ F^m}_{j-veces}(x)=x+jk$, y tenemos que $$\lim_{j\to \infty} \frac{|F^{jm}(x)|}{jm}=\lim_{j\to \infty} \left(\frac{x}{jm}-\frac{k}{m}\right)=\frac{k}{m}.$$

Adem\'as, escribimos cualquier entero $n$ en la forma $n=jm+r$ donde $0\leq r <m$. Podemos ver que hay una constante $M$ tal que $|F^r(y)-y|\leq M$ para todo $y \in \mathbb{R}$ y $0\leq r <m$. Con esto tenemos que $$\frac{|F^n(x)-F^{jm}(x)|}{n}=|\frac{F^r(F^{jm}(x))-F^{jm}(x)|}{n} \leq \frac{M}{n}.$$

As\'i  $$\lim_{n\to \infty} \frac{|F^n(x)-F^{jm}(x)|}{n}=0$$

\noindent y por tanto  \begin{equation}\label{sub:roracional} \bar{\rho}(F)=\lim_{n\to \infty} \frac{|F^n(x)|}{n}=\lim_{j\to \infty} \frac{|F^{jm}(x)|}{mj}=\frac{k}{m}.\end{equation}

Veamos ahora el caso cuando $f$ no tiene puntos peri\'odicos. Sabemos que $F^n(x)-x$ nunca es un entero si $n\neq0$, entonces siempre hay un entero $k_n$ tal que $$k_n<F^n(x)-x<k_n+1$$ para todo $x \in \mathbb{R}$. Aplicando repetidas veces la desigualdad anterior para $x=0,F^n(x),F^{2n}(x),\dots,$  tenemos
\begin{align*}
k_n<&F^n(0)<k_n+1 \\
k_n<F^{2n}(0)&-F^n(0)<k_n+1 \\
\vdots& \\
k_n<F^{mn}(0)&-F^{(m-1)n}(0)<k_n+1 \\                       
\end{align*}
Sumando tenemos $mk_n<F^{mn}(0)<m(k_n+1)$ y as\'i $\frac{k_n}{n}<\frac{F^{mn}(0)}{mn}<\frac{k_n+1}{n}.$\\

De la desigualdad original obtenemos  $\frac{k_n}{n}<\frac{F^{n}(0)}{n}<\frac{k_n+1}{n}.$ Combinando las desigualdades anteriores resulta $$|\frac{F^{mn}(0)}{mn}-\frac{F^n(0)}{n}|<\frac1n.$$

Esta misma desigualdad la podemos obtener con los argumentos $m,n$ intercambiados, por lo que tambi\'en es cierta $|\frac{F^{mn}(0)}{mn}-\frac{F^m(0)}{m}|<\frac1m$. Con las dos \'ultimas desigualdades tenemos $$|\frac{F^{n}(0)}{n}-\frac{F^m(0)}{m}|\leq|\frac{F^{n}(0)}{n}-\frac{F^{mn}(0)}{mn}|+|\frac{F^{mn}(0)}{mn}-\frac{F^m(0)}{m}| \leq \frac1n+\frac1m$$

As\'i, la sucesi\'on $\{\frac{F^n(0)}{n}\}_{n>0}$ es de Cauchy y por lo tanto converge. As\'i $$\bar{\rho}(F)=\lim_{n \to \infty} \frac{F^n(0)}{n}$$ existe (ya que  veremos a continuaci\'on que $\bar{\rho}(F)$ es independiente del $x$ que tomemos). Y hemos terminado la existencia. \\

Veamos que este l\'imite no depende de la elecci\'on de $x$. Para $x,y \in \mathbb{R}$ existe un entero $m$ tal que $|x+m-y|<1$, adem\'as por el Teorema~\ref{sub:propiedadeslifting} con $d=1$,  $F^n(x+m)-(x+m)=F^{n}(x)-x$, luego
\begin{align*}
|F^n(x)-F^n(y)|&=|(F^n(x)-x)-(F^n(y)-y)+(x-y)| \\
               & \leq |(F^n(x)-x)-(F^n(y)-y)|+|x-y| \\
               &=|(F^n(x+m)-(x+m))-(F^n(y)-y)|+|x-y| \\
               &=|(F^n(x+m)-F^n(y))-(x+m-y)|+|x-y| \\
               & \leq |F^n(x+m)-F^n(y)|+|x+m-y|+|x-y|, \tag*{por ecuaci\'on (~\ref{sub:desigualdad1})}\\
               & \leq 1+1+|x-y| 
\end{align*}
As\'i, dado $\epsilon>0$, para que $\frac{|F^n(x)-F^n(y)|}{n}\leq \frac{2+|x-y|}{n}<\epsilon$ es suficiente tomar $n> \frac{2+|x-y|}{\epsilon}$. Por lo tanto $$\lim_{n\to \infty} \frac{|F^n(x)-F^n(y)|}{n}=0$$
y con esto tenemos que $\bar{\rho}(F)$ es independiente del $x$ que tomemos.
\end{proof}

Este resultado tambi\'en permite establecer que $\rho(f)$ depende continuamente de $f$.

\begin{teo}\label{sub:continuidadro}
Supongamos que $f:S^1 \rightarrow S^1$  es un difeomorfismo que preserva la orientaci\'on. Para $\epsilon>0$, existe un $\delta>0$ tal que si $g:S^1 \rightarrow S^1$ es tambi\'en un difeomorfismo y $\sup_{x \in S^1} d(f(x),g(x))<\delta$ (donde $d$ es una distancia en $S^1$), entonces $|\rho(f)-\rho(g)|<\epsilon$.
\end{teo}
\begin{proof}
Elijamos $n$ tal que $\frac2n<\epsilon$. Podemos elegir un \emph{levantamiento} $F$ de $f$ con $r-1<F^n(0)<r+1$ para todo entero $r$. Podemos tambi\'en elegir $\delta>0$ suficientemente peque\~no para que exista un \emph{levantamiento} $G$ de $g$ con $r-1<G^n(0)<r+1$. Retomando las desigualdades de la prueba del Teorema~\ref{sub:definicionro}, tenemos $$|\frac{F^{nm}(0)}{nm}-\frac{G^{nm}(0)}{nm}|<\frac2n<\epsilon$$ 
para todo $m \in \mathbb{R}$. Sabiendo que $$\lim_{m\to \infty} |\frac{F^{nm}(0)}{nm}|=\bar{\rho}(F), $$ se tiene que $|\bar{\rho}(F)-\bar{\rho}(G)|<\epsilon$ y por tanto $|\rho(f)-\rho(g)|<\epsilon$.
\end{proof}

\subsection{Rotaciones en la circunferencia.}
En una rotaci\'on de la circunferencia vemos que el \'angulo de rotaci\'on en este tipo de funciones no es m\'as que el n\'umero de rotaci\'on del que hablamos en general.
\begin{coro}
Sea $f: S^1 \rightarrow S^1$ una rotaci\'on de la circunferencia por un \'angulo $0\leq \theta \leq 1$, es decir, $f(z)=e^{2 \pi \theta i}z$ con $z \in S^1$. Entonces $\rho(f)=\theta$.
\end{coro}
\begin{proof}
Sea $F$ un \emph{levantamiento} de $f$, entonces, $$e^{2 \pi iF(x)}=\pi_1(F^n(x))=f^n(\pi_1(x))=e^{2 \pi n \theta i}e^{2 \pi xi}=e^{2 \pi (n\theta+x) i}$$
as\'i, $F^n(x)=x+n\theta$, luego $$\bar{\rho}(F)=\lim_{n\to \infty} \frac{F^n(0)}{n}=\lim_{n\to \infty} \frac{n\theta}{n}=\theta$$

\noindent por lo tanto, $\bar{\rho}(F)=\rho(f)=\theta$.
\end{proof}

\begin{teo}\label{sub:orbitasperiodicas}
Sea $f: S^1 \rightarrow S^1$ una rotaci\'on de la circunferencia por un \'angulo $0\leq \theta \leq 1$.
\begin{itemize}
\item[a)] Si $\rho(f)=\frac kn$ con $k,n$ enteros coprimos, entonces $f^n(z)=z$ para todo $z \in S^1$, es decir, $f$ es globalmente peri\'odica de per\'iodo $n$.
\item[b)] Si $\rho(f)=r$ con $r$ irracional, entonces $f$ no tiene \'orbitas peri\'odicas y cada \'orbita es densa en $S^1$. 
\end{itemize}
\end{teo}
\begin{proof}
Probemos a), $f^n(z)=e^{2 \pi n \frac kn i} z=e^{2 \pi k i} z=z$, adem\'as $n$ es el menor entero para lo que esto sucede. Por tanto $f$ es globalmente peri\'odica de per\'iodo $n$. \\

Para la primera parte de b), supongamos que existe $z_0 \in S^1$ tal que su \'orbita tiene per\'iodo $n$, entonces 
$$f^n(z_0)=e^{2 \pi n r i} z_0=z_0$$ y de esto se tiene que $2 \pi n r= 2 \pi k$ para cualquier entero $k$, luego $r=\frac kn$ que es una contradicci\'on,  por lo tanto no existen \'orbitas peri\'odicas. Para ver la densidad, sea $z \in S^1$, si $f^n(z)=f^m(z)$ entonces $e^{2 \pi n r i} z=e^{2 \pi m r i} z$ y por tanto $(n-m)r \in \mathbb{Z}$, y esto s\'olo si $n=m$, as\'i tenemos que la \'orbita de $z$ es infinita y todos los valores son distintos dos a dos. Luego por ser la \'orbita infinita, dado $\epsilon>0$ cualquiera, debe existir dos enteros positivos $n,m$ tal  que $d(f^n(z),f^m(z))=|z|d(1,e^{2 \pi (n-m) i})<\epsilon$ ($d(x,y)$ distancia  sobre el arco en $S^1$ entre $x$ y $y$ tomada en sentido de las agujas del reloj). Tomando a $k=n-m$ (suponiendo sin p\'erdida de generalidad que $n>m$) tememos que $d(f^k(z),z)=|z|d(1,e^{2 \pi (n-m) i})<\epsilon$.\\

Sabemos que $f$ preserva las longitudes en $S^1$. Por consiguiente, $f$ transforma el arco entre $f^k(z)$ y $z$ al arco que conecta a $f^{2k}(z)$ y $f^k(z)$ que tambi\'en est\'an a una distancia menor que $\epsilon$. En particular la sucesi\'on $z,f^k(z),f^{2k}(z),\dots$ es una partici\'on de $S^1$ en arcos de longitud menor que $\epsilon$. As\'i, dado  $z_0 \in S^1$ y $\epsilon>0$ siempre existir\'a $n$ tal que $d(f^n(z),z_0)<\epsilon$. Por lo tanto, la \'orbita de $z$ es densa en $S^1$.
\end{proof}
\thispagestyle{headings}
\pagestyle{empty}
\cleardoublepage
\pagestyle{empty}
\pagestyle{headings}
\section{Din\'amica de la Ecuaci\'on de Lyness}

\subsection{Introducci\'on}

El problema planteado es describir el comportamiento de la sucesi\'on $x_1,x_2,x_3 \ldots $ de n\'umeros reales positivos generados por la ecuaci\'on

 \begin{equation}\label{sub:ecuacion principal} x_{n+1}=\frac{\alpha+x_n}{x_{n-1}}, n\geq2 \end{equation}
donde el parametro $\alpha$  es un n\'umero ral, $\alpha\geq0$, y los t\'erminos iniciales $x_1,x_2$ son n\'umeros reales positivos.\\ 

Este problema se va a resolver analizando el sistema din\'amico en $\mathbb{R}^2$ que resulta de desplegar esta ecuaci\'on. Tal desplegamiento se realiza como sigue. Sea $\mathbb{R}^+$  los n\'umeros reales $\mathbb{R}$ positivos. El aumento de la dimensi\'on es porque cada t\'ermino de la sucesi\'on depende de dos t\'erminos previos. El espacio en el que se desarrolla es el cuadrante positivo $\mathbb{R}^2_+$ de $\mathbb{R}^2$, dado por $x,y>0$. La ecuaci\'on se desplega en el difeomorfismo $f:\mathbb{R}^2_+\longrightarrow \mathbb{R}^2_+$ dado por 

\begin{equation}\label{sub:sistema dinamico} f(x,y)=(y,\frac{\alpha+y}{x}).\end{equation}

Dado $P=(x_1,x_2)\in \mathbb{R}^2_+$, la sucesi\'on determinada por $P$ es $\{x_1,x_2,x_3, \ldots, x_n, \ldots\}$ en $\mathbb{R}^+$ generada por $\{x_1,x_2\}$ y usando la ecuaci\'on en diferencias. El desplegamiento de una sucesi\'on es el conjunto $\{(x_1,x_2),(x_2,x_3),(x_3,x_4),\ldots,(x_n,x_{n+1}),\ldots \}$ que pertenecen a $\mathbb{R}^2_+$, lo cual no es otra cosa que la \'orbita $O_f(P)=\{P,f(P),f^2(P),\ldots,f^n(P), \ldots \}$. La proyecci\'on de $\mathbb{R}^2_+$ sobre el eje de las abscisas transforma la \'orbita en la sucesi\'on.\\

La \'orbita es el desplegamiento de la sucesi\'on, e inversamente la sucesi\'on es la proyecci\'on de la \'orbita. La gran ventaja del desplegamiento es que las \'orbitas en dos dimensiones son m\'as f\'aciles  de manejar y visualizar matem\'aticamente, mientras que sus proyecciones en una dimensi\'on son m\'as complicadas.

\subsection{Valores particulares del parametro $\alpha=0,1,\infty$}

Denotemos por $S=\{x_n\}_{n>0}$ la soluci\'on de (~\ref{sub:ecuacion principal}), generada a partir de los valores $x_1,x_2$.

\begin{lema} \label{sub:lemaalfa0}Si $\alpha=0$ entonces  $S$ es un 6-ciclo, es decir, todas las \'orbitas para $f$ son fijas o son de per\'iodo 6.\end{lema}
\begin{proof} La ecuaci\'on (~\ref{sub:ecuacion principal}) queda $$x_{n+1}=\frac{x_n}{x_{n-1}}$$ Sea $x_1=x$ y $x_2=y$ entonces \\

\begin{align*}
x_3&=\frac{y}{x},\,\, x_4=\frac{\frac{y}{x}}{y}=\frac{1}{x},\,\,x_5=\frac{\frac{1}{x}}{\frac{y}{x}}=\frac{1}{y},
x_6=\frac{\frac{1}{y}}{\frac{1}{x}}=\frac{x}{y},\,\,x_7=\frac{\frac{x}{y}}{\frac{1}{y}}= x,\,\,x_8=\frac{x}{\frac{x}{y}}=y.
\end{align*}

\noindent Por lo tanto $S=\{\underbrace{x,y,\frac{y}{x},\frac{1}{x},\frac{1}{y},\frac{x}{y}}_{6-ciclo},x,y, \ldots\}$, y tambi\'en se tiene, tomando a $P=(x_1,x_2)$, que $f^6(P)=P$.
\end{proof}

\begin{lema} \label{sub:lemaalfa1} Si $\alpha=1$ entonces  $S$ es un 5-ciclo, es decir, todas las \'orbitas para $f$ son fijas o son de per\'iodo 5.\end{lema}
\begin{proof} Sea $x_1=x$ y $x_2=y$ entonces \\

$x_1=\frac{1+y}{x}$,\\
 
$x_4=\frac{1+\frac{1+y}{x}}{y}=\frac{1+x+y}{xy}$,\\

$x_5=\frac{1+\frac{1+x+y}{xy}}{\frac{1+y}{x}}=\frac{1+x+y+xy}{y(1+y)}=\frac{(1+x)(1+y)}{y(1+y)}=\frac{1+x}{y}$,\\

$x_6=\frac{1+\frac{1+x}{y}}{\frac{1+x+y}{xy}}=\frac{(1+x+y)x}{1+x+y}=x$,\\

$x_7=\frac{1+x}{\frac{1+x}{y}}=y$,\\

\noindent Por lo tanto $S=\{\underbrace{x,y,\frac{1+y}{x},\frac{1+x+y}{xy}=\frac{1+x}{y}}_{5-ciclo},x,y, \ldots\}$, y tambi\'en se tiene, tomando a $P=(x_1,x_2)$, que $f^5(P)=P$.
\end{proof}

\begin{lema} \label{sub:lemaalfainfinito}Si $\alpha=\infty$ entonces  $S$ es un 4-ciclo, es decir, todas las \'orbitas para $f$ son fijas o son de per\'iodo 4.\end{lema}
\begin{proof}En la ecuaci\'on (~\ref{sub:ecuacion principal}) hacemos un cambio de variable $x_n=\sqrt{\alpha}y_n$, as\'i $\sqrt{\alpha}y_{n+1}=\frac{\alpha+\sqrt{\alpha}y_n}{\sqrt{\alpha}y_{n-1}}$, luego $y_{n+1}=\frac{1+\frac{1}{\sqrt{\alpha}}y_n}{y_{n-1}}=\frac{1}{y_{n-1}}$ por ser $\alpha=\infty$. Dados $y_1=x$ y $y_2=y$ tenemos \\
\begin{align*}
y_3&=\frac{1}{x},\,\,y_4=\frac{1}{y},\,\,y_5=x,\,\,y_6=y.
\end{align*}

\noindent Por lo tanto $\{\underbrace{x,y,\frac{1}{x},\frac{1}{y}}_{4-ciclo},x,y, \ldots\}$, as\'i, $\{y_n\}_{n> 0}$ tiene per\'iodo 4. Luego, esto tambi\'en es cierto para $\{x_n\}_{n\geq 0}$,  por lo que tomando a $P=(x_1,x_2)$ se tiene $f^4(P)=P$.
\end{proof}

En resumen, conocemos el comportamiento de (~\ref{sub:ecuacion principal}) para los citados valores de $\alpha$. Ahora nos planteamos el comportamiento para $0<\alpha<1$ y $1<\alpha<\infty$.

\subsectionmark{Comportamiento del desplegamiento de la ecuaci\'on}
\subsection{Comportamiento del desplegamiento de la ecuaci\'on de Lyness en $\mathbb{R}^2_+$.}
\subsectionmark{Comportamiento del desplegamiento de la ecuaci\'on}
La funci\'on $V$ del siguiente resultado fue introducida por G. Ladas en \cite{cox5}.
\begin{teo} 
La funci\'on $V:\mathbb{R}^2_+\longrightarrow \mathbb{R}$ dada por \begin{equation}\label{sub:invariante} V(x,y)=\frac{(x+1)(y+1)(x+y+\alpha)}{xy}\end{equation} es un invariante de $f$, es decir, $V(f(x,y))=V(x,y)$.
\end{teo}

\begin{proof}

\begin{align*} 
V(f(x,y)) & = V(y,\frac{\alpha+y}{x})      \\ 
          & = \frac{(y+1)(\frac{\alpha+y}{x}+1)(y+\frac{\alpha+y}{x}+\alpha)}{y(\frac{\alpha+y}{x})} \\
          & = \frac{(y+1)(\alpha+y+x)(xy+\alpha+y+\alpha x)}{xy(\alpha+y)} \\
          & = \frac{(y+1)(\alpha+y+x)(x+1)(\alpha+y)}{xy(\alpha+y)}   \\
          & = V(x,y)
\end{align*}
Por lo tanto, $V$ es un invariante de $f$.
\end{proof}

\begin{coro}\label{sub:corolariocurvasdenivel}Cada \'orbita de $f$ est\'a sobre una curva de nivel de $V$.\end{coro}
\begin{proof} Esto es una inmediata consecuencia del teorema anterior. Sea $C$ una curva de nivel de $V$(digamos $V=c,c\in \mathbb{R}$) y $(a,b)\in C$, entonces se verifica $V(a,b)=V(f(a,b))=V(b,\frac{\alpha+b}{a})=c$, por lo que $f(a,b)=(b,\frac{\alpha+b}{a})\in C$. Por lo tanto, cada \'orbita de $f$ est\'a sobre alguna curva de nivel de $V$. 
\end{proof}

\begin{figure}[h!] 
\begin{center}
\includegraphics[scale=.9]{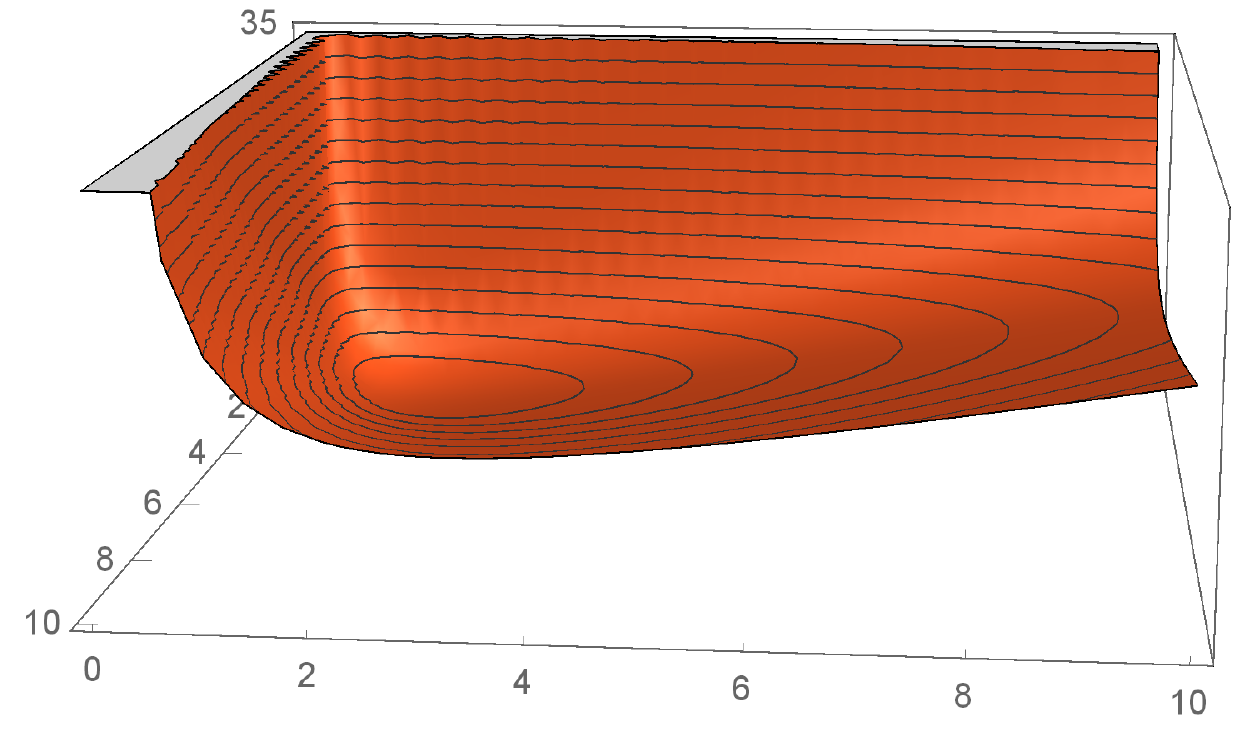}
\end{center}
\caption[Invariante $V$]{Gr\'afica de $V(x,y)$, con $(x,y) \in [0.002,200]\times[0.002,400]$ \label{sub:graficoinvariante}}
\end{figure}

\begin{teo}\label{sub:teorematazonondo}
La gr\'afica de $V$ tiene forma de <<taz\'on hondo>> (Ver figura~\ref{sub:graficoinvariante})  con  un m\'inimo en $F=(\omega,\omega)$, donde $$ \omega=\frac{1+\sqrt{1+4\alpha}}{2} \qquad   y \qquad   V(F)=V(\omega,\omega)=\frac{(\omega+1)^3}{\omega}.$$ Adem\'as, $V$ no tiene otro punto cr\'itico, y tiende a infinito en la frontera de $\mathbb{R}^2_+$.
\end{teo}

\begin{proof} 
Diferenciemos $V$:
$$\frac{\partial V}{\partial x}=\frac{(y+1)(x^2-y-\alpha)}{x^2y}\qquad   y \qquad \frac{\partial V}{\partial y}=\frac{(x+1)(y^2-x-\alpha)}{xy^2}$$
Los puntos cr\'iticos de $V$ son obtenidos por $$\frac{\partial V}{\partial x}=\frac{\partial V}{\partial y}=0.$$ Siendo $x,y>0$, se tiene que: $$x^2-y-\alpha=y^2-x-\alpha=0$$ Agrupando nos queda: 
$$x^2-y^2+x-y=(x-y)(x+y)+x-y=(x-y)(x+y+1)=0$$ por lo tanto $x=y$, ya que $x+y+1>0$. As\'i, $F=(x,y)=(\omega,\omega)$, donde $\omega^2-\omega-\alpha=0$.\\

Las ra\'ices del polinomio anterior son: $\omega=\frac{1\pm \sqrt{1+4 \alpha}}{2}$, y como $\omega>0$ porque $F \in \mathbb{R}^2_+$, resulta $$\omega=\frac{1+ \sqrt{1+4 \alpha}}{2}$$ y se tiene que $F$ es el \'unico punto cr\'itico de $V$. En $F$ utilizando que $\omega^2=\omega+\alpha$ tenemos:
\begin{align*} 
V(F)   &=\frac{(\omega+1)(\omega+1)(\omega+\omega+1)}{\omega\omega}=\frac{(\omega+1)^2(2\omega+1)}{\omega^2}=\frac{(\omega+1)^3}{\omega} \\
\frac{\partial^2V}{\partial x^2} \Big |_F &= \frac{2(y+1)(y+\alpha)}{x^3y} \Big |_F=\frac{2(\omega+1)(\omega+\alpha)}{\omega^4}=\frac{2(\omega+1)}{\omega^2}\\
\frac{\partial^2V}{\partial y^2} \Big |_F &= \frac{2(x+1)(x+\alpha)}{xy^3} \Big |_F=\frac{2(\omega+1)(\omega+\alpha)}{\omega^4}=\frac{2(\omega+1)}{\omega^2}\\
\frac{\partial^2V}{\partial x\partial y} \Big |_F &=\frac{\alpha-x^2-y^2}{x^2y^2}\Big |_F=\frac{\alpha-2\omega^2}{\omega^4}=-\frac{\omega+1}{\omega^3}
\end{align*} 

Calculamos la {\bf  Matriz Hessiana} de $V$
\begin{displaymath}
\mathbf{H(V)} =
\left( \begin{array}{cc}
\frac{\partial^2V}{\partial x^2} & \frac{\partial^2V}{\partial x\partial y}   \\
\frac{\partial^2V}{\partial y\partial x} & \frac{\partial^2V}{\partial y^2}  \\
\end{array} \right)
\end{displaymath}

\begin{displaymath}
\mathbf{H(V)\Big |_F} =
\left( \begin{array}{cc}
\frac{2(\omega+1)}{\omega^2} & -\frac{\omega+1}{\omega^3}   \\
-\frac{\omega+1}{\omega^3} & \frac{2(\omega+1)}{\omega^2}  \\
\end{array} \right)
\end{displaymath}

$$\big |H(V)|_F \big |  =\left(\frac{2(\omega+1)}{\omega^2}\right)^2-\left(-\frac{\omega+1}{\omega^3}\right)^2=\frac{(\omega+1)^2}{\omega^4}\left(4-\frac{1}{\omega^2}\right)$$
Pero $\omega^2=\frac{1+2\sqrt{1+4 \alpha}+1+4\alpha}{4}=\frac{1}{4}+\frac{1+2\sqrt{4 \alpha}+1+4\alpha}{4}>\frac14$, por lo tanto
$$\big |H(V) |_F \big |>0$$
y de aqu\'i que en $F$ hay un m\'inimo de $V$. Adem\'as, es f\'acil ver que 
$$
\lim_{\shortstack{$\scriptstyle x\rightarrow \infty$\\ $\scriptstyle y \rightarrow \infty$}}V(x, y), \lim_{\shortstack{$\scriptstyle x \rightarrow 0$}}V(x, y), \lim_{\shortstack{$\scriptstyle y \rightarrow 0$}}V(x, y),\lim_{\shortstack{$\scriptstyle x \rightarrow \infty$}}V(x, y),\lim_{\shortstack{$\scriptstyle y \rightarrow \infty$}}V(x, y)
$$
tienden a infinito, es decir, en la frontera de $\mathbb{R}^2_+$ la funci\'on $V$ tiende a infinito.
\end{proof}

\begin{figure}[h!]
\begin{center}
\includegraphics[width=8cm,height=8cm]{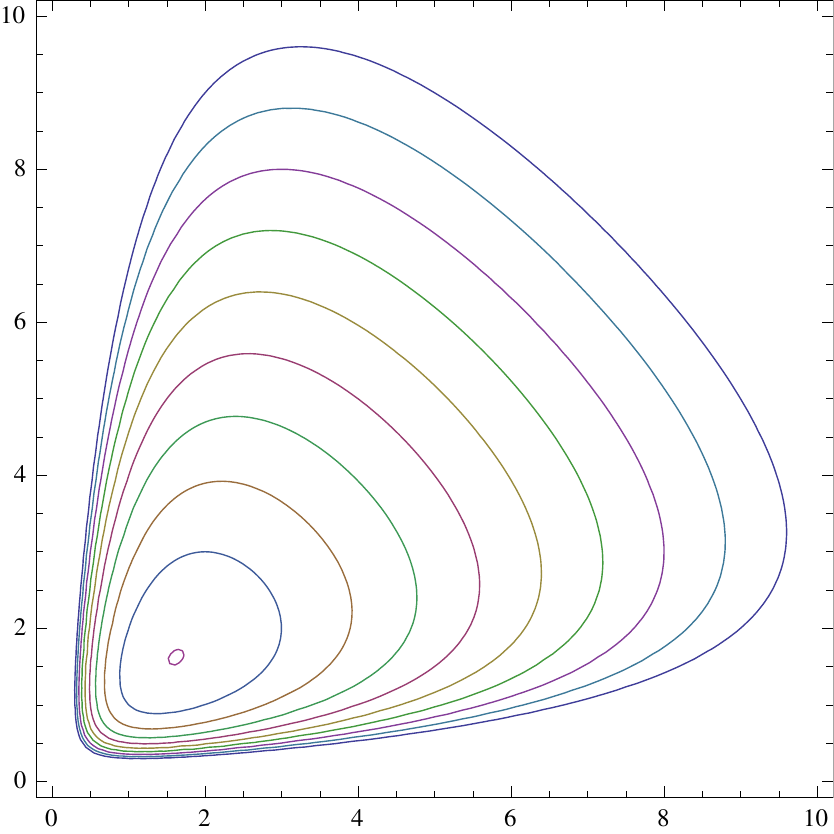}
\end{center}
\caption[Curvas de nivel de $V$]{Curvas de nivel de $V(x,y)$. \label{sub:curvasdenivel}}
\end{figure}

\begin{coro} Las curvas de nivel de la funci\'on $V$ forman un familia de curvas cerradas que contienen a $F$ y que \emph{llenan} a $\mathbb{R}^2_+$ (ver figura~\ref{sub:curvasdenivel}).\end{coro}

\begin{proof}
Es una simple consecuencia del teorema anterior, ya que $V$ tiene solo un punto cr\'itico $F$ que es un m\'inimo y tiende a infinito en la frontera de  $\mathbb{R}^2_+$. As\'i para cada curva de nivel $V=c$ con $c\geq V(F)$, se genera una curva cerrada  que contiene a $F$ y que \emph{llena} a  $\mathbb{R}^2_+$ (esto \'ultimo se puede comprender mejor porque cada \'orbita de $f$ est\'a sobre una curva de nivel de $V$ por Corolario~\ref{sub:corolariocurvasdenivel}, y entonces para cada $(x,y)\in \mathbb{R}^2_+ $ hay una \'unica curva de nivel de $V$).
\end{proof}

\begin{coro} El difeomorfismo $f$ tiene como \'unico punto fijo  a $F$ y transforma cada curva de nivel en s\'i misma.\end{coro}
\begin{proof}
$F$ es un punto fijo de $f$ $$f(\omega,\omega)=(\omega,\frac{\alpha+\omega}{\omega})=(\omega,\frac{\omega^2}{\omega})=(\omega,\omega).$$
Adem\'as, este es el \'unico punto fijo de $f$ porque si $(x,y)$ es punto fijo, entonces $f(x,y)=(y,\frac{y+\alpha}{x})=(x,y)$ de aqu\'i que se tiene$$x=y=\frac{\alpha+y}{x} $$ Obteniendo que $$x^2-x-\alpha=0 $$ por lo que $x=y=\omega$. Cada curva es transformada sobre si misma por el Corolario ~\ref{sub:corolariocurvasdenivel}.
\end{proof}

\begin{defi}
Dos homeomorfismos $\phi: A \rightarrow A$ y $\phi':B\rightarrow B$ son llamados \emph{conjugados}, y escribimos $\phi \sim \phi'$, si existe un homeomorfismo $h:A\rightarrow B$ tal que el diagrama siguiente conmuta
\end{defi}

\begin{center}
\begin{tabular}{ccc} 
$A$ & $\stackrel{\phi}{\rightarrow}$ & $A$ \\ 
 $\,\,\,\,\downarrow h$ &  & $\,\,\,\,\downarrow h$ \\ 
$B$ & $\stackrel{\phi'}{\rightarrow}$ & $B.$ 
\end{tabular}
\end{center}

\begin{defi}
Un homeomorfismo $\phi:A\rightarrow A$ de una curva cerrada $A$ es llamado \emph{casi-rotaci\'on} si es conjugado a una rotaci\'on eucl\'idea de la circunferencia a trav\'es de un \'angulo $2\pi\rho$, para alg\'un $\rho$, $0\leq \rho \leq 1$. \end{defi}

Al n\'umero $\rho$ lo llamamos \emph{n\'umero de rotaci\'on} y lo denotamos por $\rho(\phi)$ para hacer ver que depende de la funci\'on $\phi$. \\

{\bf Notaci\'on}:
Para hacer \'enfasis en el parametro $\alpha$ a\~nadiremos un sub\'indice como sigue: dado $\alpha$ sea $f_\alpha$, $F_\alpha$, $V_\alpha$ el correspondiente difeomorfismo, punto fijo y funci\'on invariante. Por el Teorema~\ref{sub:teorematazonondo} $F_\alpha=(\omega_\alpha,\omega_\alpha)$  donde  $\omega_\alpha=\frac{1+\sqrt{1+4\alpha}}{2}$.\\

Denotemos con $\rho_\alpha$  el n\'umero de rotaci\'on de la linealizaci\'on de $f_\alpha$ en $F_\alpha$. Sea $$v_\alpha=\min(V_\alpha)=V_\alpha (F_\alpha)=\frac{(\omega_\alpha+1)^3}{\omega_\alpha}.$$

Para $v>v_\alpha$, sea $C^v_\alpha$ la curva de nivel $V_\alpha=v$, y sea $\rho^v_\alpha$ el n\'umero de rotaci\'on de $f_\alpha|_{C^v_\alpha}$. Por el Teorema~\ref{sub:continuidadro} se tiene que  el n\'umero de rotaci\'on de una familia continua de difeomorfismos es continuo (para m\'as informaci\'on ver \cite{devaney} y \cite{LJ}), entonces $\rho^v_\alpha$ es una funci\'on continua de $v$ en el dominio $v_\alpha<v<\infty$.\\
  
El siguiente teorema  es de gran importancia porque es el que nos permite caraterizar como se comportan las \'orbitas de la funci\'on $f$. Este resultado es el que nos facilita el desarrollo satisfactorio del an\'alisis de la ecuaci\'on planteada al inicio. Como los casos para $\alpha=0,\infty$ son tratados en el Teorema~\ref{sub:teoremadeparametros}, supondremos que $0<\alpha<\infty$ y $v_\alpha<v<\infty$. 

\begin{teo}
\label{sub:casirotacion}Si $C$ es una curva de nivel de $V_\alpha$ entonces $f_\alpha|_C$ es una \emph{casi-rotaci\'on}.
\end{teo}

\begin{proof}
Sea $f=f_\alpha$. Si $C$ es una curva de nivel de $V_\alpha$ entonces $C$ es una curva c\'ubica en $\mathbb{R}^2_+$, obtenida por $$(x+1)(y+1)(x+y+1)-vxy=0.$$

Sea  $\bar{C}$ la completaci\'on de $C$ en el \emph{Plano Proyectivo Complejo}\footnote{El  \emph{Plano Proyectivo Complejo} $\mathbb{C}\mathbb{P}^2$ es el conjunto de l\'ineas a trav\'es del origen en $\mathbb{C}^3$. Un punto en $\mathbb{C}\mathbb{P}^2$ es representado por un tr\'io de coordenadas $(x,y,z)$ de n\'umeros complejos, no todos ceros, o equivalentemente por cualquier m\'ultiplo no cero de este tr\'io.} $\mathbb{C}\mathbb{P}^2$. La ecuaci\'on de  $\bar{C}$ es obtenida sustituyendo $x,y$ por $x/z,y/z$ en la ecuaci\'on de $C$ y multiplicando por $z^3$, tenemos la ecuaci\'on c\'ubica homog\'enea $$(x+z)(y+z)(x+y+\alpha z)-vxyz=0.$$

\noindent La completaci\'on  $\bar{f}$ de $f$ es obtenida por la f\'ormula $$\bar{f}(x,y,z)=\bar{f}(\frac x z,\frac y z, 1)=(\frac y z, \frac{\alpha+\frac y z}{\frac x z},1)=(xy, (y+\alpha z)z,xz).$$ Esta f\'ormula define $\bar{f}$ sobre todo $\mathbb{C}\mathbb{P}^2$ excepto en los tres puntos $(1,0,0)$, $(0,1,0)$ y $(0,-\alpha,0)$, donde no est\'a definido porque la f\'ormula da $(0,0,0)$. Daremos, sin embargo, una representaci\'on geom\'etrica de $\bar{f}|_{\bar{C}}$ que extiende la definici\'on de $\bar{f}$ a estos tres puntos, y hace de $\bar{f}|_{\bar{C}}$ un difeomorfismo de la curva $\bar{C}$ en s\'i misma.\\

Una \emph{involuci\'on} en un espacio es una funci\'on del espacio en s\'i mismo de per\'iodo 2 (ver ejemplo~\ref{sub:invo1} y~\ref{sub:invo2}). La completaci\'on $\bar{f}|_{\bar{C}}$ es la composici\'on de dos involuciones sobre $\bar{C}$, que son la completaci\'on de las dos involuciones para $f|_C$ del ejemplo~\ref{sub:invo2}
 
$$\xrightarrow[(x,y,z)\,\stackrel{\bar{a}^*_1}{\rightarrow}\,((y+\alpha z)z,xy,xz)\, \stackrel{\bar{a}^*_2}{\rightarrow}\,(xy,(y+\alpha z)z,xz)]{\bar{f}|_{\bar{C}}}$$

\noindent que son inducidas por los puntos $\bar{a}_1=(1,0,0)$ y $\bar{a}_2=(1,-1,0)$ de $\bar{C}$, porque la completaci\'on de l\'ineas horizontales en $\mathbb{R}^2_+$ son l\'ineas en $\mathbb{C}\mathbb{P}^2$ a trav\'es de $\bar{a}_1$, mientras que la completaci\'on de l\'ineas perpendiculares a la diagonal en $\mathbb{R}^2_+$ son l\'ineas a trav\'es de $\bar{a}_2$ en  $\mathbb{C}\mathbb{P}^2$. Por lo tanto, $\bar{f}|_{\bar{C}}$ es un difeomorfismo de $\bar{C}$ que preserva la orientaci\'on.\\

Denotemos con  $\mathbb{Z}^2$ a todos los puntos con coordenadas enteras de $\mathbb{R}^2$. Por el siguiente Lema~\ref{sub:curvaeliptica2} sabemos que $\bar{C}$ es una curva el\'iptica que es un grupo abeliano \cite{Sil}, y la estructura de este grupo puede representarse del siguiente modo
$$\mathbb{R}^2 \xrightarrow[proyecci\acute{o}n]{\pi} \underbrace{\mathbb{R}^2/\mathbb{Z}^2}_{toro} \xrightarrow[difeomorfismo]{h} \bar{C}$$

\noindent donde $\pi$ denota la proyecci\'on  en el grupo cociente, y $h$ es un difeomorfismo.\\

Tres puntos  $\bar{a},\bar{b},\bar{c} \in \bar{C}$ son colineales en $\mathbb{C}\mathbb{P}^2$ si y s\'olo si $\bar{a}+\bar{b}+\bar{c}=O$, donde $O$ es elemento neutro del grupo $\bar{C}$ (esto se verifica f\'acilmente por la definici\'on de la operaci\'on del grupo $\bar{C}$). Luego, $h^{-1}(O)=h^{-1}(\bar{a}+\bar{b}+\bar{c})=[0] \in \mathbb{R}^2/\mathbb{Z}^2$. Si $a,b,c$ son las correspondientes preim\'agenes de $\bar{a},\bar{b},\bar{c}$ por medio de $h \circ \pi$, entonces $\pi(a+b+c)=[0]$ y por tanto $a+b+c \in \mathbb{Z}^2$. As\'i, con la representaci\'on $h\circ \pi$ tenemos la siguiente propiedad

\begin{center}$\bar{a},\bar{b},\bar{c} \in \bar{C}$ son colineales en $\mathbb{C}\mathbb{P}^2$ si y s\'olo si $a+b+c \in \mathbb{Z}^2$.\end{center}

Por lo tanto la involuci\'on 
\begin{center}
\begin{tabular}{ccc} 
$\bar{a}^*:\bar{C}$ & $\rightarrow$ & $\bar{C}$ \\  
\,\,\,\,\,\,\,$\bar{b}$           & $\rightarrow$ & $\bar{c}$ 
\end{tabular}
\end{center}
inducida por $\bar{a} \in \bar{C}$ conmuta con la funci\'on lineal 
\begin{center}
\begin{tabular}{ccc} 
$a^*:\mathbb{R}^2$ & $\rightarrow$ & $\mathbb{R}^2$ \\  
\,\,\,\,\,\,\,$b$           & $\rightarrow$ & $c=-a-b$ 
\end{tabular}
\end{center}

\noindent por medio de $h\circ \pi$, es decir,

\begin{center}
\begin{tabular}{ccc} 
$\mathbb{R}^2$ & $\stackrel{a^*}{\rightarrow}$ & $\mathbb{R}^2$ \\ 
 $\,\,\,\,\downarrow h\circ \pi$ &  & $\,\,\,\,\downarrow h\circ \pi$ \\ 
$\bar{C}$ & $\stackrel{f}{\rightarrow}$ & $\bar{C}.$ 
\end{tabular}
\end{center}

As\'i, sabiendo que $\bar{f}|_{\bar{C}}=\bar{a}^*_2\circ \bar{a}^*_1$, \'esta conmuta con la funci\'on lineal $\varphi=a^*_2\circ a^*_1$ 
$$\xrightarrow[b\,\stackrel{a^*_1}{\rightarrow}\,-a_1-b\, \stackrel{a^*_2}{\rightarrow}\,-a_2-(-a_1-b).]{\varphi}$$

Ahora $-a_2-(-a_1-b)=b+T$, donde $T$ es el vector fijo $T=a_1-a_2$. Por lo tanto, $\varphi$ no es m\'as que la traslaci\'on de $\mathbb{R}^2$ por el vector $T$.\\

Si $(x,y,z) \in C\subset \bar{C}$ entonces $z=1$, luego $\bar{f}(C)=f(C)=C$, as\'i $f|_C=\bar{f}|_C$ conmuta con $\varphi$ por medio de $h \circ \pi$. $C$ es difeomorfa por $h$ a una circunferencia y \'esta es difeomorfa a una recta $L$ en $\mathbb{R}^2$ por $\pi$, as\'i

$$L \xrightarrow[proyecci\acute{o}n]{\pi} \underbrace{(L+\mathbb{Z}^2)/\mathbb{Z}^2}_{circunferencia} \xrightarrow[difeomorfismo]{h} C.$$

Si $L$ es  paralela al vector $T$, $\varphi(L)=L$. Para $a+sT \in L$ con $s \in \mathbb{R}$ tenemos $\varphi(a+sT)=T+a+sT=a+(s+1)T \in L$, es decir, $\varphi|_L$ es una traslaci\'on sobre $L$. As\'i

\begin{center}
\begin{tabular}{ccc} 
$L$ & $\stackrel{\varphi}{\rightarrow}$ & $L$ \\ 
$\downarrow \pi$ &  & $\downarrow  \pi$ \\ 
$(L+\mathbb{Z}^2)/\mathbb{Z}^2$ & $\stackrel{\gamma}{\rightarrow}$ & $(L+\mathbb{Z}^2)/\mathbb{Z}^2$ \\
$\downarrow h$ &  & $\downarrow h$ \\ 
$C$ & $\stackrel{f}{\rightarrow}$ & $C.$ 
\end{tabular}
\end{center}

Donde $\gamma$ es una rotaci\'on eucl\'idea  de la circunferencia inducida por la traslaci\'on $\varphi|_L$. Por lo tanto, $f|_C$ es conjugado  por $h$ a una rotaci\'on eucl\'idea de la circunferencia. 
\end{proof}

\begin{figure}[h!]\centering
\includegraphics[width=6cm,height=6cm]{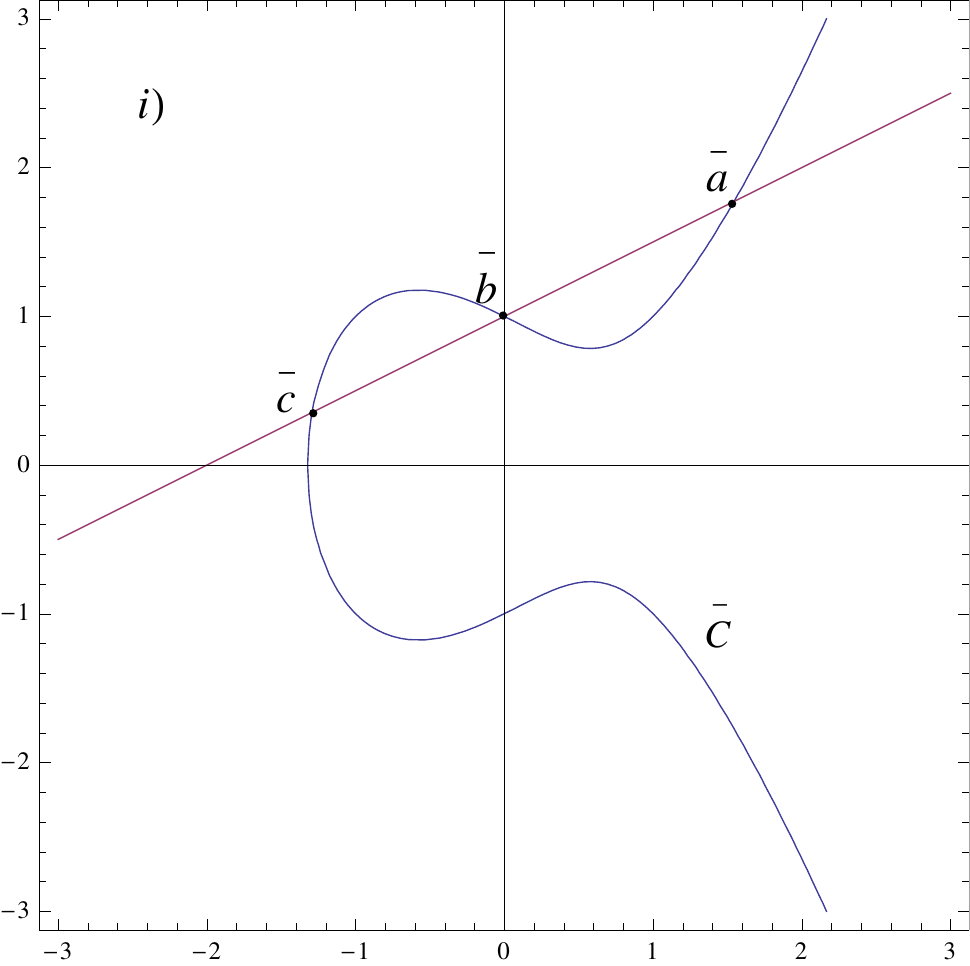} \quad \quad \includegraphics[width=6cm,height=6cm]{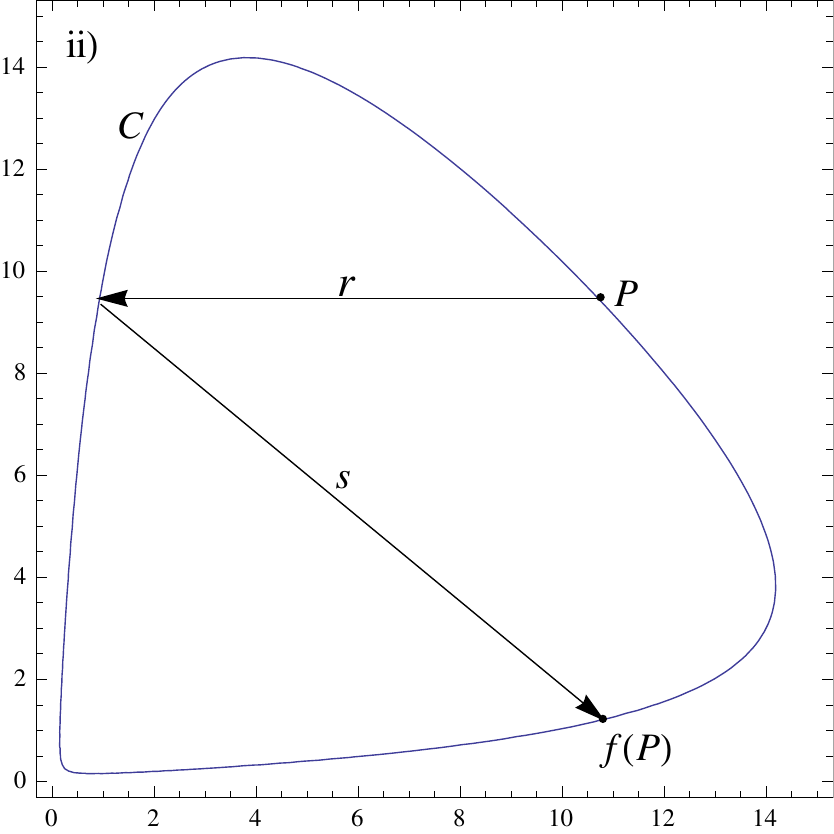} 
\caption[Involuciones.]{i) Involuci\'on sobre $\bar{C}$ inducida por $\bar{a}$. ii) $f|_C$ es la composici\'on de dos involuciones. \label{sub:figeliptica}}
\end{figure}

\begin{eje} \label{sub:invo1}
$\bar{C}$ es una curva el\'iptica por el siguiente Lema~\ref{sub:curvaeliptica2}, as\'i cada punto $\bar{a} \in \bar{C}$ induce una involuci\'on de $\bar{C}$ 
\begin{center}
\begin{tabular}{ccc} 
$\bar{a}^*:\bar{C}$ & $\rightarrow$ & $\bar{C}$ \\  
\,\,\,\,\,\,\,$\bar{b}$           & $\rightarrow$ & $\bar{c}$ 
\end{tabular}
\end{center}
donde  $\bar{c}$ es obtenido por la tercera intersecci\'on\footnote{La tercera intersecci\'on es infinito cuando la recta no corta la curva m\'as de dos veces, y si solo la corta una vez este punto se considera doble y el tercero sigue siendo infinito.} con $\bar{C}$ de la recta que pasa por $\bar{a}$ y $\bar{b}$ (ver figura ~\ref{sub:figeliptica}, i)).
\end{eje}

\begin{eje} \label{sub:invo2}
Otro ejemplo es el que muestra la figura~\ref{sub:figeliptica}, ii). En el que  $f|_C$ es la composici\'on de las siguientes dos funciones

$$\xrightarrow[(x,y)\,\stackrel{r}{\rightarrow}\,(\frac{y+\alpha}{x},y)\, \stackrel{s}{\rightarrow}\,(y,\frac{y+\alpha}{x})]{f}$$
     
\noindent donde $r$ es obtenida al intercambiar los extremos de la cuerda horizontal, y $s$ es obtenida por reflexi\'on en la diagonal. Cada una de estas funciones son involuciones y por tanto $f$ es la composici\'on de dos involuciones. 
\end{eje}

\begin{lema}  \label{sub:curvaeliptica2}
$\bar{C}$ es curva el\'iptica.
\end{lema}
\begin{proof}
Las curvas c\'ubicas que no son singulares se denominan \emph{el\'ipticas} (ver \cite{Sil}). Los casos excepcionales de c\'ubicas singulares es cuando son ``descomponibles'' o ``racionales''. Una curva c\'ubica es llamada descomponible si el polinomio que la define se puede factorizar en tres factores lineales distintos (representado tres l\'ineas) o en el producto de un factor lineal y uno cuadr\'atico (representado una l\'inea y una c\'onica). Una curva c\'ubica no descomponible se llama racional si tiene una singularidad que es un nodo (un punto doble) o una c\'uspide (un punto triple).\\

Para la prueba del lema, es suficiente demostrar que $\bar{C}$ no es ni descomponible ni racional. Esto lo haremos mostrando que el haz de curvas c\'ubicas obtenidas por la variaci\'on de $v$ como parametro en la l\'inea proyectiva compleja da lugar a cinco casos excepcionales donde $\bar{C}$ es singular, y que correponden a cinco valores de $v$ que est\'an fuera del dominio real ($v_\alpha<v<\infty$) ya establecido.\\

Veamos en dos casos que $\bar{C}$ no es descomponible:
\begin{itemize}
\item[a)] Descomposici\'on de $\bar{C}$ en tres l\'ineas. Si $v=0,\infty$, entonces $\bar{C}$ se descompone en tres l\'ineas como puede verse f\'acilmente de la ecuaci\'on de $\bar{C}$, pero ambos valores est\'an fuera del dominio. No puede haber  otra c\'ubica que se descomponga en tres l\'ineas, si esto se diera, entonces $\bar{f}$  permutar\'ia los tres puntos de la intersecci\'on, y es f\'acil ver que $\bar{f}$ no tiene otras \'orbitas de per\'iodo tres.

\item[b)] Descomposici\'on de $\bar{C}$ en una l\'inea m\'as una c\'onica. Si $v=\alpha-1$, entonces $\bar{C}$  se descompone en una l\'inea m\'as una c\'onica $$(x+y+z)(xy+yz+xz+\alpha z^2)=0.$$

Pero el valor $v=\alpha-1$ est\'a fuera del dominio de $v$ porque
\begin{align*}
v_\alpha&=\frac{(\omega+1)^3}{\omega}, \tag*{donde $\omega=\frac{1+\sqrt{1+4 \alpha}}{2}$,}\\
        &>(\omega+1)^2, \tag*{ya que $\frac{\omega+1}{\omega}>1$,}\\
        &>\frac{1+4 \alpha}{4},\tag*{$\omega+1>\omega>\frac{\sqrt{1+4 \alpha}}{2}$,}\\
        &>\alpha-1.
\end{align*}
No puede haber otra curva en el haz que se descomponga en una l\'inea m\'as una c\'onica, si esto se diera, entonces $\bar{f}$ intercambiar\'ia los  dos puntos de las intersecciones, y es f\'acil verificar que $\bar{f}$ no tiene otras \'orbitas de per\'iodo 2.
\end{itemize}

Si la curva c\'ubica en el haz es racional entonces estas singularidades son puntos fijos para $\bar{f}$, y es f\'acil verificar que $\bar{f}$ tiene solo dos puntos fijos. Uno es $\bar{F_1}=(\omega,\omega,1)$  en $v=v_\alpha$, que est\'a fuera del dominio. El otro es $\bar{F_2}=(\omega',\omega',1)$, donde $\omega'=\frac{1-\sqrt{1+4 \alpha}}{2}$, que es la otra ra\'iz de $\omega^2=\omega+\alpha.$ El correspondiente valor $v=v'_\alpha$ tambi\'en est\'a fuera del dominio porque $$v'=\frac{(\omega'+1)^3}{\omega'}=\frac{(\omega+1)^3}{\omega}-\frac{(1+4\alpha)^{3/2}}{\alpha}<v_\alpha.$$

Esto completa la lista de los cinco casos excepcionales  de las curvas c\'ubicas en el haz que no son el\'ipticas, pero $\bar{C}$ no tiene ninguno de estos casos porque como mostramos sus correspondientes valores no est\'an en el dominio para $v$. Por la tanto, $\bar{C}$ es una curva el\'iptica. 
\end{proof}

\begin{coro}
\label{sub:orbitasperiodicas}Si el n\'umero de rotaci\'on $\rho(f_\alpha|_C)$ es un n\'umero racional $p/q$(con $p$ y $q$ coprimos) entonces toda \'orbita de $f$ sobre $C$ es peri\'odica con per\'iodo $q$. Si este n\'umero es irracional entonces todas las \'orbitas de $f$ sobre $C$ son no peri\'odicas y son densas en $C$.
\end{coro}
\begin{proof}
En una rotaci\'on de la circunferencia las \'orbitas son todas peri\'odicas o densas seg\'un el n\'umero de rotaci\'on sea racional o irracional por el Teorema~\ref{sub:orbitasperiodicas}.
\begin{itemize}
\item[a)] Si $\rho(f|_C)=p/q$ con $p,q \in \mathbb{N}$ y coprimos entonces $f|_C$ es conjugado a una rotaci\'on de la circunferencia $\phi$ peri\'odica de per\'iodo $q$, luego  existe un homeomorfismo $h:S^1\rightarrow C$ tal que $f\circ h=h\circ \phi$ (recordemos que $S^1$ es la circunferencia undidad). Sea $\bar{x} \in C$ entonces existe $x_0 \in S^1$ tal que $h(x_0)=\bar{x}$. Sea $x_0,x_1,\dots,x_{q-1}$ la \'orbita de $x_0$ en $S_1$ por medio de $\phi$. Luego, 
\begin{align*} 
f(h(x_0))    & = h(\phi(x_0))=h(x_1)     \\ 
f(h(x_0))    & =h(\phi(x_1))=h(x_2) \\
             &  \vdots \\
f(h(x_{q-2}))& = h(\phi(x_{q-2}))=h(x_{q-1}) \\
f(h(x_{q-1}))& = h(\phi(x_{q-1}))=h(x_0)
\end{align*}
Por lo tanto, la \'orbita de $\bar{x}=h(x_0)$ es peri\'odica de per\'iodo $q$ en $C$ por medio de $f$. 
\item[b)] Si $\rho(f|_C)$ es irracional entonces $f|_C$ es conjugado a una rotaci\'on de la circunferencia $\phi$ densa en $S^1$, luego  existe un homeomorfismo $h:S^1\rightarrow C$ tal que $f\circ h=h\circ \phi$. Sea $\bar{x} \in C$ y $v$ un entorno cualquiera en $C$ (este entorno es tomado en la topolog\'ia heredada como subespacio de $\mathbb{R}^2$). Existe un entorno $u$ en $S^1$ tal que $h(u)=v$ y $x \in S^1$ tal que $h(x)=\bar{x}$. Como $\{\phi(x)\}_{n\geq0}$ es densa en $S^1$, entonces existe $N_0 \in \mathbb{N}$ tal que $\phi^{N_0}(x) \in u$, luego $h(\phi^{N_0}) \in v$, as\'i, $f^{N_0}(\bar{x})=f^{N_0}(h(x))=h(\phi^{N_0}(x)) \footnote{Si $f\circ h=h\circ \phi$, como h es homeomorfismo $f=h\circ \phi \circ h^{-1}$, as\'i, $f^n=\underbrace{(h\circ \phi \circ h^{-1})(h\circ \phi \circ h^{-1})\dots (h\circ \phi \circ h^{-1})}_{n-veces}=h\circ \phi^n \circ h^{-1}$ agrupando adecuadamente.} \in v$. Por lo tanto, $\{f^n(\bar{x})\}_{n\geq0}$ es densa en $C$. 
\end{itemize}

\end{proof}

\begin{figure}[h!]
\begin{center}
\includegraphics[width=5cm,height=5cm]{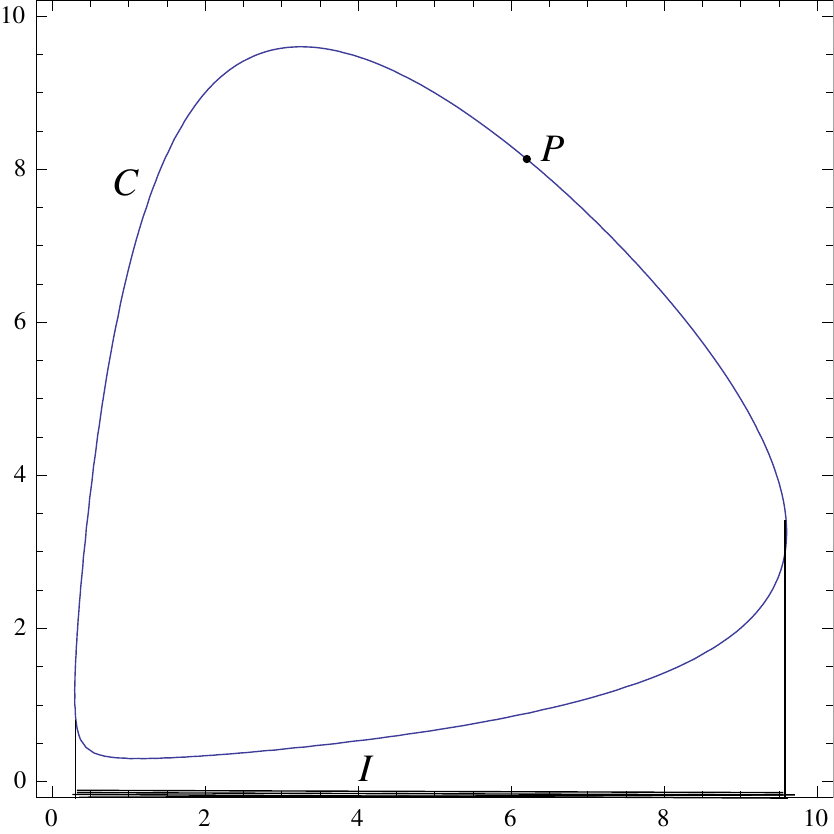}  \includegraphics[width=7cm,height=5cm]{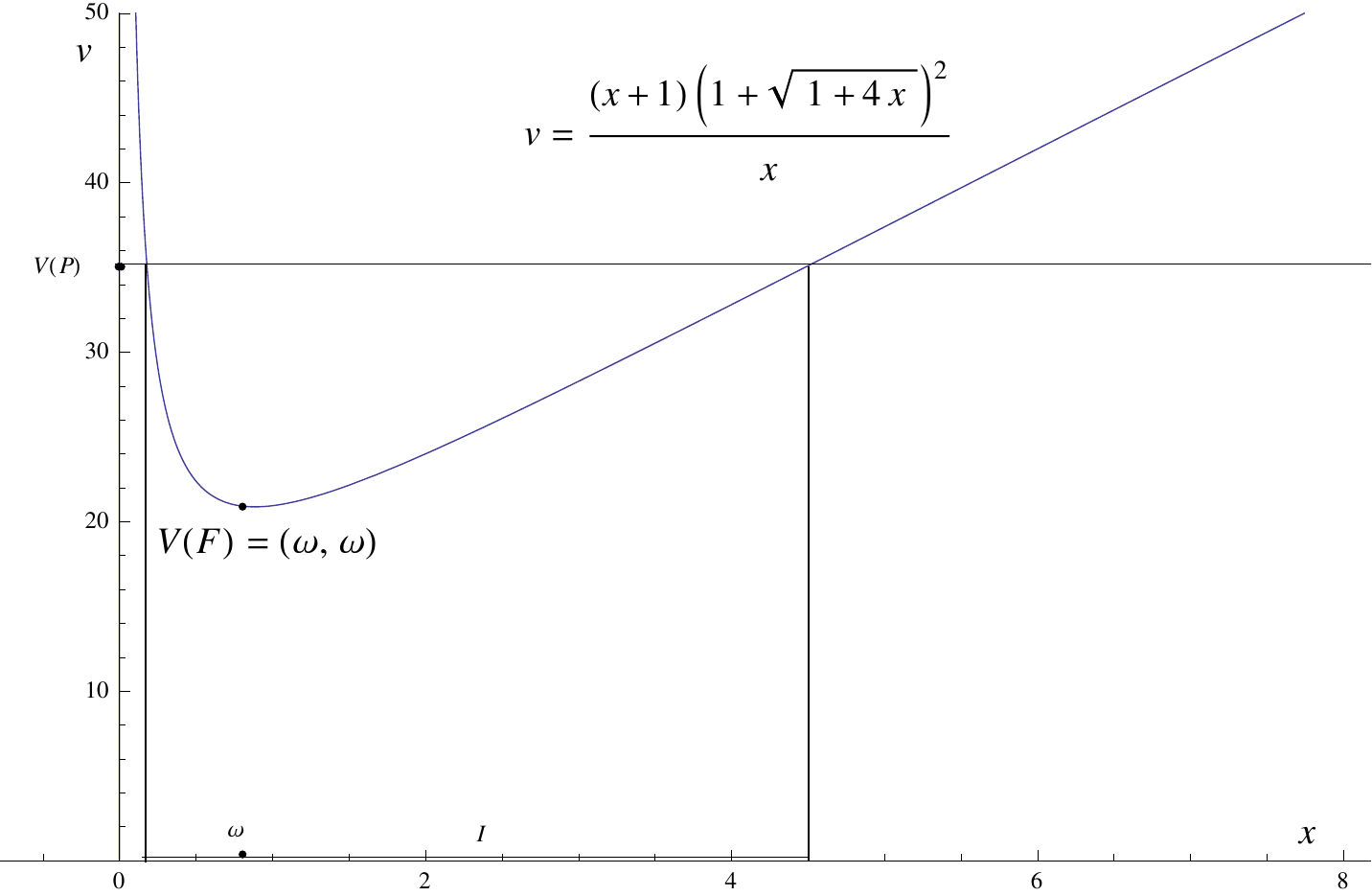}
\end{center}
\caption[Intervalo I]{El intervalo $I$ acota la sucesi\'on proyectada sobre el eje de las abscisas de la \'orbita determinada por $P$ en $\mathbb{R}^2_+$. \label{sub:intervaloi}}
\end{figure}

\begin{teo}
Dado $P\in C \subset \mathbb{R}^2_+$, denotemos por $S$ la sucesi\'on en $\mathbb{R}_+$ determinada por la proyecci\'on de la \'orbita de $P$. Sea $I$ la proyecci\'on de $C$  en el eje de las abscisas, que es un intervalo cerrado y acotado(ver figura~\ref{sub:intervaloi}). Entonces $S\subset I$, y $S$ es peri\'odica o densa en $I$. Adem\'as, los extremos de $I$ son el $\sup S$ e $\inf S$ para las sucesiones densas y estos extremos son obtenidos por la ecuaci\'on $$V(P)=\frac{(x+1)(1+\sqrt{x+\alpha})^2}{x}.$$
\end{teo}
\begin{proof}
Dado $P \in C \subset \mathbb{R}^2_+$, la \'orbita de $P$ est\'a en $C$ por el Corolario~\ref{sub:corolariocurvasdenivel}, la sucesi\'on $S$ es la proyecci\'on de la \'orbita determinada por $P$ y el intervalo $I$ es la proyeci\'on de $C$, por lo tanto $S \subset I$. La periodicidad o densidad es consecuencia del Corolario~\ref{sub:orbitasperiodicas}, ya que si una \'orbita es peri\'odica o densa en $C$ su proyeci\'on debe ser peri\'odica o densa en $I$, respectivamente.\\

Probemos ahora que $I=[a,b]$ nos da el $\sup S$ e $\inf S$. Supongamos que existen $r_1,r_2 \in \mathbb{R}$ tal que $\forall s \in S$, $r_1\leq s \leq r_2$ y $a<r_1<r_2<b$, entonces existen entornos $c_1,c_2$ en $C$ (en la topolog\'ia heredada como subespacio de $\mathbb{R}^2$) tal que sus proyecciones son los intervalos $]a,r_1[$ y $]r_2,b[$ respectivamente. Como $]a,r_1[\cap S$ y $]r_2,b[ \cap S$ son vac\'ios entonces $c_1 \cap O_f(P)$ y $c_2 \cap O_f(P)$ tambi\'en son vac\'ios, contradiciendo el hecho que $O_f(P)$ es densa en $C$. Por lo tanto, $a=\inf S$ y $b=\sup S$.\\

Establezcamos la ecuaci\'on para obtener los extremos de $I$; la normal a $C$ es 
 $$\nabla V=\left(\frac{\partial V}{\partial x},\frac{\partial V}{\partial y}\right)=\left(\frac{(y+1)(x^2-y-\alpha)}{x^2y},\frac{(x+1)(y^2-x-\alpha)}{xy^2}\right)$$
y es horizontal si $\frac{\partial V}{\partial y}=0$, es decir, $y^2=x+\alpha$. Luego, los puntos $P=(x,y) \in C$ donde la tangente es vertical son obtenidos cuando $y^2=x+\alpha$ (ya que $x,y>0$, $x+1\neq 0$). Sustituyendo en la ecuaci\'on de $C$ tenemos 
\begin{align*} 
    V(P)& = V(x,y)     \\ 
        & =\frac{(x+1)(y+1)(x+y+\alpha)}{xy} \\
        & =\frac{(x+1)(y+1)(y^2+y)}{xy} \\
        & = \frac{(x+1)(y+1)^2}{x} \\
        & = \frac{(x+1)(1+\sqrt{x+\alpha})^2}{x}.
\end{align*}
 
\end{proof}

\begin{teo}\label{sub:teoremadeparametros} Para los casos $\alpha=0,1,\infty$ el difeomorfismo $f_\alpha$ es per\'iodico de peri\'odo 6, 5, 4, respectivamente, y conjugado a una rotaci\'on del plano. Por lo tanto, toda las \'orbitas son peri\'odicas con peri\'odos 6,5,4, respectivamente, y todas las curvas de nivel tienen n\'umero de rotaci\'on $1/6$, $1/5$, $1/4$, respectivamente.
\end{teo}
\begin{proof}
Probar que $f_\alpha$ con $\alpha=0,1,\infty$ es peri\'odica con per\'iodo 6,5,4, respectivamente es consecuencia inmediata de los Lemas~\ref{sub:lemaalfa0},~\ref{sub:lemaalfa1} y~\ref{sub:lemaalfainfinito}. Faltar\'ia  mostrar que los difeomorfismos son conjugados a una rotaci\'on del plano. Haremos con detalle la prueba para $f_1$ y los cambios que hay que hacer en \'esta para obtener los otros dos resultados.\\

Para $f=f_1$: Necesitamos construir un homeomorfismo $h:\mathbb{C} \rightarrow \mathbb{R}^2_+$ tal que el diagrama siguiente conmute
\begin{center}
\begin{tabular}{ccc} 
$\mathbb{C}$ & $\stackrel{\phi}{\rightarrow}$ & $\mathbb{C}$ \\ 
 $\,\,\,\,\downarrow h$ &  & $\,\,\,\,\downarrow h$ \\ 
$\mathbb{R}^2_+$ & $\stackrel{f}{\rightarrow}$ & $\mathbb{R}^2_+$ 
\end{tabular}
\end{center}
donde $\phi$ denota la rotaci\'on $\phi(z)=ze^{\frac{-2\pi i}{5}}$ (el signo menos en la f\'ormula no es indispensable, pero se lo pondremos para que rote igual que $f$ en el sentido de las agujas del reloj y as\'i tener una mejor vizualizaci\'on). \\

\begin{figure}[h!]\centering
\includegraphics[width=6cm,height=6.5cm]{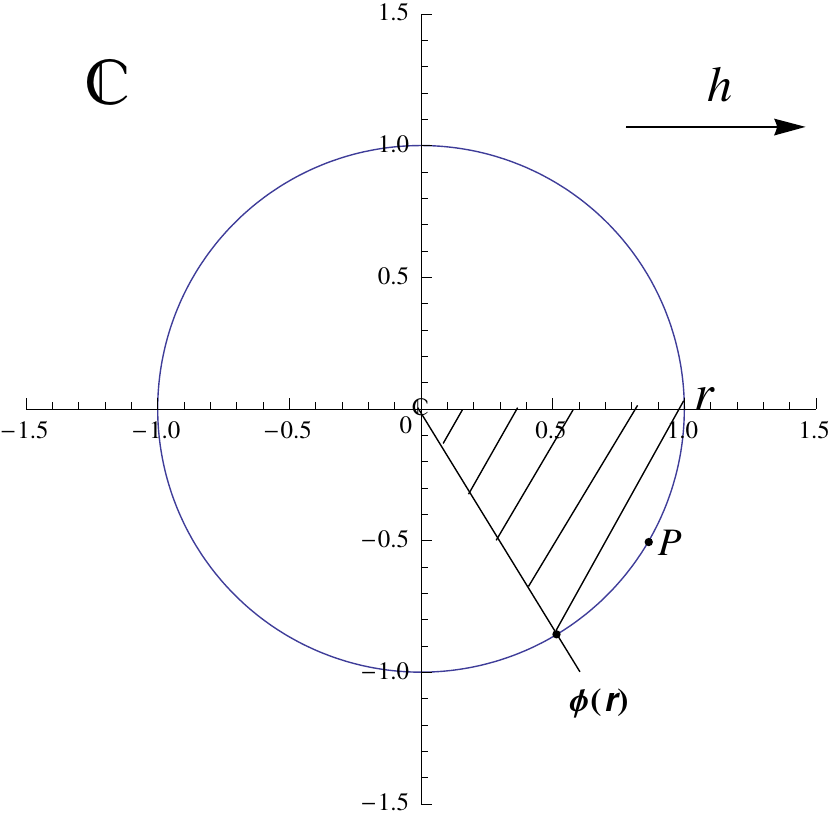}  \includegraphics[width=6cm,height=6.5cm]{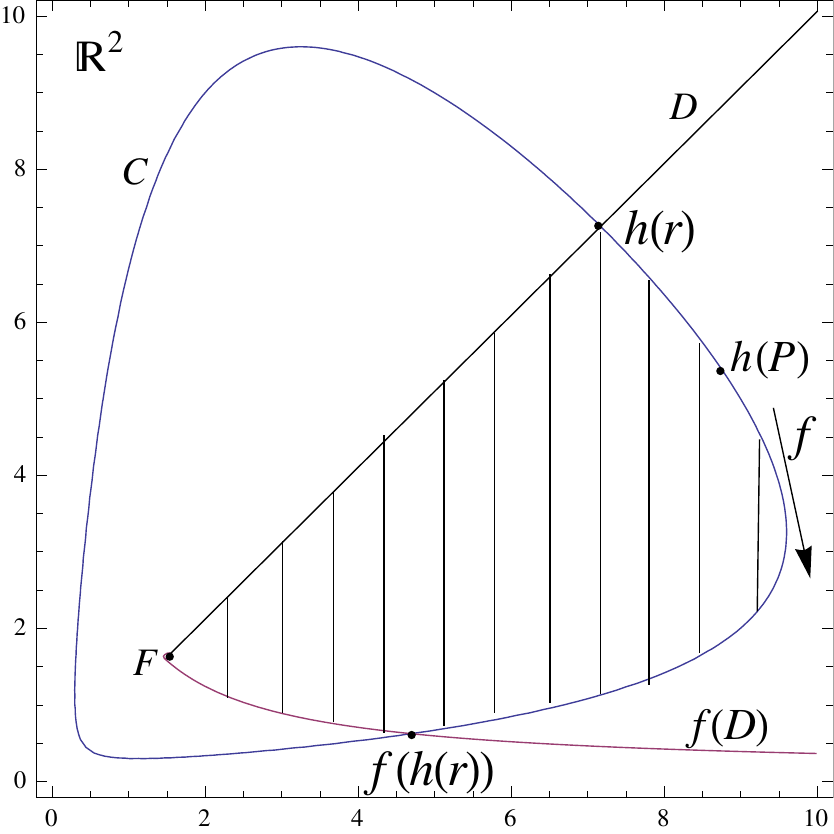}
\caption[Conjugaci\'on entre $\mathbb{C}$ y $\mathbb{R}^2_+$.]{Conjugaci\'on entre $f$ y la rotaci\'on $\phi$, $P=re^{-i\theta}$. }
\end{figure}

\begin{itemize}
\item[a)] Sea $h(0)=F=(\omega,\omega)$, donde $\omega=\frac{1+\sqrt{5}}{2}$, para $r\geq0$ y $\theta=0$, $h(r)=(\omega+r,\omega+r)$ y $h(\phi(r))=f(h(r))$, as\'i tenemos definida a $f$ sobre las semirectas que se generan con $\theta=0,\frac{2\pi}{5}$.
\item[b)] Para $r\geq0$ y $0< \theta < \frac{2\pi}{5}$, sea $h(P)=h(re^{-i\theta})=(x_\theta,y_\theta)\in C$ (donde $C$ es la curva de nivel de $V$ que pasa por el punto $(\omega+r,\omega+r)$), tal que $d_C((x_\theta,y_\theta),h(r))=\frac{5\theta}{2\pi}d_C(f(h(r)),h(r))$ elegido en la direcci\'on de las agujas del reloj a partir de $h(r)$, donde $d_C$ es la distancia sobre la curva $C$.  Como $0\leq \frac{5\theta}{2\pi} \leq 1$, as\'i $(x_\theta,y_\theta) \in \widehat{h(r)f(h(r))}$ (arco sobre $C$). El punto $(x_\theta,y_\theta)$ es \'unico porque la distancia depende linealmente del parametro $\theta$ y se ha elegido un sentido sobre la curva.
\item[c)] Para $r\geq 0$ y $\frac{2\pi q}{5}\leq \theta \leq \frac{2\pi (q+1)}{5}$ con $q=1,2,3,4$, sea $h(re^{-i\theta})=(f^q\circ h\circ  \phi^{-q})(re^{-i\theta})$, donde $\phi^{-1}(z)=ze^{\frac{2\pi i}{5}}$.
\end{itemize}

As\'i, para $P=re^{-i\theta}$, con $r\geq0$, se tiene:
$$ 
h(P) = \begin{cases} 
         h_0(P)=(x_\theta,y_\theta), & 0\leq \theta \leq \frac{2\pi}{5}\,\,y\,\,(x_\theta,y_\theta)\,\,como\,\,se\,\,defini\'o\,\,en\,\,b)\\ 
         h_1(P)=(f\circ h_0\circ \phi^{-1})(P), & \frac{2\pi}{5}\leq \theta \leq \frac{4\pi}{5} \\
         h_2(P)=(f^2\circ h_0\circ \phi^{-2})(P), & \frac{4\pi}{5}\leq \theta \leq \frac{6\pi}{5}\\
         h_3(P)=(f^3\circ h_0\circ \phi^{-3})(P), & \frac{6\pi}{5}\leq \theta \leq \frac{8\pi}{5}\\
         h_4(P)=(f^4\circ h_0\circ \phi^{-4})(P), & \frac{8\pi}{5}\leq \theta \leq \frac{10\pi}{5}\\
       \end{cases} 
$$ 
Por tanto, $h$ est\'a bien definida en todo $\mathbb{C}$, ya que en cada parte depende por medio de las funciones $f$ y $\phi$ de la definici\'on en a) y b) donde si est\'a bien definida. \\

Veamos que $f\circ h=h\circ \phi$  o equivalentemente que $h=f\circ h\circ \phi^{-1}$:
\begin{itemize}
\item Si $\frac{2\pi}{5}\leq \theta \leq \frac{4\pi}{5}$ entonces $h_1(P)=(f\circ h_0\circ \phi^{-1})(P)$ y cumple trivialmente la propiedad.
\item Si $\frac{4\pi}{5}\leq \theta \leq \frac{6\pi}{5}$ entonces $h_2(P)=(f^2\circ h_0\circ \phi^{-2})(P)=(f\circ(f\circ h_0\circ \phi^{-1})\circ \phi^{-1})(P)=(f\circ h_1\circ \phi^{-1})(P)$.
\item Si $\frac{6\pi}{5}\leq \theta \leq \frac{8\pi}{5}$ entonces $h_3(P)=(f^3\circ h_0\circ \phi^{-3})(P)=(f\circ (f^2\circ h_0\circ \phi^{-2})\circ \phi^{-1})(P)=(f\circ h_2\circ \phi^{-1})(P)$.
\item Si $\frac{8\pi}{5}\leq \theta \leq \frac{10\pi}{5}$ entonces $h_4(P)=(f^4\circ h_0\circ \phi^{-4})(P)=(f\circ (f^3\circ h_0\circ \phi^{-3})\circ \phi^{-1})(P)=(f\circ h_3\circ \phi^{-1})(P)$.
\item Si $0\leq \theta \leq \frac{2\pi}{5}$ entonces $h_0(P)=(f^5\circ h_0\circ \phi^{-5})(P)$ porque tanto $f$ como $\phi$ son de per\'iodo 5, luego $h_0(P)=(f\circ (f^4\circ h_0\circ \phi^{-4})\circ \phi^{-1})(P)=(f\circ h_4\circ \phi^{-1})(P)$.
\end{itemize}

As\'i, la funci\'on $h$ cumple las propiedades requeridas para garantizar que $f \sim \phi$ y adem\'as nos da el n\'umero de rotaci\'on que es justamente $1/5$ para todas la curvas de nivel generadas para $\alpha=1$.\\

Para terminar la prueba,  es suficiente decir que para $f_0$ y $f_\infty$ se toman $\phi(z)=ze^{\frac{-2\pi i}{6}}$ y $\phi(z)=ze^{\frac{-2\pi i}{4}}$ respectivamente, y se hacen exactamente los mismos pasos que para $f_1$.
\end{proof}

{\bf Importante:} los casos $\alpha=0,1,\infty$ son excepcionales porque en todos los otros casos $f_\alpha$ es una transformaci\'on que gira a $\mathbb{R}^2_+$, pero rotando diferentes curvas a diferentes velocidades, como lo comprobaremos en los posteriores teoremas. Si $0<\alpha<1$ entonces $f_\alpha$ rota las curvas m\'as r\'apido, y si $1<\alpha<\infty$ entonces las rota m\'as despacio. Mientras que es los casos $\alpha=0,1,\infty$ $f_\alpha$ rota a todas las curvas a la misma valocidad, porque el n\'umero de rotaci\'on es el mismo para todas las curvas de nivel.

\begin{teo}
\begin{enumerate}\label{sub:cotasro}
\item[a)]Si $0<\alpha<1$ entonces $1/6<\rho^v_\alpha<1/5$, para todo $v$, $v_\alpha<v<\infty$. 
\item [b)]Si $1<\alpha<\infty$ entonces $1/5<\rho^v_\alpha<1/4$, para todo $v$, $v_\alpha<v<\infty$.
\end{enumerate}
\end{teo}
\begin{proof}
Supongamos primero $1<\alpha<\infty$. Sea $C=C^v_\alpha$, $\rho=\rho^v_\alpha$. Sea $\bar{C}$ la completaci\'on de $C$ en $\mathbb{C}\mathbb{P}^2$, como en la demostraci\'on del Teorema~\ref{sub:casirotacion}. Sea $\Delta$ los tres ejes en $\mathbb{C}\mathbb{P}^2$ obtenidos por la condici\'on $xyz=0$, y $\Delta'$ las tres l\'ineas obtenidas por $$(x+z)(y+z)(x+y+\alpha z)=0.$$

Entonces el haz de curvas c\'ubicas todas pasan a trav\'es de $\Delta\cap\Delta'$ que consiste de los siete puntos $X=(1,0,0),Y=(0,1,0),Z=(1,-1,0),M=(0,-1,1),N=(-1,0,1),S=(0,-\alpha,1),T=(-\alpha,0,1)$, como ilustra la figura ~\ref{sub:planoproyectivo}.

\begin{figure}[h!]
\begin{center}
\includegraphics[width=8cm,height=8cm]{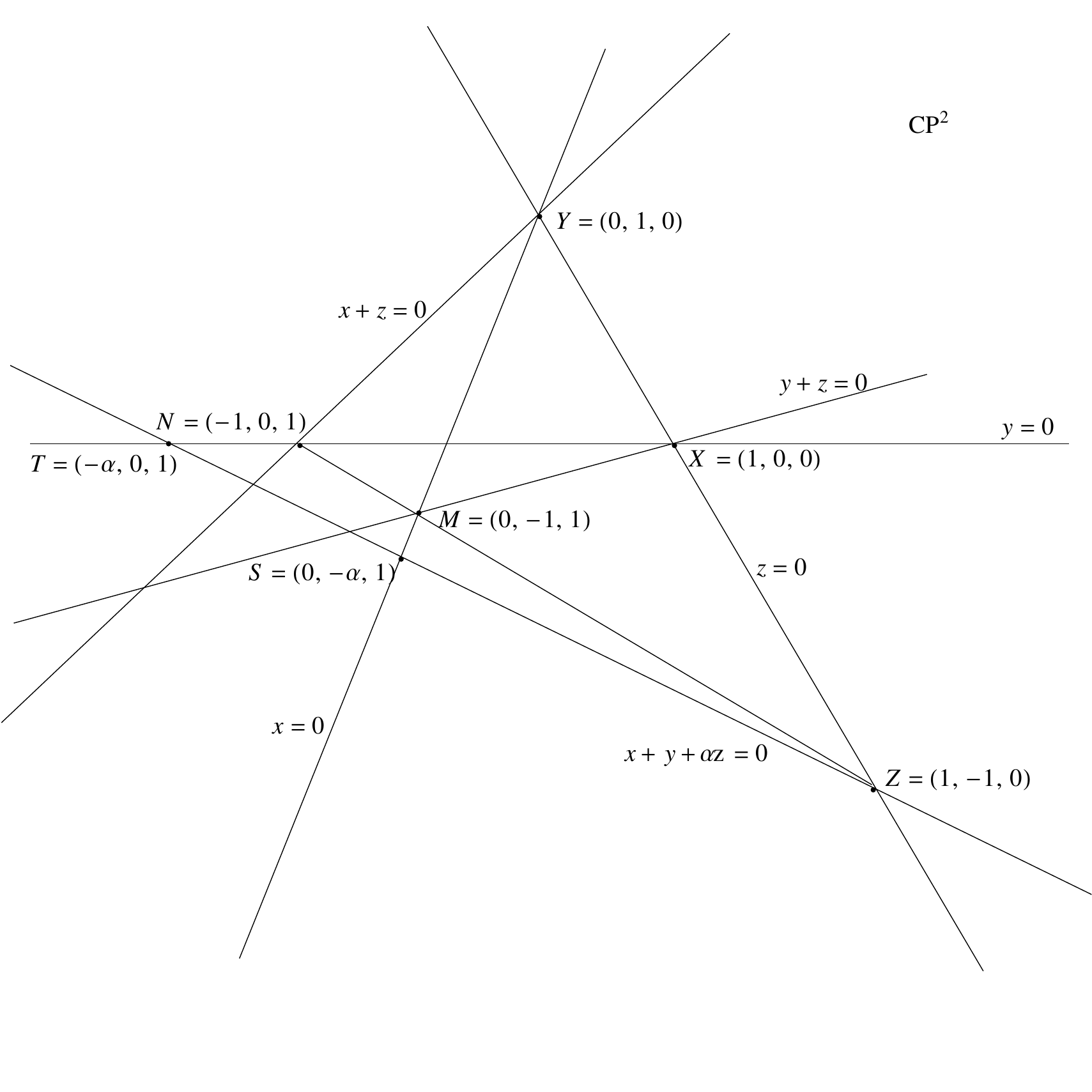} 
\end{center} 
\caption[Plano Proyectivo Complejo $\mathbb{C} \mathbb{P}^2$.]{Los siete puntos comunes del haz de curvas c\'ubicas. \label{sub:planoproyectivo}}
\end{figure}

Adem\'as, $X,Y$ son puntos dobles de $\Delta\cap\Delta'$, las curvas c\'ubicas todas tocan a $MX$ en $X$ y a $NY$ en $Y$. Pretendemos que la \'orbita de $T$ por medio de $\bar{f}$ contega los siete puntos en le siguiente orden $$T\rightarrow M \rightarrow Y \rightarrow Z \rightarrow X \rightarrow N \rightarrow S \rightarrow \ldots.$$

Para probar esto usamos la representaci\'on de $\bar{f}|_{\bar{C}}$ como la composici\'on $Z^* \circ X^*$ de las dos involuciones inducidas sobre $\bar{C}$ por $X$ y $Z$ (en la prueba del Teorema~\ref{sub:casirotacion} las llamamos $\bar{a}^*_1$ y $\bar{a}^*_2$ respectivamente):

\begin{align*}
T & \stackrel{X^*}{\rightarrow} & N & \stackrel{Z^*}{\rightarrow} & M, \\
M &\rightarrow                  & X & \rightarrow              & Y, \\
Y &\rightarrow                  & Z & \rightarrow              & Z,  \tag*{siendo  $Z$ un punto de inflexi\'on de $\bar{C}$,} \\
Z &\rightarrow                  & Y & \rightarrow              & X, \\
X &\rightarrow                  & M & \rightarrow              & N,  \tag*{sabiendo que $MX$ toca a $\bar{C}$ en $X$,}\\
N &\rightarrow                  & T & \rightarrow              & S.
\end{align*}

Con la f\'ormula para $\bar{f}$ se verifican las im\'agenes de $T,M,Z,N$.  La intersecci\'on de $\bar{C}$ con el plano proyectivo real $\mathbb{R}\mathbb{P}^2$ es la uni\'on de dos curvas disjuntas $C$ y  $C'$, de lo anterior $C$ es la que ya conocemos alrededor de $\mathbb{R}^2_+$, y la otra es  $C'$ que est\'a contenida en el complemento $\mathbb{R}\mathbb{P}^2-\mathbb{R}^2_+$ que contiene la \'orbita de $T$. Siguiendo la notaci\'on de la prueba del Teorema~\ref{sub:casirotacion} podemos encontrar una l\'inea $L'$ paralela $L$ en $\mathbb{R}^2$ tal que la traslaci\'on $\varphi|_{L'}$ conmute con $\bar{f}|_{C'}$ por medio de  $h\circ \pi$. La traslaci\'on $\varphi$ es la misma sobre ambas l\'ineas $L,L'$ y por consiguiente induce la misma rotaci\'on sobre ambos c\'irculos $\pi(L),\pi(L')$, entonces los difeomorfismos $f|_C$ y $\bar{f}|_{C'}$ son conjugados a la misma rotaci\'on y por tanto tienen el mismo n\'umero de rotaci\'on $\rho$. \\

Relacionamos a $\rho$ con la \'orbita de $T$, como sigue. La l\'inea en el infinito $z=0$ en $\mathbb{R}\mathbb{P}^2$ pasa por $C'$ en los tres puntos $X,Y,Z$. Eliminamos la l\'inea en el infinito de $\mathbb{R}\mathbb{P}^2$ y nos da $\mathbb{R}^2$; eliminamos $X,Y,Z$ de $C'$ y nos da $C' \cap \mathbb{R}^2$, que consiste en tres componentes que tienden al infinito en las direcciones $X,Y$ y $Z$ a lo largo de las as\'intotas $y=-1,x=-1$ y $x+y+\alpha=v$, respectivamente, como se ilustra en la figura~\ref{sub:planoreal} i).\\

\begin{figure}[h!]
\begin{center}
\includegraphics[width=8cm,height=8cm]{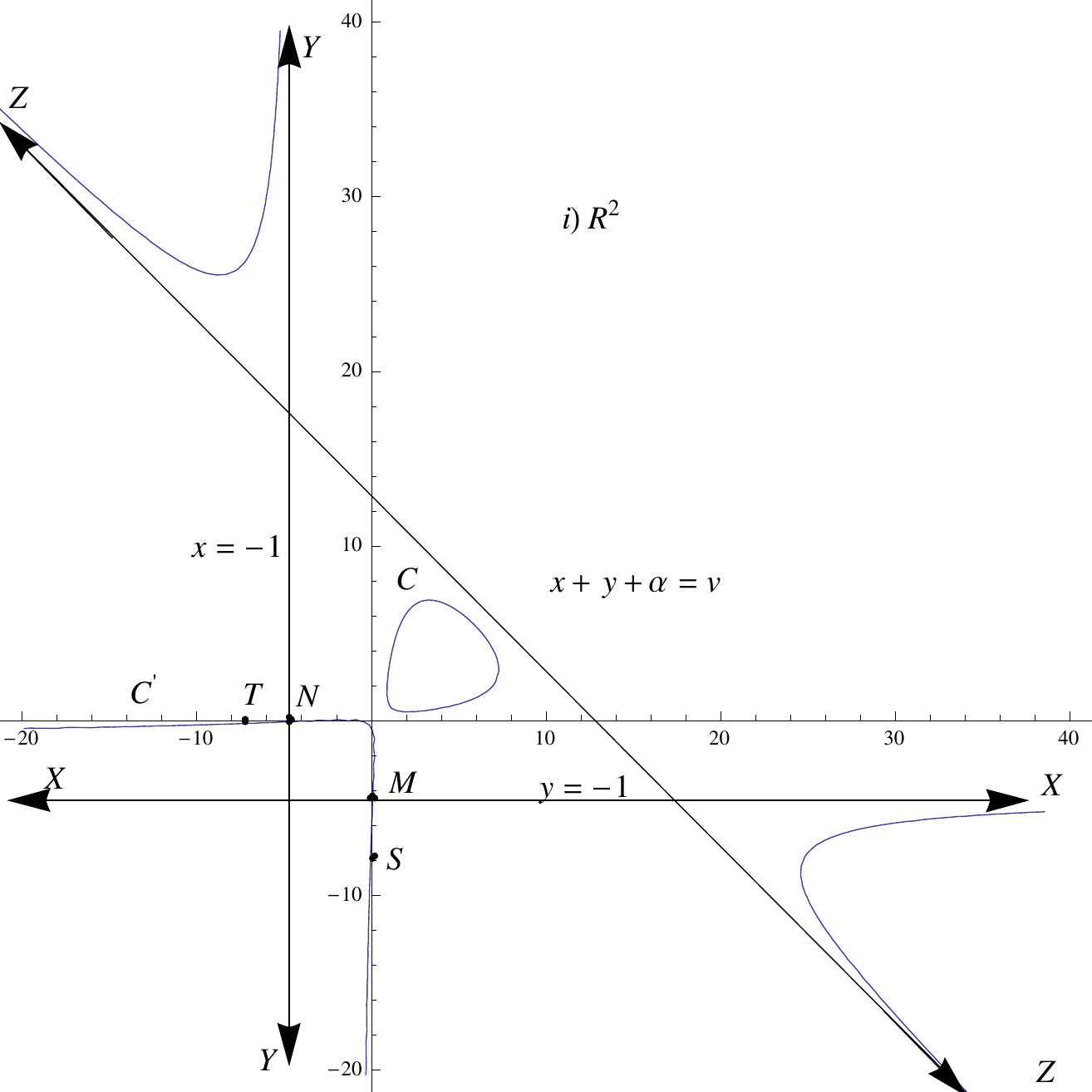}  \includegraphics[width=4cm,height=4cm]{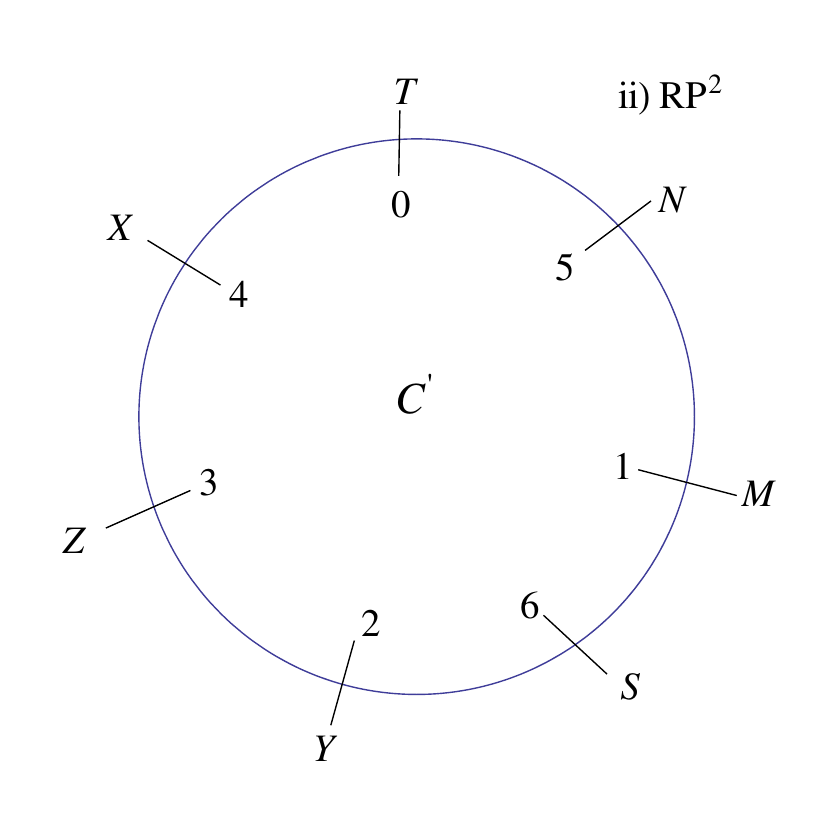}
\end{center} 
\caption[Curva $C$ en $\mathbb{R}^2$ cuando $1<\alpha<\infty$.]{ La \'orbita de $T$ sobre $C'$ en i) $\mathbb{R}^2$, y ii) $\mathbb{R}\mathbb{P}^2$, cuando $1<\alpha<\infty$. \label{sub:planoreal}}
\end{figure}

Reemplazamos las tres direcciones de las as\'intotas por $X,Y,Z$ para obtener de nuevo $C'$ en  $\mathbb{R}\mathbb{P}^2$, e identificamos el orden $TNMSYZX$ de los siete puntos alrededor de la curva $C'$, como se muestra en la figura~\ref{sub:planoreal} ii). Los n\'umeros dentro de la curva indican el orden $TMYZXNS$ de los puntos de la \'orbita de $T$. Vemos que  $T$ est\'a entre  $X=\bar{f}^4(T)$ y $Y=\bar{f}^5(T)$ sobre $C'$, entonces $4\rho<1<5\rho$. Por lo tanto $\frac15<\rho<\frac14$, como deseabamos.\\

En el otro caso, $0<\alpha<1$, las posiciones de $T,N$ son intercambiadas en el orden de los puntos alrededor de $C'$, de forma semejante las posiciones de $S,M$ son tambi\'en intercambiadas, como ilustra en la figura~\ref{sub:planoreal2} i).\\

\begin{figure}[h!]
\begin{center}
\includegraphics[width=8cm,height=8cm]{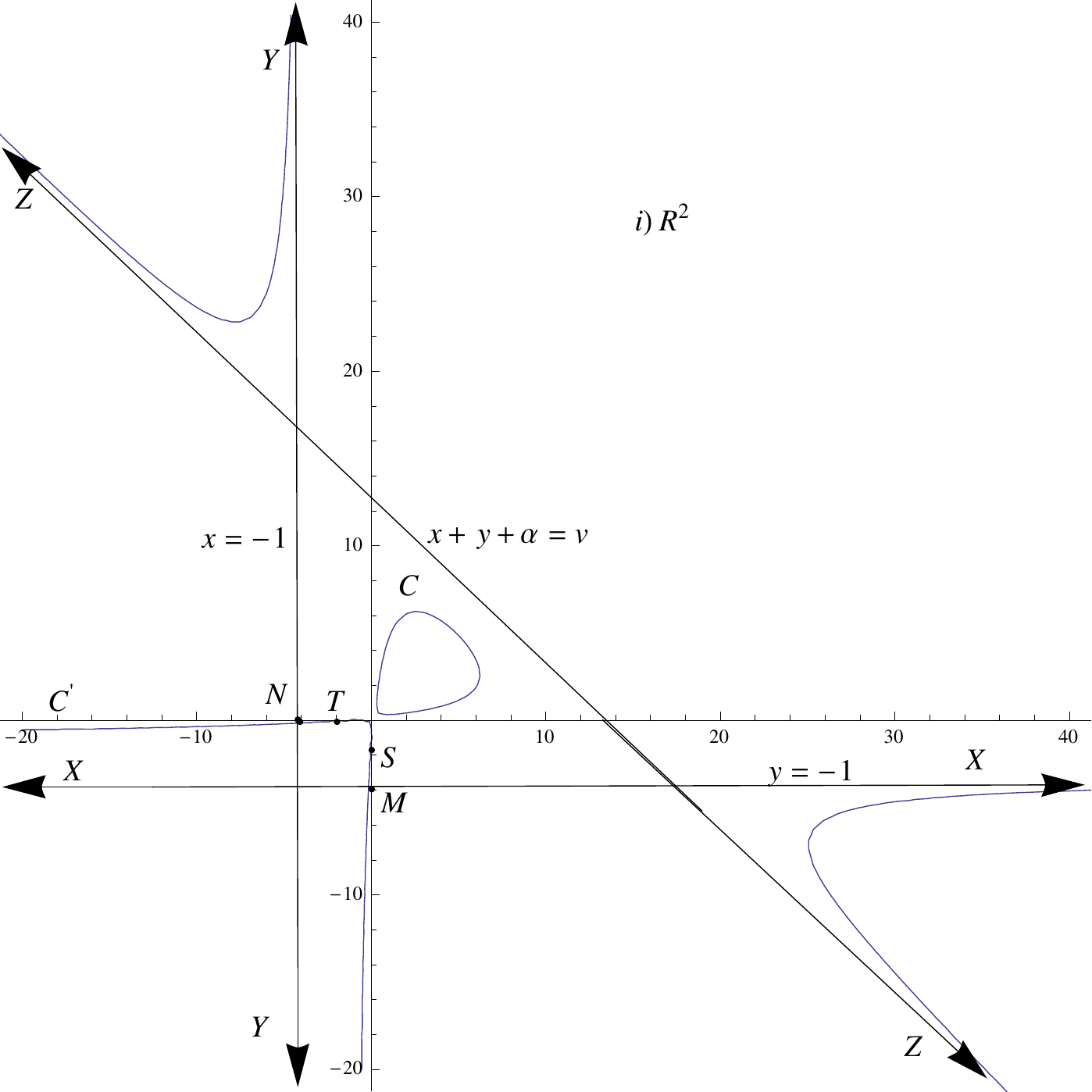}   \includegraphics[width=4cm,height=4cm]{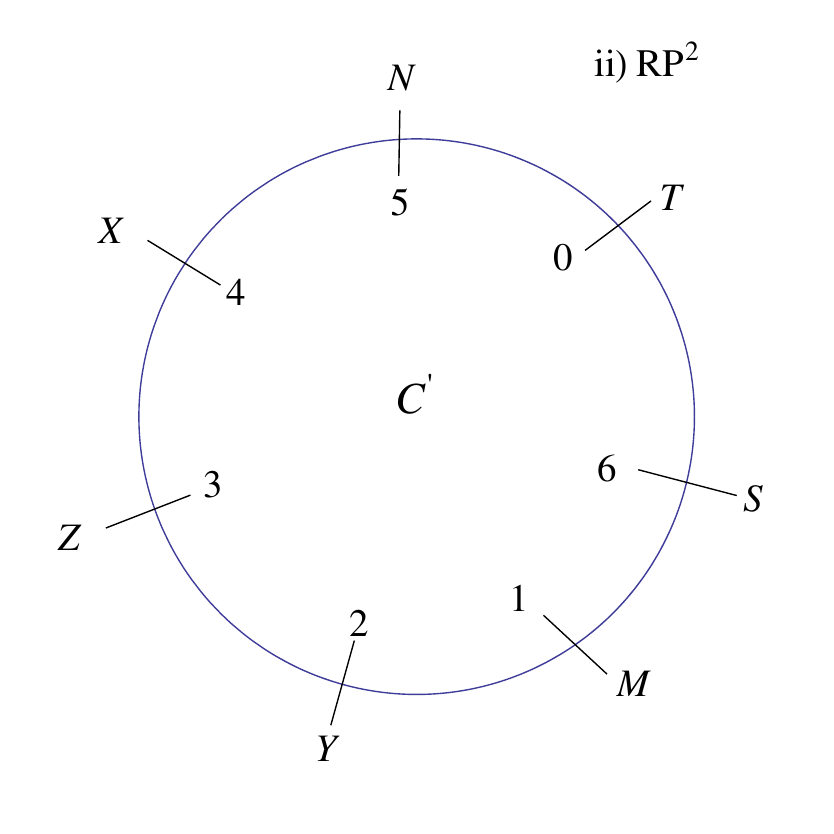}
\end{center}  
\caption[Curva $C$ en $\mathbb{R}^2$ cuando $0<\alpha<1$.]{ La \'orbita de $T$ sobre $C'$ en i) $\mathbb{R}^2$, y ii) $\mathbb{R}\mathbb{P}^2$, cuando $0<\alpha<1$. \label{sub:planoreal2}}
\end{figure}

En este caso $T$ est\'a  entre $N=\bar{f}^5(T)$ y $S=\bar{f}^6(T)$ sobre $C'$, entonces $5\rho<1<6\rho$. Por lo tanto, $\frac16<\rho<\frac15$.
\end{proof}

Ahora examinaremos el l\'imite de $\rho^v_\alpha$ cuando $v$ se aproxima a su valor m\'inimo en su dominio.

\begin{figure}[h!]
\begin{center}
\includegraphics[scale=.7]{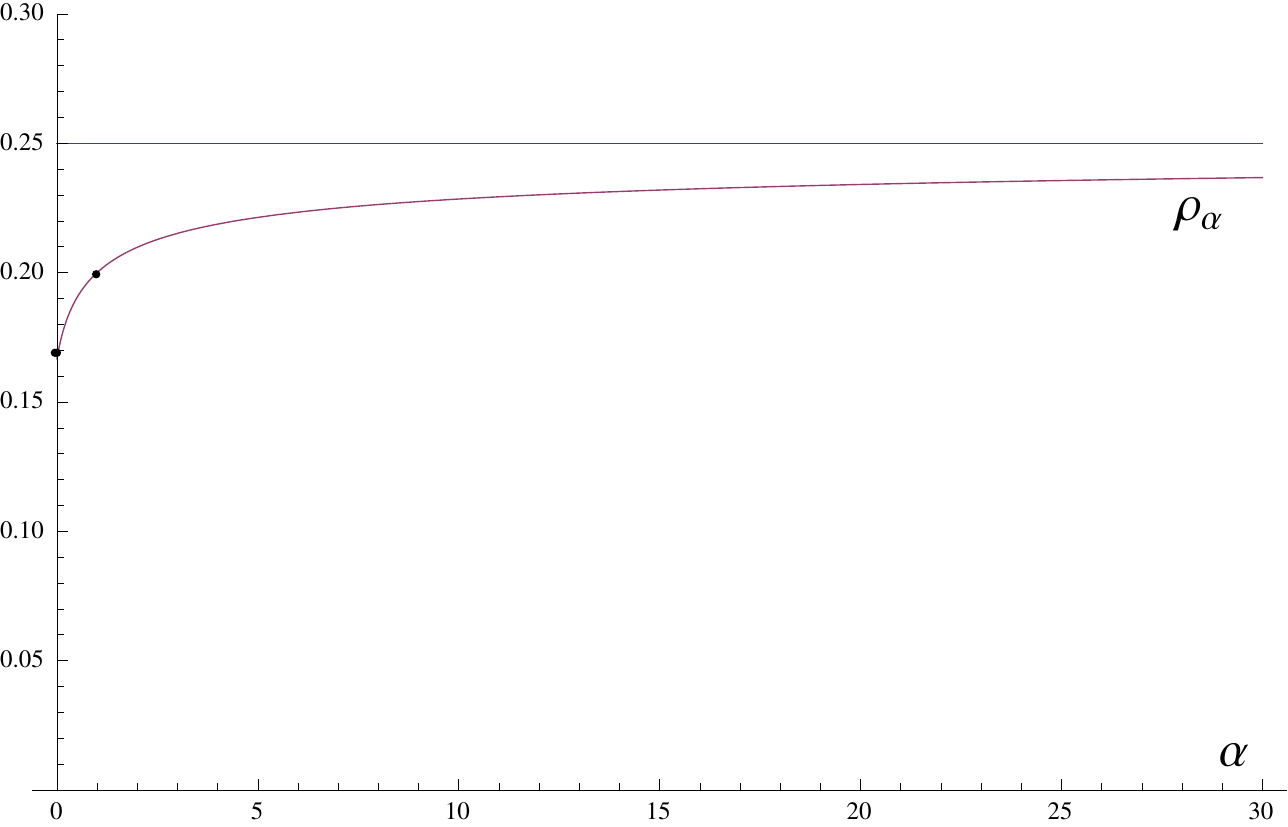}
\end{center}
\caption[Funci\'on $\rho_\alpha$]{Gr\'afica de $\rho_\alpha$. \label{sub:funcionro}}
\end{figure}

\begin{teo}
\label{sub:valorro}
Si $0\leq \alpha \leq \infty$  entonces $\rho^v_\alpha \rightarrow \rho_\alpha$ cuando $v \rightarrow v_\alpha$, y (ver figura~\ref{sub:funcionro}) $$\rho_\alpha=\frac{1}{2\pi}\arccos \left({\frac{1}{1+\sqrt{1+4\alpha}}}\right).$$
\end{teo}

\begin{proof}
Sea $f=f_\alpha$, $\omega=\omega_\alpha$, $F=F_\alpha=(\omega_\alpha,\omega_\alpha)$. Denotemos con $f^*$ la aproximaci\'on lineal de $f$ en $F$, es decir, $f^*(X)=f(F)+Df(F)(X-F)$, donde $X=(x,y)$.  Calculemos el n\'umero de rotaci\'on $\rho_\alpha$ de $f^*$. Sea $\xi,\eta$ coordenadas locales de $F=(\omega,\omega)$ ($\xi=x-\omega,\eta=y-\omega$) entonces 

$$f^*\begin{pmatrix} \xi \\ \eta  \\ \end{pmatrix}=f(F)+\begin{pmatrix}0 & 1\\ -1&1/\omega  \\ \end{pmatrix} \begin{pmatrix} \xi \\ \eta  \\ \end{pmatrix}$$

\noindent y haciendo una traslaci\'on  por medio de $\bar{\xi}=\xi-1+\omega$ y $\bar{\eta}=\eta-\omega$ tenemos la funci\'on lineal

$$f^*\begin{pmatrix} \bar{\xi} \\ \bar{\eta}  \\ \end{pmatrix}  = \begin{pmatrix}0 & 1\\ -1&1/\omega  \\ \end{pmatrix} \begin{pmatrix} \bar{\xi} \\ \bar{\eta}  \\ \end{pmatrix}.$$

Los valores propios $\lambda$ de la matriz $A$ asociada a $f^*$ son obtenidos por $\lambda(\lambda-1/\omega)+1=0$. Por lo tanto, $\lambda=e^{\pm i\theta}$, donde $\cos \theta=\frac{1}{2\omega}$. Calculamos los vectores propios de $A$ y obtenemos una matriz 

$$P=\begin{pmatrix}\cos \theta & -\sin \theta\\ 1&0  \\ \end{pmatrix}$$

\noindent tal que $P^{-1}AP$ es la matriz de rotaci\'on por un \'angulo $\theta$. Asi, $f^*$ es conjugada por un cambio lineal de coordenadas a una rotaci\'on del plano a trav\'es de un \'angulo $\theta$. Por lo tanto, el n\'umero de rotaci\'on $\rho_\alpha$ de $f^*$  es dado por 
$$\rho_\alpha=\rho(f^*)=\frac{\theta}{2\pi}=\frac{1}{2\pi}\arccos\left(\frac{1}{2\omega}\right)=\frac{1}{2\pi}\arccos \left(\frac{1}{1+\sqrt{1+4\alpha}}\right).$$

Demostremos ahora que $\rho^v_\alpha \rightarrow \rho_\alpha$ cuando $v \rightarrow v_\alpha$. Sea $S$ el c\'irculo de rayos a trav\'es de $F$ (aqu\'i un rayo es una semirecta en el plano iniciando en $F$). La aproximaci\'on lineal  $f^*$ induce\footnote{Para un rayo $s \in S$, $f^*(s)$ no es m\'as que la tangente a la curva $f(s)$ en $F$, es decir, $g$ se define como $g(s)=f^*(s)$ y por tanto tendr\'an el mismo n\'umero de rotaci\'on.} un difeomorfismo $g: S\rightarrow S$ tal que $g(s)$ es la tangente en $F$ de la curva diferenciable $f(s)$ en $\mathbb{R}^2_+$ para cada $s \in S$ (ver figura~\ref{sub:rayos}). Por tanto, esta funci\'on $g$ es conjugada  a la rotaci\'on  por el \'angulo $\theta$, entonces $\rho(g)=\frac{\theta}{2\pi}=\rho_\alpha$.

\begin{figure}[h!]
\begin{center}
\includegraphics[width=8cm,height=8cm]{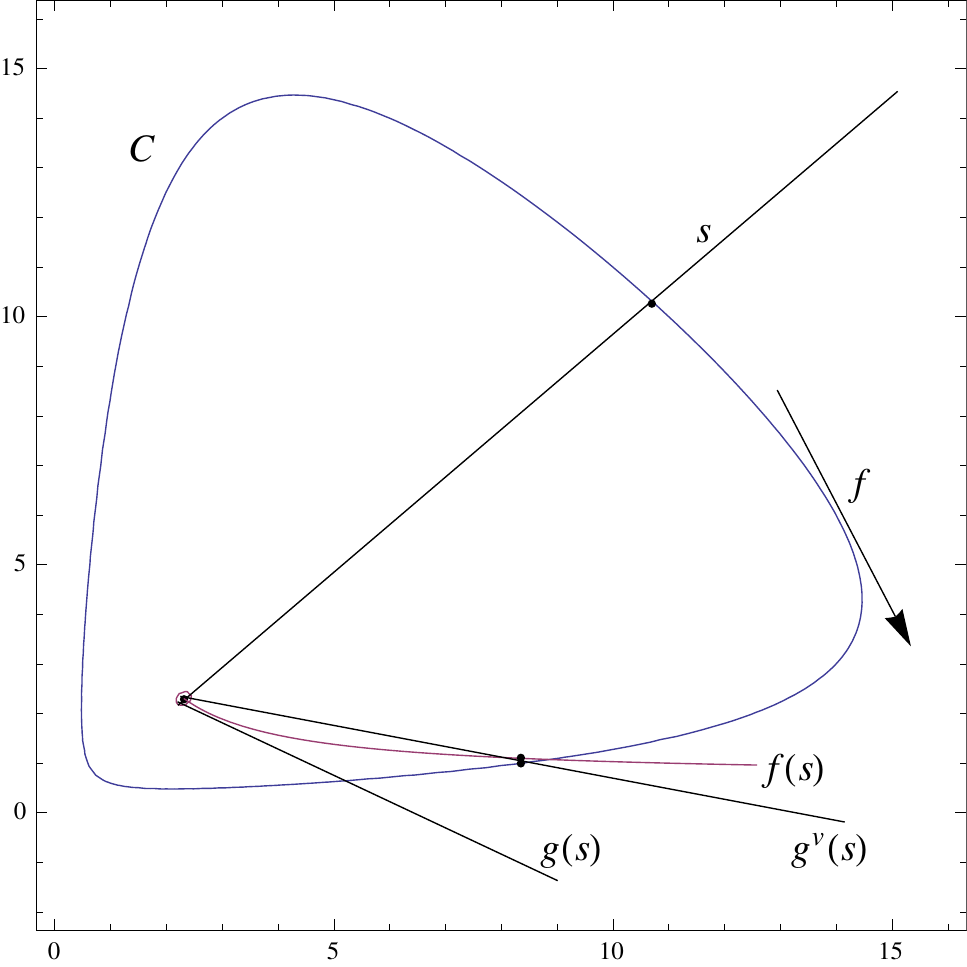}
\end{center}  
\caption[C\'irculo de rayos.]{Difeomorfismo del c\'irculo de rayos a trav\'es de $F$. \label{sub:rayos}}
\end{figure}

Sea $g^v:S\rightarrow S$ el rayo que pasa por el punto $f(s)\cap C$. Para cada $s \in S$, $g^v(s)$ est\'a bien definida por ser  $C$  una curva convexa. Entonces la proyecci\'on 
$$\pi:S \rightarrow C: s \mapsto f(s)\cap C $$

\noindent induce una conjugaci\'on entre $g^v$ y $f|_C$, ya que el siguiente diagrama claramente conmuta

\begin{center}
\begin{tabular}{ccc} 
$S$ & $\stackrel{g^v}{\rightarrow}$ & $S$ \\ 
 $\,\,\,\,\downarrow \pi$ &  & $\,\,\,\,\downarrow \pi$ \\ 
$C$ & $\stackrel{f}{\rightarrow}$ & $C.$ 
\end{tabular}
\end{center}

Por lo tanto, $\rho(g^v)=\rho(f|_C)=\rho^v_\alpha$. Como $v \rightarrow v_\alpha$ entonces el rayo $g^v(s)$ tiende a la tangente $g(s)$, por la diferenciabilidad de $f(s)$. As\'i, $g^v \rightarrow g$ cuando  $v \rightarrow v_\alpha$ por compacidad. Luego, $\rho(g^v) \rightarrow \rho(g)$ por la continuidad del n\'umero de rotaci\'on. Por lo tanto, $\rho^v_\alpha \rightarrow \rho_\alpha$.
\end{proof}

\begin{teo}\label{sub:aproximacionro}
Si $0<\alpha <\infty $ entonces $\rho^v_\alpha \rightarrow 1/5$ cuando $v\rightarrow \infty$. Si $v$ es grande entonces una aproximaci\'on para el n\'umero de rotaci\'on es dada por $$\rho \approx \frac{\ln{v}}{5\ln{v}-\ln{\alpha}}.$$
\end{teo}
\begin{proof}
Cuando $\alpha=1$ el resultado es inmediato porque $\rho^v_1=1/5$ para todo $v$, por el Teorema~\ref{sub:teoremadeparametros}. Supongamos primero que $1<\alpha$, entonces por el Teorema~\ref{sub:cotasro} $1/5<\rho^v_\alpha<1/4$. Por lo tanto, necesitamos mostrar que, dado un $\epsilon>0$, $$\exists v_0,\forall v, v>v_0\Rightarrow \rho^v_\alpha<1/5+\epsilon.$$

Sea $C=C^v_\alpha$ y $P=(\lambda,\lambda)$ un punto de la diagonal sobre $C$, tomando el que est\'a m\'as lejos del origen, entonces $\lambda\rightarrow\infty$ cuando $v\rightarrow\infty$. As\'i, es suficiente demostrar que $$\exists\lambda_0,\forall\lambda,\lambda>\lambda_0\Rightarrow \rho^v_\alpha<1/5+\epsilon.$$

Elijamos un entero $Q>\frac{1}{25\epsilon}$, y sea $N=5Q+6$. Sea $\{x_n\}_{n\geq0}$ la sucesi\'on determinada por la proyecci\'on de la \'orbita de $P$. Ahora introducimos otra sucesi\'on $\{y_n\}_{n\geq0}$ que comience igual y que est\'e pr\'oxima a $\{x_n\}_{n\geq0}$ para los primeros $N$ t\'erminos, pero que sea m\'as f\'acil de analizar. Definamos $y_0=1,y_1=y_2=\lambda$, y para $n\geq2$:
$$ 
y_{n+1} = \begin{cases} 
         \frac{\alpha+y_n}{y_{n-1}},& si\,\, n\equiv 0,3 \,\,(mod\, 5)\\ 
         \frac{y_n}{y_{n-1}},& si\,\, n\equiv1,2 \,\,(mod\, 5) \\
         \frac{\alpha}{y_{n-1}},& si\,\, n\equiv4 \,\,(mod\, 5).
       \end{cases} 
$$ 

Con $\{y_n\}_{n\geq0}$ as\'i definida  verifiquemos por inducci\'on que las f\'ormulas siguientes son ciertas para $q\geq0$:

\begin{align} \label{sub:induccion1}
    y_{5q}   & = \alpha^q     \\ 
    y_{5q+1} & =\frac{2\alpha^q\lambda}{1+\alpha^q} \\
    y_{5q+2} & =\frac{2\lambda}{1+\alpha^q} \\
    y_{5q+3} & =\frac{1}{\alpha^q} \\
    y_{5q+4} & =\frac{(1+\alpha^{q+1})(1+\alpha^q)}{2\alpha^q\lambda} 
\end{align}

Para $q=0,1$ es f\'acil ver que se cumplen. Supongamos que se cumplen para todo  $p\leq q$. Ahora, probemos para $q+1$.
\begin{align*} 
    y_{5(q+1)}   & = \frac{\alpha}{y_{5q+3}}=\frac{\alpha}{\frac{1}{\alpha^q}}=\alpha^{q+1}     \\ 
    y_{5(q+1)+1} & =\frac{\alpha+y_{5(q+1)}}{y_{5q+4}}=\frac{\alpha+\alpha^{q+1}}{\frac{(1+\alpha^{q+1})(1+\alpha^q)}{2\alpha^q\lambda}}=\frac{2\alpha^{q+1}\lambda}{1+\alpha^{q+1}} \\
    y_{5(q+1)+2} & =\frac{y_{5(q+1)+1}}{y_{5(q+1)}}=\frac{\frac{2\alpha^{q+1}\lambda}{1+\alpha^{q+1}}}{\alpha^{q+1}}=\frac{2\lambda}{1+\alpha^{q+1}} \\
    y_{5(q+1)+3} & =\frac{y_{5(q+1)+2}}{y_{5(q+1)+1}}=\frac{\frac{2\lambda}{1+\alpha^{q+1}}}{\frac{2\lambda\alpha^{q+1}}{1+\alpha^{q+1}}}=\frac{1}{\alpha^{q+1}} \\
    y_{5(q+1)+4} & =\frac{\alpha+y_{5(q+1)+3}}{y_{5(q+1)+2}}=\frac{\alpha+\frac{1}{\alpha^{q+1}}}{\frac{2\lambda}{1+\alpha^{q+1}}}=\frac{(1+\alpha^{q+2})(1+\alpha^{q+1})}{2\alpha^{q+1}\lambda} 
\end{align*}
El razonamiento anterior para la construci\'on de $\{y_n\}_{n\geq0}$ es como sigue. Suponiendo que $\lambda\gg1$, hemos simplificado el numerador $\alpha+y_n$ en la difinici\'on de $y_{n+1}$ ignorando a $\alpha$ si $\alpha\ll y_n$ e ignorando a $y_n$ si $y_n\ll\alpha$. Estos cambios son relativamente peque\~nos para los primeros $N$ t\'erminos de $\{y_n\}_{n\geq0}$ que es siempre pr\'oxima a $\{x_n\}_{n\geq0}$, cuando $\lambda$ es suficientemente grande. M\'as precisamente, podemos mostrar que, dado $\delta>1$,
\begin{equation}\label{sub:asterisco}
\exists\lambda_1,\forall\lambda,\forall n,\,\,\lambda>\lambda_1 \& 1\leq n\leq N \Rightarrow \delta^{-1}<\frac{x_n}{y_n}<\delta.
\end{equation}
Para $0\leq q\leq Q$ sea 

\begin{align*}
X_q&=(x_{5q+5},x_{5q+6})=f^{5q+4}(P)\\
Y_q&=(y_{5q+5},y_{5q+6})=\left(\alpha^{q+1},\frac{2\alpha^{q+1}\lambda}{1+\alpha^{q+1}}\right).
\end{align*}
Sea $\xi_q,\eta_q$ las pendientes de los vectores $0X_q$, $0Y_q$. Entonces $$\eta_q=\frac{2\lambda}{1+\alpha^{q+1}}$$ as\'i $\eta_0>\eta_1> \dots >\eta_Q$, por ser $\alpha>1$. Tomando $\lambda>\alpha^{Q+1}$ tenemos $$\eta_Q>\frac{2\alpha^{Q+1}}{1+\alpha^{Q+1}}>1.$$

Por lo tanto, $Y_0,Y_1,\dots,Y_Q$ est\'an todos por encima de la diagonal. Utilizamos ahora la ecuaci\'on (~\ref{sub:asterisco}) para mostrar que esto mismo es cierto para $X_n,0\leq n \leq Q$ tomando a $\lambda$ lo suficientemente grande, como sigue. Sea $\delta=\sqrt{\eta_Q}$, por ser $\delta>1$ existe $\lambda_1$ satisfaciendo la ecuaci\'on (~\ref{sub:asterisco}). Sea $\lambda_0=\max\{\alpha^{Q+1},\lambda_1\}$. As\'i, por la ecuaci\'on (~\ref{sub:asterisco}), si $\lambda>\lambda_0$ y para $0\leq q \leq Q$ tenemos $$\xi_q=\frac{x_{5q+5}}{x_{5q+6}}>\frac{\delta^{-1}y_{5q+6}}{\delta y_{5q+5}}=\frac{\eta_q}{\delta^2}\geq\frac{\eta_Q}{\delta^2}=1.$$
Por lo tanto, $X_0,X_1,\dots,X_Q$ est\'an todos por encima de la diagonal (ver figura~\ref{sub:xarriba}).

\begin{figure}[h!]
\begin{center}
\includegraphics[width=8cm,height=8cm]{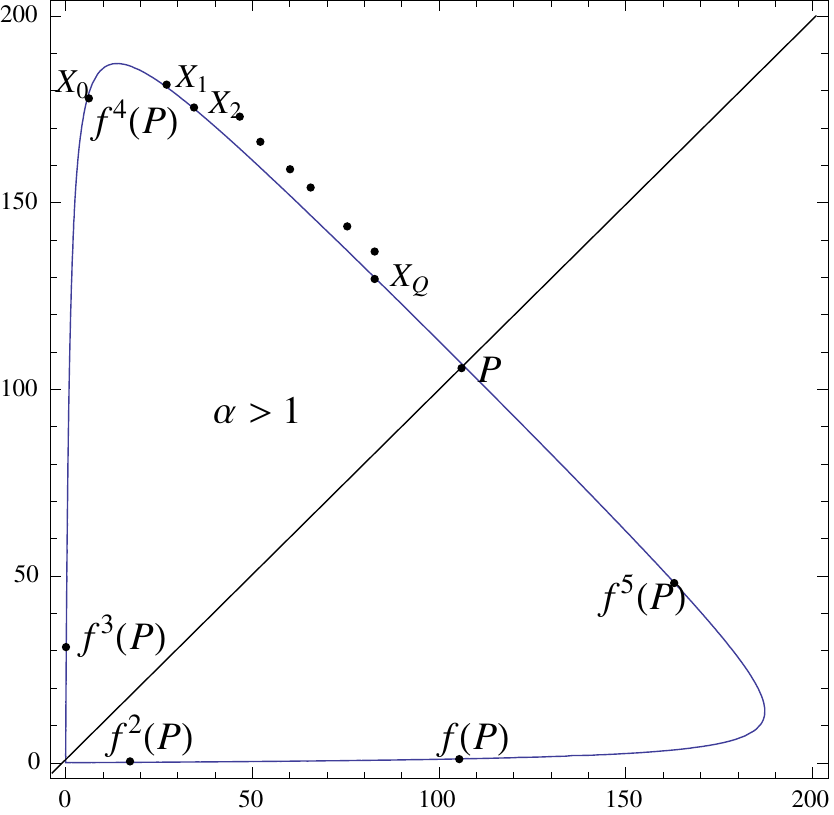}
\end{center}
\caption[Curva cuando $v\rightarrow\infty$ con  $\alpha>1$.]{Los $X_n,0\leq n \leq Q$ est\'an por encima de la diagonal cuando $\alpha>1$. \label{sub:xarriba}}
\end{figure}

La \'orbita de $P$ contiene a todos los $X_i$, y hace cinco iteraciones  en todo $C$ entre dos  $X_i$ consecutivos. Cuando la \'orbita alcanza $X_Q=f^{5Q+4}(P)$ tienen que haberse hecho ligeramente menos de $Q+1$ rotaciones alrededor de $C$ por que $X_Q$ est\'a arriba de la diagonal y hemos comenzado en $P$ que est\'a sobre la diagonal. Por lo tanto, $(5Q+4)\rho^v_\alpha<Q+1$. As\'i, 
\begin{align*} 
    \rho^v_\alpha & < \frac{Q+1}{5Q+4}    \\ 
                  & =\frac15\left(\frac{1+\frac1Q}{1+\frac{4}{5Q}}\right) \\
                  & <\frac15\left(1+\frac1Q-\frac{4}{5Q}\right), \tag*{porque $a>b>0 \Rightarrow \frac{1+a}{1+b}<1+a-b$} \\
                  & =\frac15\left(1+\frac{1}{5Q}\right)\\
                  & =\frac15+\frac{1}{25Q}\\
                  & <\frac15+\epsilon,\tag*{porque $Q>\frac{1}{25\epsilon}$} 
\end{align*}
Por lo tanto, $\rho^v_\alpha \rightarrow \frac15$ cuando $v\rightarrow\infty$, como se quer\'ia.\\

Para el otro caso $\alpha<1$ se utiliza un argumento similar. En este caso $1/6<\rho^v_\alpha<1/5$, por el Teorema~\ref{sub:cotasro}. Ahora necesitamos demostrar que dado $\epsilon>0$, $$\exists\lambda_0,\forall\lambda,\lambda>\lambda_0\Rightarrow \rho^v_\alpha>1/5-\epsilon.$$
Sea $Q,N$ como antes. Ahora, sea $X_q=(x_{5q+2},x_{5q+3})=f^{5q+1}(P)$, para $0\leq q \leq Q$, y  $Y_q=(y_{5q+2},y_{5q+3})$ su correspondiente t\'ermino en $\{y_n\}_{n\geq0}$. Sea $\xi_q,\eta_q$ como antes, entonces $$\eta_q=\frac{\frac{1}{\alpha^q}}{\frac{2\lambda}{1+\alpha^q}}=\frac{1+\alpha^q}{2\lambda\alpha^q}=\frac{\frac{1}{\alpha^q}+1}{2\lambda}$$  por tanto $\eta_Q>\eta_{Q-1}> \dots >\eta_1>\eta_0$ por ser $\alpha<1$. Tomando a $\lambda>\alpha^{-Q}$ se tiene que

 $$\eta_Q=\frac{1+\alpha^q}{2\lambda\alpha^q}<\frac{\alpha^Q(1+\alpha^Q)}{2\alpha^Q}=\frac{1+\alpha^Q}{2}<1.$$ 
 
As\'i, todos los $Y_i$ est\'an abajo de la diagonal. Por ecuaci\'on (~\ref{sub:asterisco}) esto mismo es cierto para los $X_i$ (ver figura~\ref{sub:xabajo}), esto se ve como sigue. Sea $\delta=\sqrt{\frac1\eta_Q}$, siendo $\delta>1$, existe $\lambda_1$ tal que la ecuaci\'on (~\ref{sub:asterisco}) se satisface. Sea $\lambda_0=\max\{\lambda_1,\alpha^{-Q}\}$, entonces por la misma ecuaci\'on, si $\lambda>\lambda_0$ se tiene que para $0\leq q\leq Q$ $$\xi_q=\frac{x_{5q+3}}{x_{5q+2}}<\frac{\delta(y_{5q+3})}{\delta^{-1}(x_{5q+2})}=\delta^2\eta_q<\delta^2\eta_Q=1.$$

\begin{figure}[h!]
\begin{center}
\includegraphics[width=8cm,height=8cm]{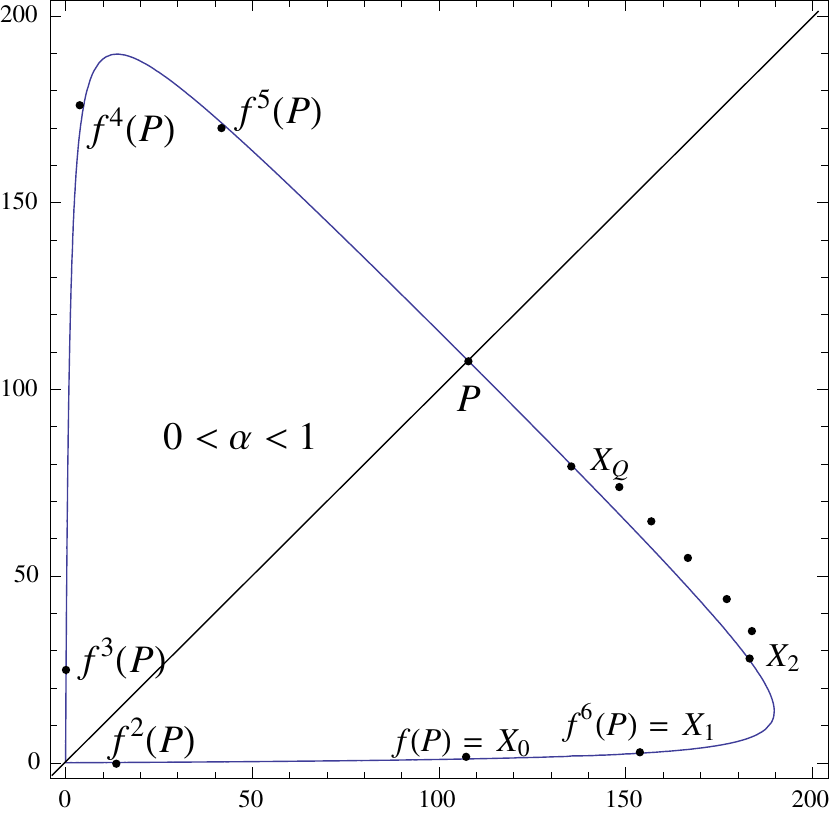}
\end{center}
\caption[Curva cuando $v\rightarrow\infty$ con  $0<\alpha<1$.]{Los $X_n,0\leq n \leq Q$ est\'an abajo de la diagonal cuando $0<\alpha<1$. \label{sub:xabajo}}
\end{figure}

Por lo tanto, cuando la \'orbita alcanza a $X_Q=f^{5Q+1}(P)$ \'esta ha dado poco m\'as de $Q$ rotaciones. As\'i, $(5Q+1)\rho^v_\alpha>Q$, y tenemos $$\rho^v_\alpha>\frac{Q}{5Q+1}=\frac15\left(\frac{1}{1+\frac{1}{5Q}}\right)>\frac15\left(1-\frac{1}{5Q}\right)=\frac15-\frac{1}{25Q}>\frac15-\epsilon,$$
como deseabamos. Por tanto, $\rho^v_\alpha\rightarrow1/5$ cuando $v\rightarrow\infty$.\\

Para la segunda parte del teorema, tomemos $\rho=\rho^v_\alpha$ y veamos primero que si $\alpha=1$ obtenemos por medio de la aproximaci\'on el valor correcto $\rho^v_\alpha=1/5$ para todo $v$. Supongamos ahora que $1<\alpha$, cuando $v$ es grande entonces $C$ tiene  aproximadamente la forma de un tri\'angulo con hipotenusa $x+y=v$ (ver figura~\ref{sub:pendientes}), en lo que sigue $v$ lo tomaremos siempre como un valor grande tanto como queramos.\\

\begin{figure}[h!]
\begin{center}
\includegraphics[width=8cm,height=8cm]{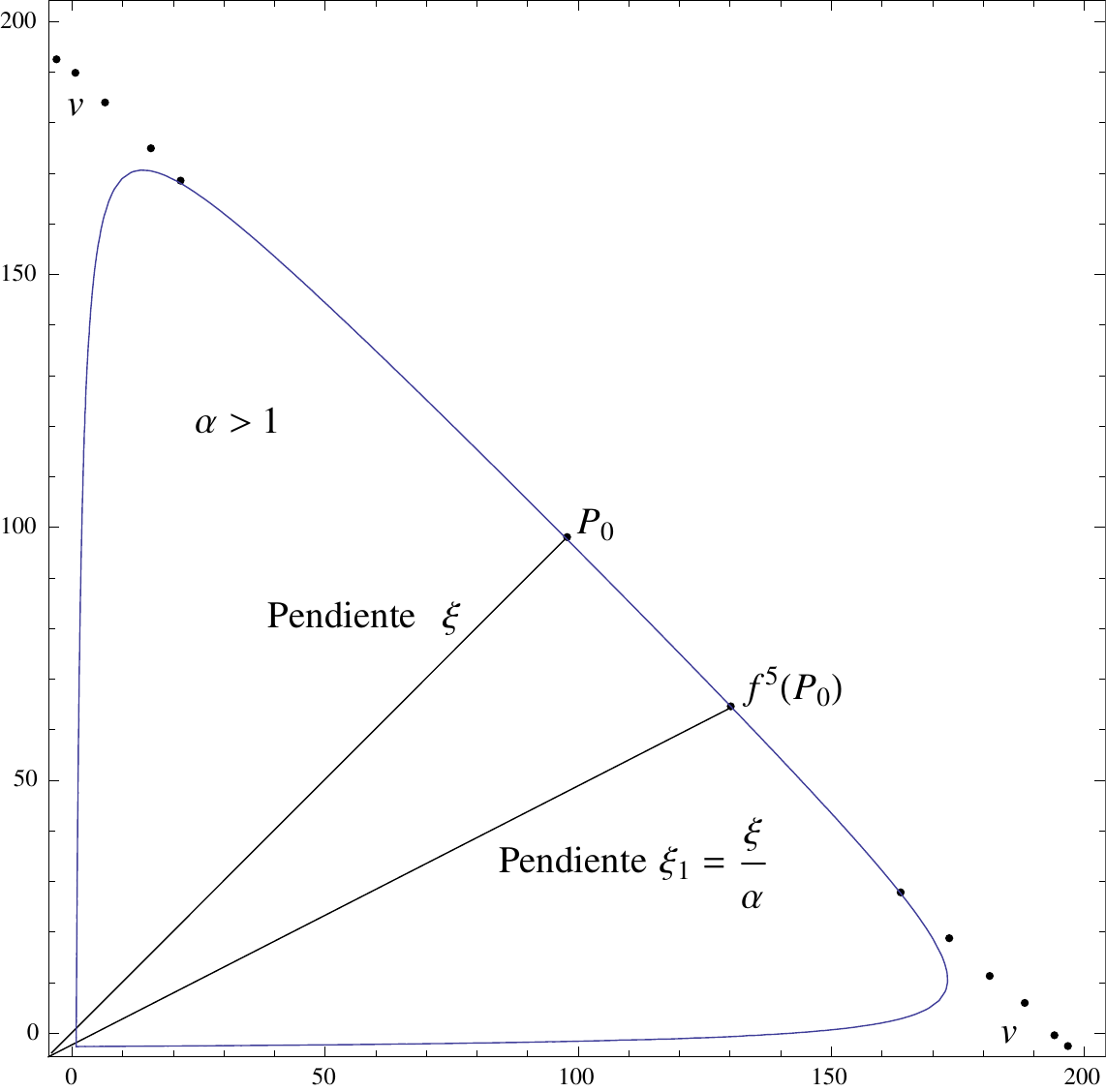}
\end{center}
\caption[Cuando $v$ es grande, $C$ tiene forma de tri\'angulo.]{Si $\alpha>1$ las pendientes decrecen por el factor $\frac1\alpha$. \label{sub:pendientes}}
\end{figure}

Sea $\{y_n\}_{n\geq0}$ la sucesi\'on iniciando en $P_0=(y_1,y_2)=(tv,tv)$, donde $s+t=1$ y $s,t>0$. En el calculo de $y_{n+1}$ para $n\geq2$ hacemos la misma aproximaci\'on como antes, ignorando a $\alpha$ si $\alpha \ll y_n$ e ignorando a $y_n$ si $y_n \ll \alpha$. Ahora la sucesi\'on  es $$y_3=\frac{s}{t},y_4=\frac{s+\alpha t}{stv},y_5=\frac{\alpha t}{s},y_6=\frac{\alpha tv}{s+\alpha t},y_7=\frac{sv}{s+\alpha t}, \dots$$
haciendo el mismo proceso de inducci\'on que en las ecuaciones (~\ref{sub:induccion1}) y siguientes, podemos probar f\'acilmente que para $q\geq1$ se cumple:
\begin{align*} 
    y_{5q}   & = \frac{\alpha^qt}{s}     \\ 
    y_{5q+1} & =\frac{\alpha^qtv}{s+\alpha^qt} \\
    y_{5q+2} & =\frac{sv}{s+\alpha^qt} \\
    y_{5q+3} & =\frac{s}{\alpha^qt} \\
    y_{5q+4} & =\frac{(s+\alpha^{q+1}t)(s+\alpha^qt)}{\alpha^qtvs} 
\end{align*}

\noindent por lo tanto, $$f^{5q}(P_0)\approx(y_{5q+1},y_{5q+2})=\left(\frac{\alpha^q tv}{s+\alpha^q t},\frac{sv}{s+\alpha^q t}\right).$$

Sea 
\begin{align*} 
    \xi \,\,  & = \,la\,\, pendiente\,\,del \,\,vector\,\, 0P_0=\frac{s}{t},    \\ 
    \xi_1 & = \,la\,\, pendiente\,\, del\,\, vector\,\, 0f^5(P_0)\approx\frac{s}{\alpha t},    \\
   \vdots &   \\
    \xi_n & = \,la\,\, pendiente\,\, del\,\, vector\,\, 0f^{5n}(P_0)\approx\frac{y_{5q+2}}{y_{5q+1}}=\frac{s}{\alpha^nt}.               
\end{align*}

As\'i, $\xi_1\approx\frac{\xi}{\alpha}$ y de forma semejante $\xi_n\approx\frac{\xi}{\alpha^n}$, por ser $\alpha>1$ la sucesi\'on $\xi_n$ es decreciente. Hemos aproximado las pendientes para $f^{5n}(P_0)$ por medio de los $y_n$, ya que por la ecuaci\'on (~\ref{sub:asterisco}) la sucesi\'on $\{y_n\}_{n>0}$ est\'a siempre pr\'oxima a la \'orbita de $P_0$ y por tanto sus pendientes. \\

Sea $P=(\sqrt{v},v)$, entonces $f(P)\approx(v,\sqrt{v})$ ($f(P)=(v,\frac{\alpha+v}{\sqrt{v}})\approx(v,\frac{v}{\sqrt{v}})$ ignorando a $\alpha$ porque $v$ es grande), es importante subrayar que las aproximaciones de arriba se hicieron para $P_0=(tv,sv)$ con $s+t=1$  y ahora tomamos a $P=(\sqrt{v},v)$, en este caso $s=1$ y $t=\frac{1}{\sqrt{v}}$, pero como $v$ es tan grande como queramos entonces $t\approx0$ y por tanto $s+t\approx1$.\\

\begin{figure}[h!]
\begin{center}
\includegraphics[width=8cm,height=8cm]{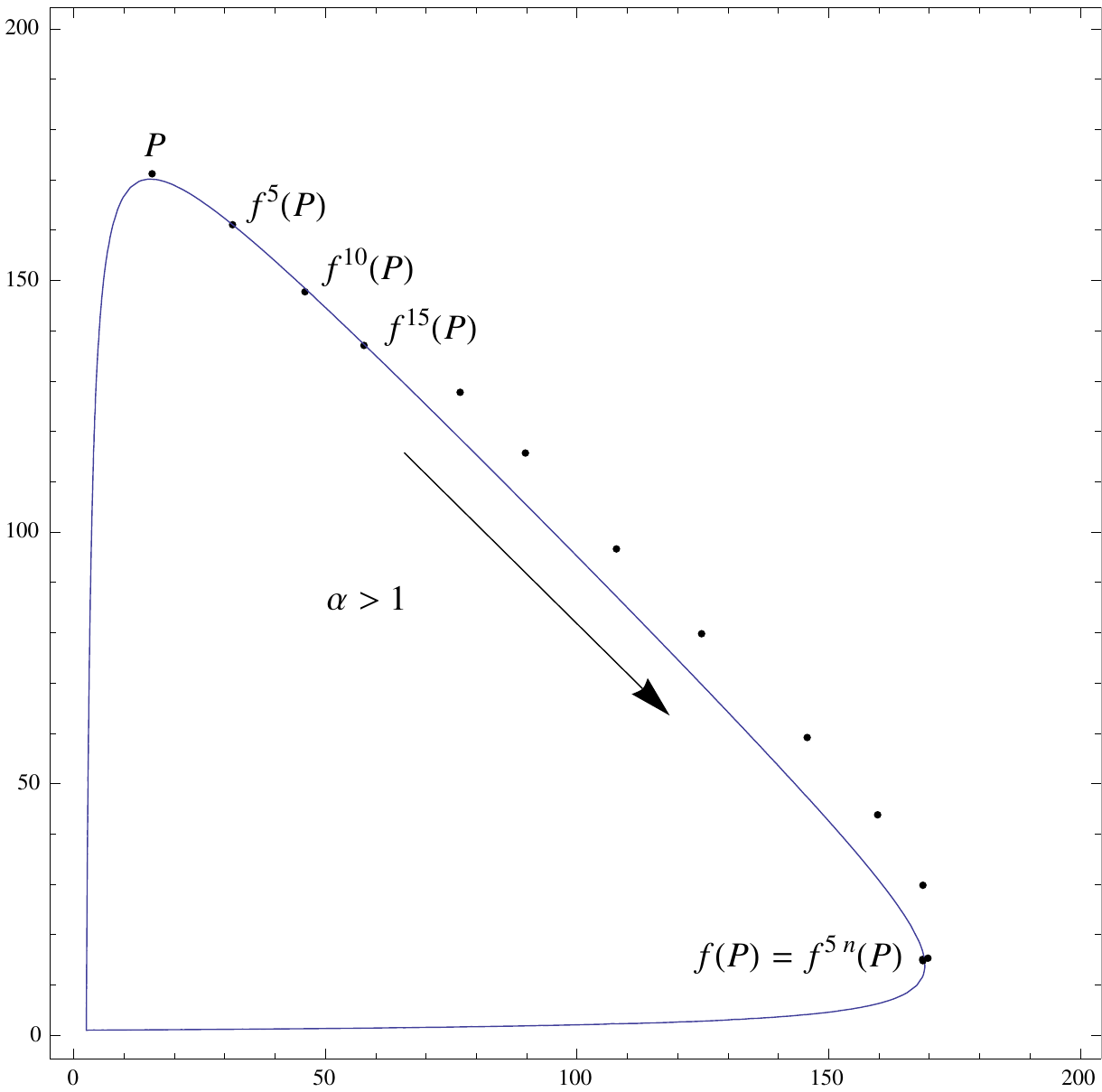}
\end{center}
\caption[Sucesi\'on de $f^{5n}(P)$ sobre $C$.]{Los puntos $f^{5n}(P)$ van desplaz\'andose lentamente sobre la ``hipotenusa del tri\'angulo'' desde $P$ hasta $f(P)$. \label{sub:pendientes2}}
\end{figure}
Las potencias $P,f^{5n}(P),f^{10n}(P),f^{15n}(P),\ldots$ van desplaz\'andose lentamente (dependiendo del valor $\alpha$) por la ``hipotenusa'' de $C$ desde $P$ hasta $f(P)$ (ver figura~\ref{sub:pendientes2}) por ser las pendientes siempre decreciente y porque la pendiente para $P$ es $\sqrt{v}$ que es un valor grande  y la pendiente de $f(P)$ es  aproximadamente $\frac{1}{\sqrt{v}}$ que es un valor muy peque\~no, as\'i, para alg\'un $n$ se tendr\'a que $f^{5n}(P)\approx f(P)$, entonces $\xi=\sqrt{v}$ y $\xi_n\approx\frac{1}{\sqrt{v}}$. Por lo tanto, $$\frac{1}{\sqrt{v}}\approx\xi_n=\frac{\xi}{\alpha^n}=\frac{\sqrt{v}}{\alpha^n}.$$ As\'i, $\alpha^n\approx v$, luego $n\ln \alpha\approx\ln v$. Y se tiene que $$n\approx\frac{\ln v}{\ln \alpha}.$$
Durante cada cinco iteraciones la \'orbita rota una sola vez sobre todo $C$ m\'as un poco. Despu\'es de $5n$ iteraciones la \'orbita alcanza a $f(P)$ aproximadamente, as\'i, $5n\rho\approx n+\rho$. Por lo tanto $$\rho\approx\frac{n}{5n-1}\approx\frac{\ln v}{5\ln v-\ln \alpha }.$$ 

Finalmente, consideremos el caso cuando $\alpha<1$. Sea $P=(\sqrt{v},v)$ y $f(P)\approx (v,\sqrt{v})$ como antes. Calculemos las pendientes de $f^{5n+1}(P_0)$ como las aproximaciones de antes.
\begin{align*} 
    \xi \,\,  & = \,la\,\, pendiente\,\,del \,\,vector\,\, 0f(P_0)\approx\frac{1}{tv},    \\ 
    \xi_1 & = \,la\,\, pendiente\,\, del\,\, vector\,\, 0f^{5+1}(P_0)\approx\frac{s}{\alpha t v}+\frac1v,    \\
   \vdots &   \\
    \xi_n & = \,la\,\, pendiente\,\, del\,\, vector\,\, 0f^{5n+1}(P_0)\approx\frac{y_{5q+3}}{y_{5q+2}}=\frac{s}{\alpha^ntv}+\frac1v.               
\end{align*}

Ahora, con la aproximacion para $P$ sabemos que $v$ es muy grande y por tanto despreciamos sumando $\frac1v$ y tambi\'en sabemos que $s=1,t=\frac{1}{\sqrt{v}}$, as\'i, nos queda que $\xi\approx\frac{1}{\sqrt{v}}$, $\xi_1\approx\frac{1}{\alpha\sqrt{v}}$ y $\xi_n\approx\frac{1}{\alpha^n\sqrt{v}}$. Ahora,  $\frac{\xi}{\alpha}>\xi$ y por lo tanto las pendientes son crecientes. As\'i, durante cada cinco iteraciones la \'orbita rota sobre todo $C$  una sola vez menos un poco (ver figura~\ref{sub:pendientes3}). Por tanto,  las potencias $f(P),f^6(P),f^{11}(P),\ldots$ van desplaz\'andose hacia atr\'as a lo largo de la ``hipotenusa'' desde $f(P)$ hasta $P$. Por el mismo razonamiento que antes, existir\'a un $n$ tal que $f^{5n+1}(P)\approx P$ y por tanto sus pendientes ser\'an aproximadamente iguales. As\'i $$\sqrt{v}\approx\xi_n\approx\frac{\xi}{\alpha^n}\approx\frac{1}{\alpha^n\sqrt{v}}.$$

\begin{figure}[h!]
\begin{center}
\includegraphics[width=8cm,height=8cm]{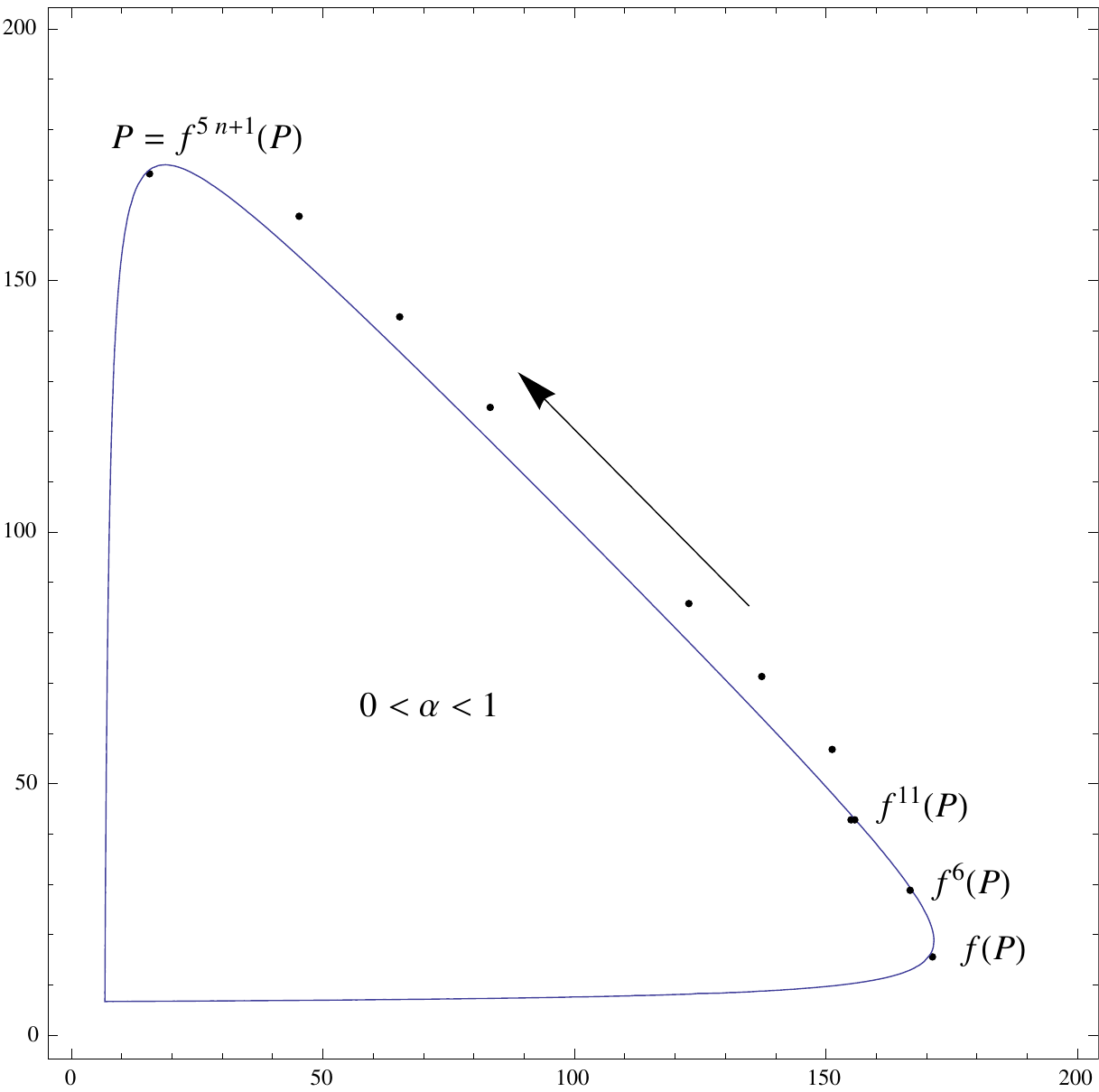}
\end{center}
\caption[Sucesi\'on de $f^{5n+1}(P)$ sobre $C$.]{Los puntos $f^{5n+1}(P)$ se desplazan desde $f(P)$ a $P$ sobre la ``hipotenusa del tri\'angulo''. \label{sub:pendientes3}}
\end{figure}

Por lo tanto $\alpha^n\approx\frac1v$, luego, $n \ln \alpha\approx -\ln v$, y entonces $n\approx -\frac{\ln v}{\ln \alpha}$. Ahora, despu\'es de $5n+1$ iteraciones la \'orbita alcanza aproximadamente a $P$, y as\'i,  $(5n+1)\rho\approx n$. Por lo tanto $$\rho\approx \frac{n}{5n+1}\approx\frac{-\frac{\ln v}{\ln \alpha}}{5(-\frac{\ln v}{\ln \alpha})+1}=\frac{\ln v}{5\ln v-\ln \alpha}.$$

\end{proof}

\begin{figure}[h!]
\begin{center}
\includegraphics[width=6cm,height=5cm]{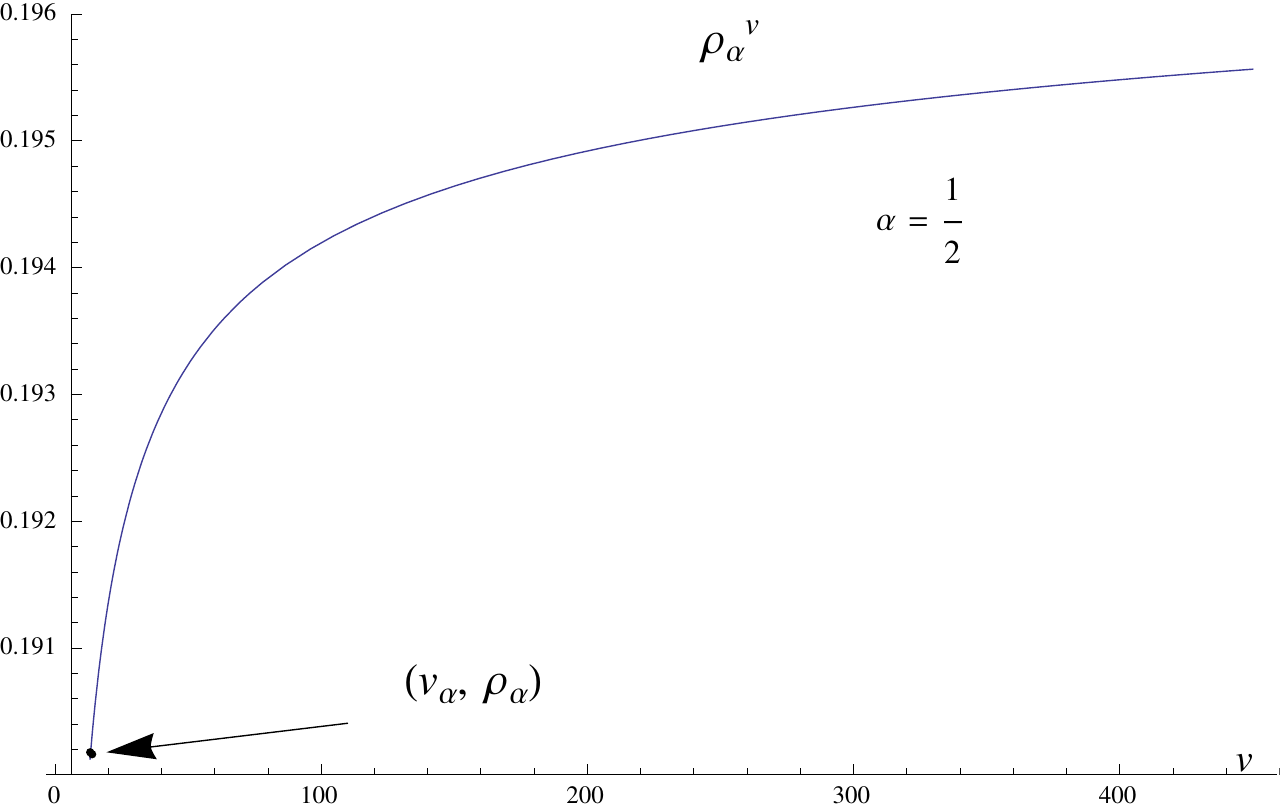} \quad \includegraphics[width=6cm,height=5cm]{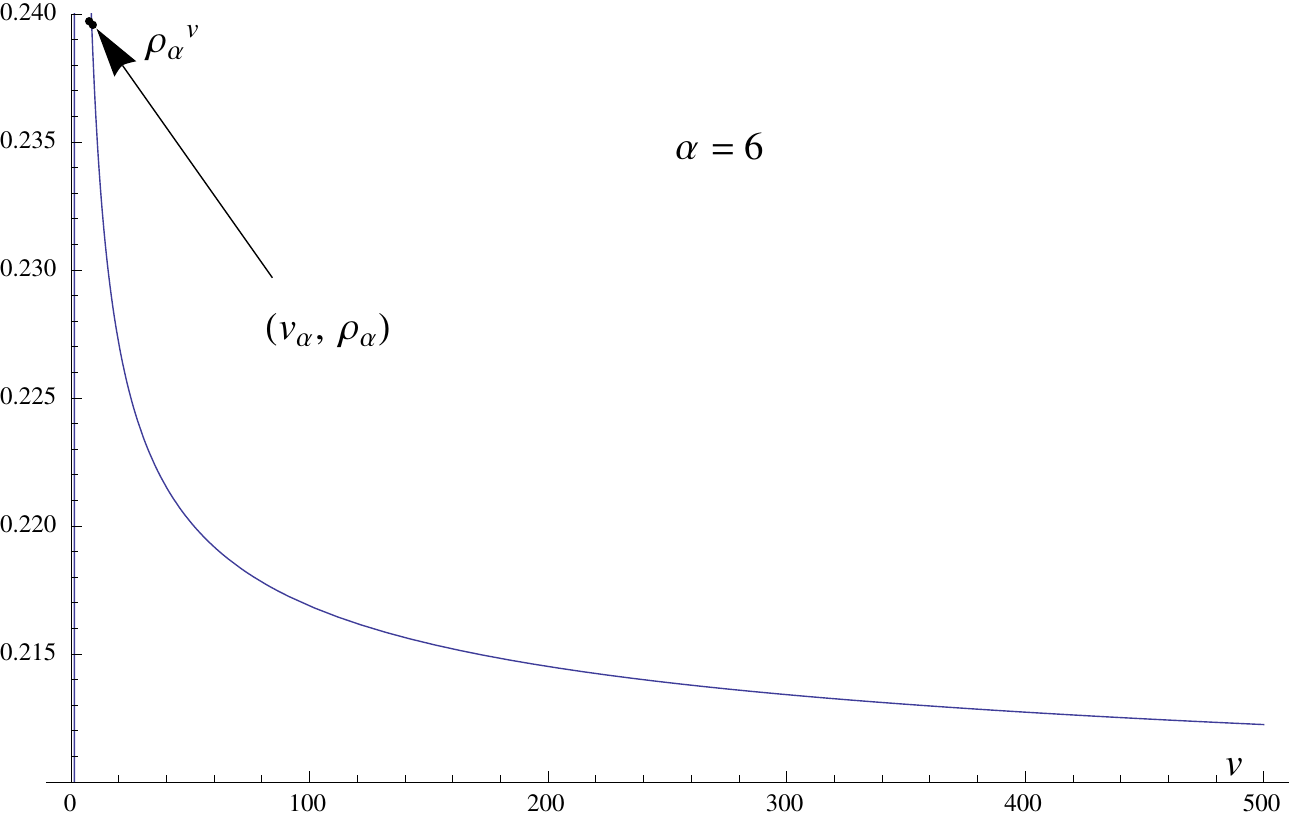}
\end{center}
\caption[Gr\'afica de $\rho^v_\alpha$ ]{Gr\'aficas de $\rho^v_\alpha$ como funci\'on de $v$ para $\alpha= \frac12,6$ con la aproximaci\'on del Teorema~\ref{sub:aproximacionro}.  \label{sub:rhodev}}
\end{figure}

Con  la aproximaci\'on del teorema anterior, E. C. Zeeman en 1996\cite{art}, represent\'o gr\'aficamente  la funci\'on  $\rho^v_\alpha$ para varios valores de $\alpha$ (ver figura~\ref{sub:rhodev}). Con \'esto f\'ormulo la siguiente conjetura

\begin{conje}
\label{sub:conjetura} La funci\'on $\rho^v_\alpha$ es diferenciable respecto a $v$ y se tiene:
\begin{itemize}
\item[a)] Si $0<\alpha<1$ entonces $\rho^v_\alpha$ es estrictamente creciente con relaci\'on a $v$.
\item[b)] Si $1<\alpha<\infty$ entonces $\rho^v_\alpha$  es estrictamente decreciente con relaci\'on a $v$.
\end{itemize}
\end{conje}

Esta conjetura fue probada al a\~no siguiente por F. Beukers y R. Cushman \cite{cje}. El argumento utilizado es como sigue: se considera un  sistema hamiltoniano correspondiente a la funci\'on hamiltoniana $V$. El flujo del campo vectorial hamiltoniano cuyas curvas integrales son las soluciones de 

\begin{align*}
\dot{x}&=xy\frac{\partial V}{\partial y}=(x+1)\left(y+\frac{x+\alpha}{y}\right)\\
\dot{x}&=-xy\frac{\partial V}{\partial x}=-(y+1)\left(x-\frac{y+\alpha}{x}\right)
\end{align*}

\noindent porque $V$ es una funci\'on con todos sus puntos cr\'iticos no degenerados y que es un m\'inimo correspondiente al valor $v_\alpha$. Sabemos que $V^{-1}(v)=C^v_\alpha$ es difeomorfa a la circunferencia unidad. As\'i, todas las \'orbitas del campo vectorial hamiltoniano con $v>v_\alpha$ son peri\'odicas de per\'iodo $T(v)$. Como probamos en el Teorema~\ref{sub:casirotacion}  $f|_{C^v_\alpha}$ es conjugada a una rotaci\'on de la circunferencia por medio de un \'angulo $\rho^v_\alpha$ que depende continuamente  de $v$. Luego, el tiempo $\tau(v)$  que  tarda la curva integral del campo vectorial hamiltoniano iniciando en $(x,y) \in C^v_\alpha$ en llegar a $f(x,y) \in C^v_\alpha$ no depende del punto inical $(x,y)$, solo depende del valor de $v$. Por lo tanto, el n\'umero de rotaci\'on de la funci\'on $f|_{C^v_\alpha}$ es $\tau(v)/T(v)$. El resultado que probaron es:
\begin{teo}\label{sub:teoconje}
La funci\'on $$\rho:[v_\alpha,\infty) \rightarrow \mathbb{R}: v \mapsto \frac{\tau(v)}{T(v)}$$ es anal\'itica real y es estrictamente creciente si $0<\alpha<1$ y estrictamente decreciente si $1<\alpha<\infty$.
\end{teo}

La demostraci\'on no se desarrolla en este texto por tener algunos aspectos te\'oricos de mayor dificultad, pero est\'a desarrollada y puede verse en \cite{cje}.

\subsection{\'Orbitas peri\'odicas.}
Ya conocemos que si $\alpha=0,1,\infty$  entonces $f_\alpha$ tiene \'orbitas de los per\'iodos 6,5,4, respectivamente. Ahora clasificaremos los per\'iodos de \'orbitas peri\'odicas para otros valores de $\alpha$.

\begin{teo}
Si $\alpha \neq 0,1,\infty$ entonces $f_\alpha$ tiene solo \'orbitas peri\'odicas de todos los per\'iodos  $q$ para los que exista un coprimo $p$ tal que $p/q$ est\'e entre $\rho_\alpha$ y $1/5$. El conjunto de todos estos per\'iodos generados de la familia $\{f_\alpha:\alpha \neq0,1,\infty \}$ es: $$9,11,13,14,16,17,19,\,\,y\,\,\forall n \in \mathbb{N}, n\geq 21,\,\, excepto\,\,42.$$ 
\end{teo}
\begin{proof}
En el Teorema~\ref{sub:valorro} $\rho^v_\alpha\rightarrow \rho_\alpha $ cuando $v\rightarrow v_\alpha$, en el Teorema~\ref{sub:aproximacionro} $\rho^v_\alpha\rightarrow 1/5 $ cuando $v\rightarrow \infty$, por la continuidad del n\'umero de rotaci\'on y por la monoton\'ia de $\rho^v_\alpha$, el conjunto $\{\rho^v_\alpha:v_\alpha<v<\infty\}$ es igual al intervalo abierto entre $\rho_\alpha$ y $1/5$. Si el n\'umero racional $p/q$ ($p,q$ coprimos) est\'a en el intervalo abierto entre $\rho_\alpha$ y $1/5$ entonces $p/q=\rho^v_\alpha$ para alg\'un $v$ en el dominio, luego, por el Corolario~\ref{sub:orbitasperiodicas} todas las \'orbitas de $f_\alpha$ sobre $C^v_\alpha$ tienen per\'iodo $q$, as\'i, se satisface la primera parte del enunciado.\\

Ahora, tenemos que encontrar el conjunto $Q$ de per\'iodos de \'orbitas peri\'odicas de $f_\alpha$, cuando $\alpha\neq 0,1,\infty$. Por lo expuesto arriba, y por el Teorema~\ref{sub:cotasro}, el conjunto $Q$  es el de todos los $q$ (exceto $5$, ya que se toma el intervalo abierto entre $\rho_\alpha$ y $1/5$) para los cuales  existe un coprimo $p$ talque $1/6<p/q<1/4$.\\

As\'i, $q\in Q$ si y s\'olo si existe un coprimo $p$ tal que $q/6<p<q/4$, para ello es suficiente encotrar al menos un primo $p$ en dicho rango, porque si $p$ es primo entonces $p/q$ ser\'ia una fraci\'on irreducible entre $1/6$ y $1/4$. Cuando $q$ es grande, veamos que el n\'umero de primos entre $q/6$ y $q/4$ es siempre mayor que $2$ por el teorema de n\'umeros primos\footnote{La aproximaci\'on de la cantidad de n\'umeros primos menores o iguales que $x$ es $Li(x)$, donde $Li(x)=\int_2^x \frac{dy}{ln{y}}$.}. La cantidad aproximada de n\'umero primos entre $q/6$ y $q/4$ es $$g(q)=Li(q/4)-Li(q/6)=\int_2^{q/4} \frac{dy}{ln{y}}-\int_2^{q/6} \frac{dy}{ln{y}}=\int_{q/6}^{q/4} \frac{dy}{ln{y}}.$$ Derivando est\'a funci\'on obtenemos que $\frac{dg}{dq}=\frac{1}{4\ln (\frac q4)}-\frac{1}{6\ln (\frac q6)}$, luego tenemos que:
\begin{align*} 
          \frac q6 & < \frac q4     \\ 
      \left(\frac q6\right)^6 & =\frac{q^4}{2^6}\frac{q^2}{3^6} \\
                   & >\frac{q^4}{2^8}=\left(\frac{q}{4}\right)^4, \tag*{$\frac{q^4}{2^6}>\frac{q^4}{2^8}$ y tomando a $q>60$, $\frac{q^2}{3^6}>1$} \\
   \ln \left(\frac q6\right)^6 & > \ln\left(\frac q4\right)^4 \\
   \frac{1}{\ln \left(\frac q6\right)^6} & <\frac{1}{\ln\left(\frac q4\right)^4} \\
    \frac{1}{6\ln \frac q6} & <  \frac{1}{4\ln \frac q4}
\end{align*}

As\'i, la derivada de $g(q)$ es positiva y se tiene que  $g(q)$ es creciente para todo $q>60$, adem\'as, $g(61)\approx 2.00914$ y  por tanto \'este es el $q$ grande adecuado para garantizar que al menos haya un primo en el intervalo (hemos buscamos el $q$ para el que $g$ sea  mayor o igual que dos por el error que tiene la f\'ormula del teorema de n\'umeros primos).\\
 
Si $q \geq 61$ entonces por lo anterior sabemos que siempre habr\'a al menos un primo en el intervalo y se tiene que $q \in Q,\forall q \geq 61$. Los n\'umeros menores que $61$ los cubrimos en la siguiente tabla, en la cual listamos para cada $q$  el conjunto de $p$ coprimos  tal que $p/q$ est\'e en el intervalo.
\begin{center}
$
\begin{tabular}{||c|c||c|c||c|c||c|c||c|c||c|c||}
\hline
\hline
q & p & q  & p & q & p & q & p & q & p & q & p \\
\hline
\hline
1 & - & 11 & 2 & 21 & 4,5 & 31 & 6,7 & 41 & 7,8,9,10 & 51 & 10,11\\
\hline
2 & - & 12 & - & 22 & 5 & 32 & 7 & 42 & - & 52 & 9,11\\
\hline
3 & - & 13 & 3 & 23 & 4,5 & 33 & 7,8 & 43 & 8,9,10 & 53 & 9,10,11,12,13\\
\hline
4 & - & 14 & 3 & 24 & 5 & 34 & 7 & 44 & 9 & 54 & 11, 13\\
\hline
5 & - & 15 & - & 25 & 6 & 35 & 6,8 & 45 & 8,11 & 55 & 12,13 \\
\hline
6 & - & 16 & 3 & 26 & 5 & 36 & 7 & 46 & 9,11 & 56 & 11,13\\
\hline
7 & - & 17 & 3,4 & 27 & 5 & 37 & 7,8,9 & 47 & 8,9,10,11 & 57 & 10,11,13,14 \\
\hline
8 & - & 18 & - & 28 & 5 & 38 & 7,9 & 48 & 11 & 58 & 11,13\\
\hline
9 & 2 & 19 & 4 & 29 & 5,6,7 & 39 & 7,8 & 49 & 9,10,11,12 & 59 & 10,11,12,13,14\\
\hline
10 & - & 20 & - & 30 & 7 & 40 & 7,9 & 50 & 9,11 & 60 & 11,13\\
\hline
\end{tabular}
$
\end{center}
Por lo tanto $Q=\{9,11,13,14,16,17,19,\,\,y\,\,\forall n \in \mathbb{N}, n\geq 21,\,\, excepto\,\,42\}$.
\end{proof}

{\bf \'Orbitas de per\'iodo $9$:} para encontrar \'orbitas de per\'iodo $9$ debemos ver las curvas $C$ con n\'umero de rotaci\'on $2/9$, siendo \'este el \'unico m\'ultiplo de $1/9$ que est\'a entre $1/6$ y $1/4$. Luego, por los Teoremas~\ref{sub:valorro},~\ref{sub:aproximacionro} y~\ref{sub:teoconje} debemos elegir a $\alpha$ tal que $2/9$ est\'e entre $1/5$ y $\rho_\alpha$. Como $2/9>1/5$ entonces $\rho_\alpha>2/9$. Por lo tanto $\alpha>\alpha_9$, donde $\alpha_9$ es obtenido por la ecuaci\'on 
\begin{equation}\label{sub:ecper9}\frac{2}{9}=\frac{1}{2\pi}\arccos \left(\frac{1}{1+\sqrt{1+4\alpha_9}}\right)\end{equation}
de donde $$\alpha_9=\frac{1-2\cos \left(\frac{4\pi}{9}\right)}{4\cos^2 \left(\frac{4\pi}{9}\right)}\approx 5.4114741...$$
Ahora ya sabiendo cuales son los valores de $\alpha$ adecuados, damos el siguiente resultado para encontrar los valores correspondientes de $v$ y as\'i  obtener \'orbitas peri\'odicas de per\'iodo 9. 

\begin{teo}
\label{sub:orbitasperiodo9} Todas las \'orbitas sobre $C^{\bar{v}}_\alpha$ son peri\'odicas con per\'iodo 9 si y s\'olo si  $\alpha_9<\alpha<\infty$ y $\bar{v}=\frac{(\alpha-1)(\alpha^2-\alpha+1)}{\alpha}$.
\end{teo}
\begin{proof}
Sea $f=f_\alpha$ y $C=C^{\bar{v}}_\alpha$ donde $0<\alpha<\infty$ y $v_\alpha<\bar{v}<\infty$. Mostraremos que toda \'orbita sobre $C$ es peri\'odica de per\'iodo 9 si y s\'olo si $\alpha_9<\alpha<\infty$ y $\bar{v}=\frac{(\alpha-1)(\alpha^2-\alpha+1)}{\alpha}$.\\

Por el Corolario~\ref{sub:orbitasperiodicas} toda \'orbita sobre  $C$ tiene el mismo per\'iodo, as\'i, es suficiente analizar una \'unica \'orbita sobre $C$. Consideremos la  \'orbita sim\'etrica a trav\'es de un punto de la diagonal $P=(\lambda,\lambda) \in C$. Iniciemos en $P$ construyendo los ocho siguientes puntos sim\'etricamente sobre $C$ dibujando alternadamente cuerdas verticales y horizontales a partir de $P$, como se muestra en la figura~\ref{sub:orbita9}.

\begin{figure}[h!]
\begin{center}
\includegraphics[width=8cm,height=8cm]{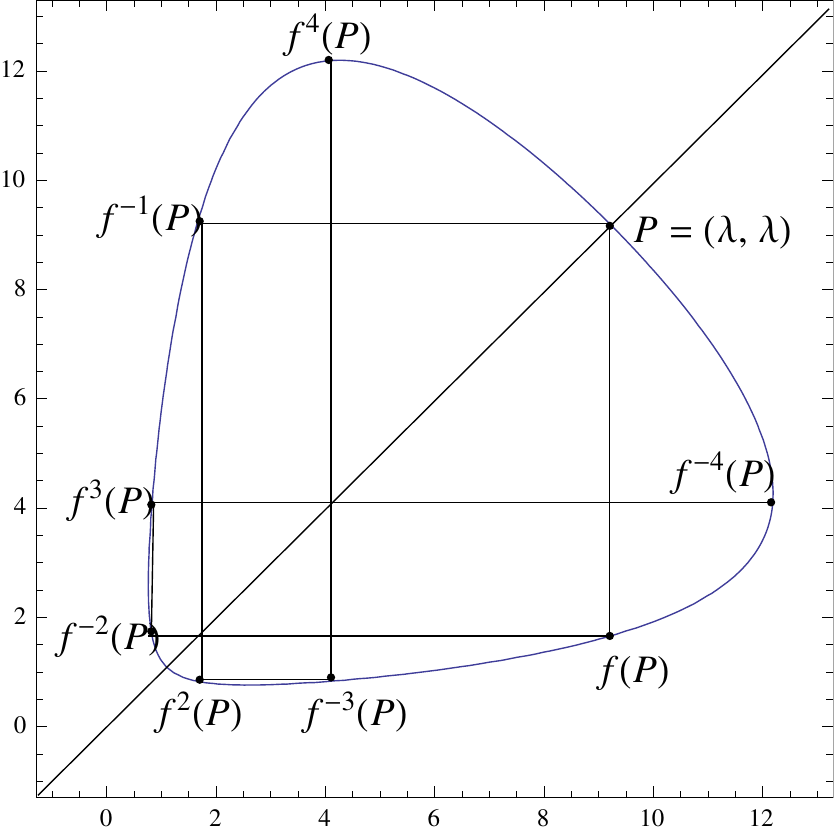}
\end{center}
\caption[\'Orbita de per\'iodo 9.]{\'Orbita sim\'etrica de per\'iodo 9. \label{sub:orbita9}}
\end{figure}

Hay que subrayar que $f|_C$ es la composici\'on de dos involucones, la primera obtenida intercambiando los puntos extremos de cada cuerda horizontal, y la segunda por la simetr\'ia en la diagonal (como se expuso en el ejemplo~\ref{sub:invo2}). Usando esta representaci\'on de $f|_C$ podemos identificar los nueve puntos como los nueve primeros  puntos $$f^{-4}(P),f^{-3}(P),f^{-2}(P),f^{-1}(P),P,f(P),f^2(P),f^3(P),f^4(P)$$ de la  \'orbita de $f^{-4}(P)$, como se muestra en  la figura~\ref{sub:orbita9}. Esta \'orbita tiene per\'iodo 9 si y s\'olo si $f$ transforma $f^4(P)$ en $f^{-4}(P)$. Por simetr\'ia justamente $f^4(P)$ se transforma en $f^{-4}(P)$. Por lo tanto, para que lo anterior suceda se tiene que la cuerda horizontal que pasa por $f^4(P)$ debe transformar a $f^4(P)$ en s\'i misma; en otras palabras esta cuerda debe ser tangente  horizontalmente a la curva $C$, es decir, $\frac{\partial V}{\partial x}=0$ en $f^4(P)$, luego las coordenadas de $f^4(P)$ satisfacen \begin{equation}\label{sub:puntof4} x^2-y-\alpha=0.\end{equation} As\'i, la sucesi\'on determinada por $P$ es
\begin{align*} 
    x_1 & = \lambda     \\ 
    x_2 & =\lambda \\
    x_3 & =\frac{\lambda+\alpha}{\lambda} \\
    x_4 & = \frac{(\alpha+1)\lambda+\alpha}{\lambda^2} \\
    x_5 & = \frac{\alpha\lambda^2+(\alpha+1)\lambda+\alpha}{\lambda(\lambda+\alpha)}\\
    x_6 & = \frac{\lambda(2\alpha\lambda^2+(\alpha^2+\alpha+1)\lambda+\alpha)}{(\lambda+\alpha)((\alpha+1)\lambda+\alpha)}
\end{align*}
Sustituyendo $f^4(P)=(x_5,x_6)$ en la ecuaci\'on (~\ref{sub:puntof4}) obtenemos: $$\left(\frac{\alpha\lambda^2+(\alpha+1)\lambda+\alpha}{\lambda(\lambda+\alpha)}\right)^2-\frac{\lambda(2\alpha\lambda^2+(\alpha^2+\alpha+1)\lambda+\alpha)}{(\lambda+\alpha)((\alpha+1)\lambda+\alpha)}-\alpha=0.$$ Multiplicando por denominador y fatorizando tenemos: $$(\lambda+1)(\lambda^2-\lambda-\alpha)\left[2\alpha \lambda^3+(-\alpha^3+3\alpha^2+2\alpha +1)\lambda^2+2\alpha(\alpha+1)\lambda+\alpha^2\right]=0.$$
Ignoramos el primer factor\footnote{La raiz $\lambda=-1$ es matem\'aticamente trascendente porque implica que $\bar{v}=0$, y as\'i la curva c\'ubica se descompone en tres l\'ineas, en las que las intersecciones forman una \'orbita peri\'odica de per\'iodo tres, que, aunque esta fuera de $\mathbb{R}^2_+$, aparecer\'ia entre las \'orbitas de per\'iodo nueve.} porque  $\lambda>0$ y el segundo tambi\'en lo descartamos porque es cero solo en el punto fijo $F$ cuando $\bar{v}=v_\alpha$, y es diferente de cero cuando $\bar{v}>v_\alpha$ que es el caso que nos interesa. Por lo tanto el tercer factor debe ser cero\begin{equation}\label{sub:ecfactor3}2\alpha\lambda^3+2\alpha(\alpha+1)\lambda+\alpha^2=(\alpha^3-3\alpha^2-2\alpha -1)\lambda^2.\end{equation}
As\'i,
\begin{align*} 
\bar{v} & = V(P)    \\ 
        & =\frac{(\lambda+1)^2(2\lambda+\alpha)}{\lambda^2} \\
        & =\frac{2\alpha\lambda^3+\alpha(\alpha+4)\lambda^2+2\alpha(\alpha+1)\lambda+\alpha^2}{\alpha\lambda^2} \\
        & = \frac{\alpha(\alpha+4)\lambda^2+(\alpha^3-3^2-2\alpha-1)\lambda^2}{\alpha \lambda^2}, \tag*{Por (~\ref{sub:ecfactor3}),} \\
        & = \frac{\alpha^3-2\alpha^2+2\alpha-1}{\alpha}\\
        & = \frac{(\alpha-1)(\alpha^2-\alpha+1)}{\alpha}.
\end{align*}

Hemos probado que $C$ contiene una \'orbita de per\'iodo $9$ si y s\'olo si esta \'ultima f\'ormula para $\bar{v}$ es cierta con $\bar{v}>v_\alpha$ (porque $C$ existe s\'olo si $\bar{v}>v_\alpha$), por el siguiente Lema~\ref{sub:lemadevalfa} $\bar{v}>v_\alpha$ si y s\'olo si $\alpha>\alpha_9$ ($\alpha\neq\infty$, ya que como vimos en este caso todas las \'orbitas tienen per\'iodo $4$). As\'i, $C$ existe y tiene per\'iodo $9$ si y s\'olo si $\bar{v}=\frac{(\alpha-1)(\alpha^2-\alpha+1)}{\alpha}$ con $\bar{v}>v_\alpha$  y \'esto \'ultimo es cierto si y s\'olo si $\alpha>\alpha_9$.
\end{proof}

\begin{lema}
\label{sub:lemadevalfa}
$\bar{v}\geq v_\alpha$ si y s\'olo si $\alpha\geq\alpha_9$ y $\bar{v}\leq v_\alpha$ si y s\'olo si $\alpha\leq\alpha_9$.
\end{lema}
\begin{proof}
Recordemos que $\alpha_9=\frac{1-2c}{4c^2}$ donde $c=\cos \left(\frac{4\pi}{9} \right)$, luego $4c^3-3c=\cos \left(3\frac{4\pi}{9}\right)=-\cos\left(\frac{\pi}{3}\right)=-\frac{1}{2}$, as\'i, $8c^3-6c+1=0$.\\

Primero mostremos que $\alpha=\alpha_9$ implica que $\bar{v}=v_\alpha$, como sigue.
\begin{align*} 
\bar{v}& =  \frac{(\alpha_9-1)(\alpha^2_9-\alpha_9+1)}{\alpha_9},\tag*{sustituyendo $\alpha=\alpha_9$,}  \\ 
       & =\frac{(1-2c-4c^2)(1-4c+8c^3+16c^4)}{16c^4(1-2c)},\tag*{con $\alpha_9=\frac{1-2c}{4c^2}$,}\\
       & =\frac{1-6c+4c^2+24c^3-64c^5-64c^6}{16c^4(1-2c)} \\
       & =\frac{4c^2+16c^3-64c^5-64c^6}{16c^4(1-2c)}, \tag*{por que $1-2c=-8c^3$,}\\
       & =\frac{(1+2c)^3}{4c^2}, \tag*{por factorizaci\'on y cancelaci\'on de t\'erminos,}\\
       & =\frac{(1+\frac 1\omega_9)^3}{(\frac 1\omega_9)^2}, \tag*{por ecuaci\'on (~\ref{sub:ecper9}) $\omega_9=\frac{1+\sqrt{1+4\alpha_9}}{2}=\frac{1}{2c}$,}\\           & =\frac{(1+\omega_9)^3}{\omega_9}\\
       & = v_\alpha.
\end{align*}

Mostraremos ahora que la gr\'afica de $$\bar{v}= \frac{(\alpha-1)(\alpha^2-\alpha+1)}{\alpha} \qquad   y \qquad v_\alpha=\frac{(1+\omega)^3}{\omega}=\frac{(3+\sqrt{1+4\alpha})^3}{4(1+\sqrt{1+4\alpha})}$$ Se cortan solo una vez, ver (figura ~\ref{sub:graficovv}).\\

\begin{figure}[h!]
\begin{center}
\includegraphics[width=8cm,height=8cm]{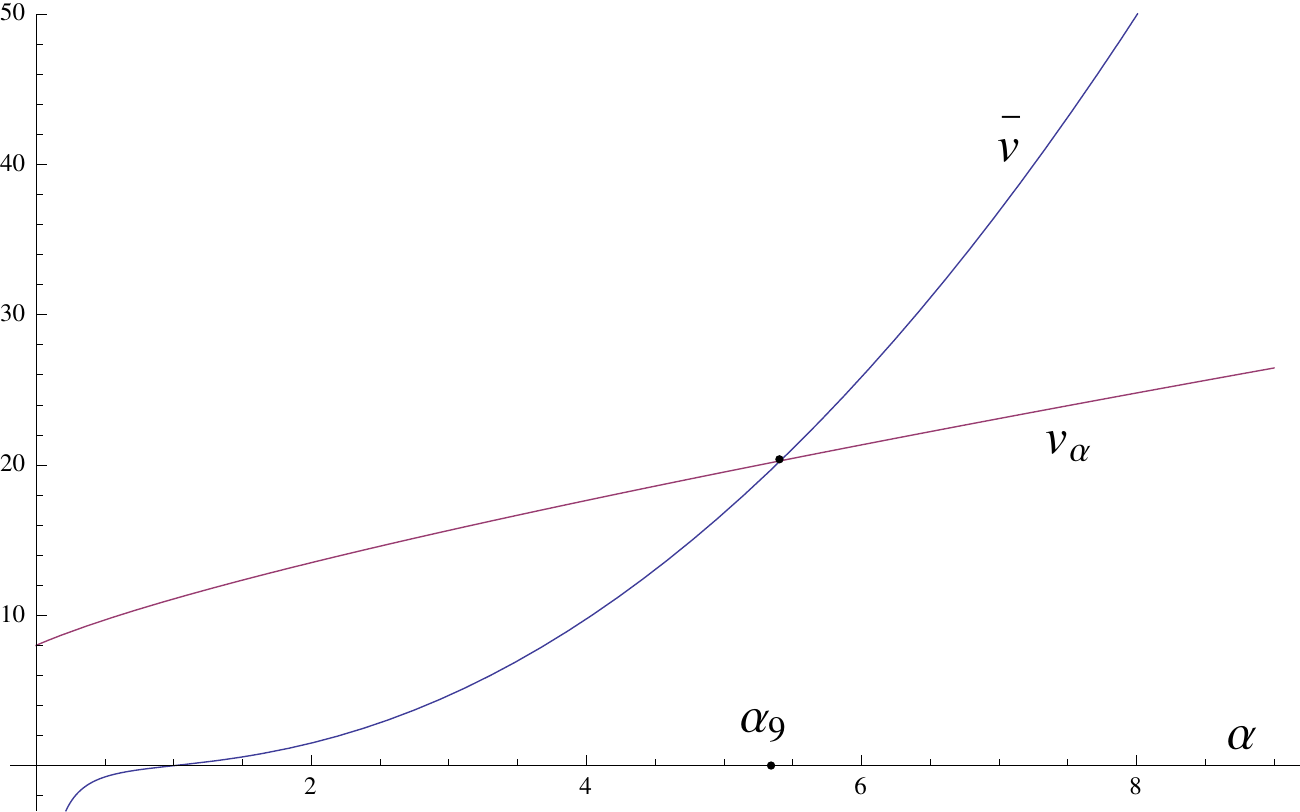}
\end{center}
\caption[Gr\'afica de $\bar{v}$ y $v_\alpha$ para $\alpha_9$]{La gr\'afica de $\bar{v}$ y $v_\alpha$ se cortan solo una vez y lo hacen en $\alpha=\alpha_9$. \label{sub:graficovv}}
\end{figure}

Diferenciando: $$\frac{d\bar{v}}{d\alpha}= \frac{2\alpha^3-2\alpha^2+1}{\alpha} \qquad   y \qquad \frac{dv_\alpha}{d\alpha}=\left(\frac{3+\sqrt{1+4\alpha}}{1+\sqrt{1+4\alpha}}\right)^2.$$ 
Por lo tanto, si $0<\alpha<\infty$  entonces $\frac{d\bar{v}}{d\alpha}>0$ (ya que $-\frac12< 2\alpha^3-2\alpha^2,\forall \alpha$ en su dominio) y $\frac{dv_\alpha}{d\alpha}>0$ entonces ambas son crecientes. Ahora, $v(3)=\frac{14}{3}<8=v_0$, as\'i, si $0<\alpha\leq3$ entonces $v<8<v_\alpha$ por lo tanto las gr\'aficas no se cortan cuando $0<\alpha\leq3$. Si $\alpha=3$ entonces $\frac{d\bar{v}}{d\alpha}>4>\frac{dv_\alpha}{d\alpha}$, m\'as a\'un, \'esto es cierto para todo $\alpha$ en el dominio $3\leq\alpha<\infty$ porque en este dominio $\frac{d\bar{v}}{d\alpha}$ es creciente y $\frac{dv_\alpha}{d\alpha}$ es decreciete, ya que $\frac{d^2\bar{v}}{d\alpha^2}=\frac{2(\alpha-1)(\alpha^2+\alpha+1)}{\alpha^3}$ es positiva y $\frac{d^2v_\alpha}{d\alpha^2}=-\frac{8(3+\sqrt{1+4\alpha})}{\sqrt{1+4\alpha}(1+\sqrt{1+4\alpha})^3}$ es negativa. Con esta informaci\'on sabemos que las gr\'aficas pueden cortarse solo una vez en su dominio, pongamos en $\alpha=\alpha_9$. Ya verificamos que: si $0<\alpha\leq3$ entonces $\bar{v}<8<v_\alpha$, sus gr\'aficas  se cortan solo una vez y si $3\leq\alpha<\infty$ entonces $\frac{d\bar{v}}{d\alpha}>4>\frac{dv_\alpha}{d\alpha}$, por lo tanto todas las implicaciones del enunciado son ciertas. 

\end{proof} 

{\bf \'Orbitas de per\'iodo 11:} de forma semejante al caso de \'orbitas de per\'iodo 9, para encontrar \'orbitas de per\'iodo 11 debemos elegir $\alpha$ tal que $2/11=0.\bar{18}$ est\'e entre $\rho_\alpha$ y $1/5$. Por lo tanto, $\alpha<\alpha_{11}$ donde $$\alpha_{11}=\frac{1-2\cos \left(\frac{4\pi}{11}\right)}{4\cos^2 \left(\frac{4\pi}{11}\right)}\approx 0.2450749...$$
\begin{teo}
\label{sub:orbitasperiodo11} Todas las \'orbitas sobre $C^{\bar{v}}_\alpha$ son peri\'odicas con per\'iodo 11 si y s\'olo si  $0<\alpha<\alpha_{11}$ y $$\bar{v}=\frac{(1-\alpha)(1-2\alpha+\sqrt{1-4\alpha^2+4\alpha^3})}{2\alpha^2}.$$
\end{teo}
\begin{proof}
La prueba es repetir los pasos de la demostraci\'on del Teorema~\ref{sub:orbitasperiodo9} con algunas peque\~nas variantes. Sea $f=f_\alpha$ y $C=C^{\bar{v}}_\alpha$ donde $0<\alpha<\alpha_{11}$  y $v_\alpha<\bar{v}<\infty$. Mostraremos que toda \'orbita sobre $C$ es peri\'odica de per\'iodo 11 si y s\'olo si $0<\alpha<\alpha_{11}$ y $\bar{v}=\frac{(1-\alpha)(1-2\alpha+\sqrt{1-4\alpha^2+4\alpha^3})}{2\alpha^2}.$\\

Ya sabemos que es suficiente ver una \'unica \'orbita sobre $C$ porque todas tienen el mismo per\'iodo. Tomemos el  punto de la diagonal $P=(\lambda,\lambda) \in C$. Iniciemos en $P$ construyendo los diez siguientes puntos sim\'etricamente sobre $C$ dibujando alternadamente cuerdas verticales y horizontales a partir de $P$, como se muestra en  la figura~\ref{sub:orbita11}.

\begin{figure}[h!]
\begin{center}
\includegraphics[width=8cm,height=8cm]{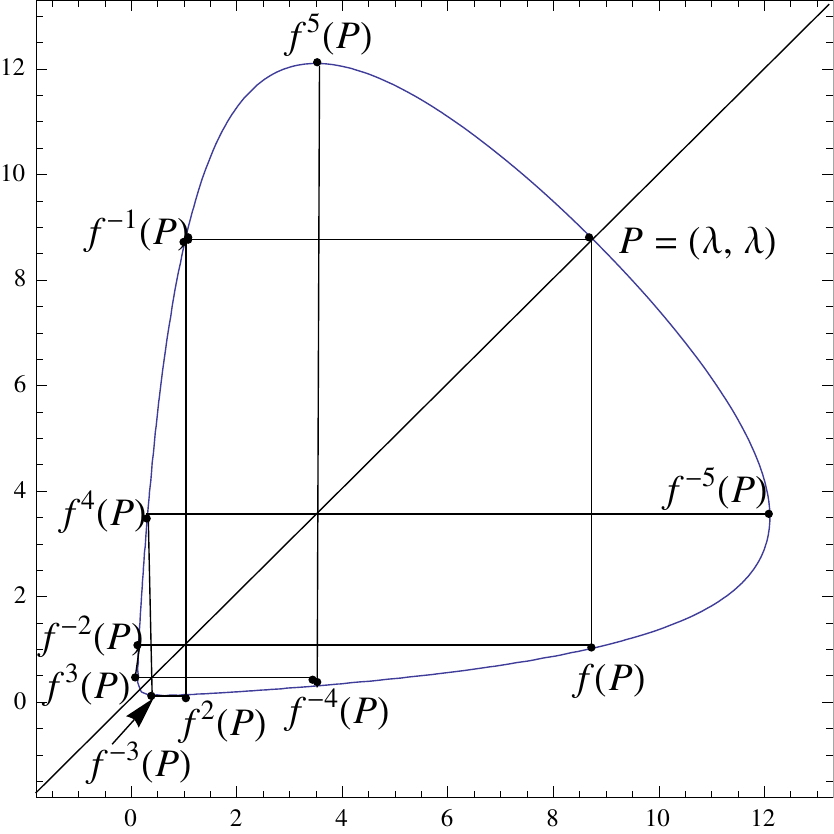}
\end{center}
\caption[\'Orbita de per\'iodo 11.]{\'Orbita sim\'etrica de per\'iodo 11. \label{sub:orbita11}}
\end{figure}

Utilizando la representaci\'on de las dos involuciones para  $f|_C$ podemos identificar los once puntos como los once primeros  puntos $$f^{-5}(P),f^{-4}(P),f^{-3}(P),f^{-2}(P),f^{-1}(P),P,f(P),f^2(P),f^3(P),f^4(P),f^5(P)$$ de la  \'orbita de $f^{-5}(P)$ (ver figura ~\ref{sub:orbita11}). Esta \'orbita tiene per\'iodo 11 si y s\'olo si $f$ transforma a $f^5(P)$ en $f^{-5}(P)$. Por simetr\'ia  $f^5(P)$ se transforma en $f^{-5}(P)$. La cuerda en $f^5(P)$ debe ser tangente  horizontalmente a la curva $C$, es decir, $\frac{\partial V}{\partial x}=0$ en $f^5(P)$, luego las coordenadas de $f^5(P)$ satisfacen \begin{equation}\label{sub:puntof5} x^2-y-\alpha=0.\end{equation} La sucesi\'on determinada por $P$ es
\begin{align*} 
    x_1 & = \lambda     \\ 
    x_2 & =\lambda \\
    x_3 & =\frac{\lambda+\alpha}{\lambda} \\
    x_4 & =\frac{(\alpha+1)\lambda+\alpha}{\lambda^2} \\
    x_5 & =\frac{\alpha\lambda^2+(\alpha+1)\lambda+\alpha}{\lambda(\lambda+\alpha)}\\
    x_6 & =\frac{\lambda(2\alpha\lambda^2+(\alpha^2+\alpha+1)\lambda+\alpha)}{(\lambda+\alpha)((\alpha+1)\lambda+\alpha)}\\
    x_7 & =\frac{\lambda\left(\lambda+\alpha^2\lambda+\alpha(1+\lambda+2\lambda^2)\right)}{(\alpha+\lambda)(\alpha+\lambda+\alpha\lambda)}\\
\end{align*}
Sustituyendo $f^5(P)=(x_6,x_7)$ en la ecuaci\'on (~\ref{sub:puntof5}) obtenemos: \begin{equation*}  \begin{split}\left(\frac{\lambda(2\alpha\lambda^2+(\alpha^2+\alpha+1)\lambda+\alpha)}{(\lambda+\alpha)((\alpha+1)\lambda+\alpha)}\right)^2 -\frac{\lambda\left(\lambda+\alpha^2\lambda+\alpha(1+\lambda+2\lambda^2)\right)}{(\alpha+\lambda)(\alpha+\lambda+\alpha\lambda)}-\alpha=0. \end{split}\end{equation*}

Reagrupando y simplificando:
\begin{equation*} \begin{split} & \left(\lambda^2-\lambda-\alpha\right)\left[2\lambda^3+(A+B)\lambda^2+2(\alpha+1)\lambda+\alpha\right]\\ & \quad \left[2\lambda^3+(A-B)\lambda^2+(\alpha+1)\lambda+\alpha\right]=0.\end{split}\end{equation*}

Donde $$A=\frac{2\alpha^3+6\alpha^2+3\alpha-1}{2\alpha^2}\qquad   y \qquad B=\frac{(1-\alpha)\sqrt{1-4\alpha^2+4\alpha^3}}{2\alpha^2}.$$

Ignoramos el primer  factor porque  es cero solo en el punto fijo $F$. \\

Por lo tanto \begin{equation}\label{sub:ecfactor311}2\lambda^3+2(\alpha+1)\lambda+\alpha=(-A\pm B)\lambda^2.\end{equation}
Luego
\begin{align*} 
\bar{v}&= V(P)  \\
       &=\frac{2\alpha\lambda^3+\alpha(\alpha+4)\lambda^2+2\alpha(\alpha+1)\lambda+\alpha^2}{\alpha\lambda^2},\tag*{como en el Teorema~\ref{sub:orbitasperiodo9},}\\
       &= (\alpha+4)+(-A\pm B), \tag*{por la ecuaci\'on (~\ref{sub:ecfactor311}),} \\
       &= \frac{\left(1-\alpha\right)\left[1-2\alpha\pm \sqrt{1-4\alpha^2+4\alpha^3}\right]}{2\alpha^2}.
\end{align*}

Si tomamos la ra\'iz negativa tenemos que $\bar{v}<v_\alpha,\forall \alpha$, que es una contradicci\'on. Por lo tanto tomamos la ra\'iz positiva. Como sabemos, $C$ existe s\'olo si $\bar{v}>v_\alpha$, y por el siguiente Lema~\ref{sub:lemadevalfa11} $\bar{v}>v_\alpha$ si y s\'olo si $\alpha<\alpha_{11}$ ($\alpha>0$, ya que si $\alpha=0$ todas sus \'orbitas tienen per\'iodo 6). Por lo tanto, $C$ contiene solo \'orbitas de per\'iodo 11 si y s\'olo si  $\bar{v}=\frac{\left(1-\alpha\right)\left[1-2\alpha + \sqrt{1-4\alpha^2+4\alpha^3}\right]}{2\alpha^2}$ y $\bar{v}>v_\alpha$, y \'esto \'ultimo es cierto si y s\'olo si $\alpha<\alpha_{11}$.
\end{proof}

\begin{lema}
\label{sub:lemadevalfa11}
$\bar{v}\geq v_\alpha$ si y s\'olo si $\alpha\leq\alpha_{11}$ y $\bar{v}\leq v_\alpha$ si y s\'olo si $\alpha\geq\alpha_{11}$.
\end{lema}
\begin{proof}
Utilizando el \emph{software Mathematica 7.0} podemos verificar que $\bar{v}(\alpha_{11})=v_\alpha(\alpha_{11})\approx 8.890367467127085$, as\'i $\bar{v}$ y $v_\alpha$ se cortan en $\alpha_{11}$ (ver figura~\ref{sub:funcionesvv11}).

\begin{figure}[h!]
\begin{center}
\includegraphics[width=8cm,height=8cm]{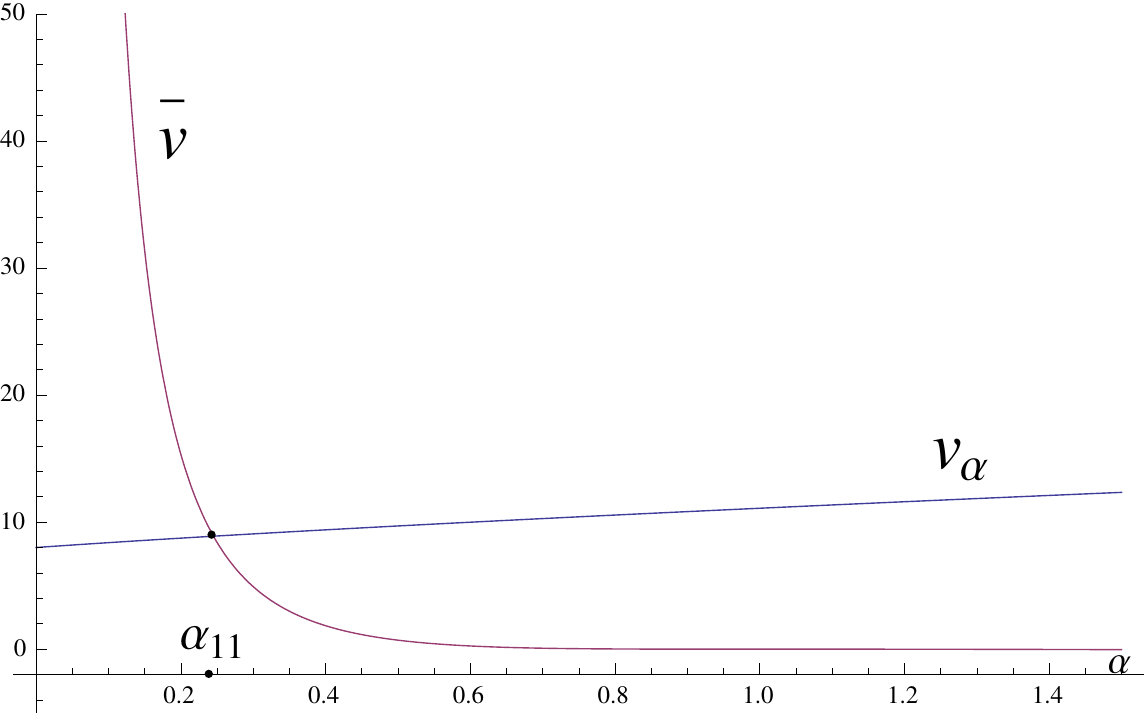}
\end{center}
\caption[Gr\'afica de $\bar{v}$ y $v_\alpha$ para $\alpha_{11}$]{La gr\'afica de $\bar{v}$ y $v_\alpha$ se cortan solo una vez y lo hacen en $\alpha=\alpha_{11}$. \label{sub:funcionesvv11}}
\end{figure}

Diferenciando:\\

$$\frac{d\bar{v}}{d\alpha}=
\frac{-2+\alpha+4\alpha^2-3\alpha^3-2\alpha^4-2\sqrt{1-4\alpha^2+4\alpha^3}+3\alpha\sqrt{1-4\alpha^2+4\alpha^3}}{2\alpha^3\sqrt{1-4\alpha^2+4\alpha^3}}$$
$$\frac{dv_\alpha}{d\alpha}=\left(\frac{3+\sqrt{1+4\alpha}}{1+\sqrt{1+4\alpha}}\right)^2$$

La \'unica ra\'iz real de $\frac{d\bar{v}}{d\alpha}=0$ es $\alpha=1$, $\frac{d\bar{v}}{d\alpha}(0.5)<0$ y $\frac{d\bar{v}}{d\alpha}(1.5)<0$ entonces $\frac{d\bar{v}}{d\alpha}<0, \forall \alpha>0 (\alpha\neq1)$ por continuidad de $\bar{v}$ para $\alpha>0$. $\frac{dv_\alpha}{d\alpha}>0,\forall \alpha>0$. Como las derivadas de $\bar{v}$ y $v_\alpha$ son de signo opuesto  entonces estas dos funciones solo pueden cortarse una vez, llamese  $\alpha=\alpha_{11}$. Adem\'as, $\bar{v}(1)=0<v_\alpha(1)$, por lo tanto todas las implicaciones del enunciado son ciertas.
\end{proof}

\bibliographystyle{plain}

\end{document}